\documentclass[11pt]{article}
\usepackage{mathrsfs}
\usepackage{amssymb}
\usepackage{amsmath}
\usepackage{amsbsy}
\usepackage{epsfig}
\usepackage{enumerate}
\usepackage{bm}
\usepackage{color}

\usepackage[colorlinks,linkcolor=blue]{hyperref}

\newtheorem{lemma}{Lemma}[section]
\newtheorem{theorem}{Theorem}[section]
\newtheorem{proposition}{Proposition}[section]

\newtheorem{corollary}{Corollary}[section]
\newtheorem{definition}{Definition}[section]
\newtheorem{remark}{Remark}[section]
\newtheorem{conjecture}{Conjecture}[section]

\def\pr{\textsf{P}} 
\def\ep{\textsf{E}} 
\def\Sbep{\widehat{\mathbb E}} 
\def\cSbep{\widehat{\mathcal E}} 
\def\Capc{\mathbb V} 
\def\cCapc{\mathcal V} 
\def\upCapc{\widehat{\mathbb V}} 
\def\lowCapc{\widehat{\mathcal V}} 
\def\outCapc{\widehat{\mathbb V}^{\ast}}
\def\outcCapc{\widehat{\mathcal V}^{\ast}} 

\def\tildeCapc{\mathbb C^{\ast}}

\def\vSbep{\breve{\mathbb E}}

\renewcommand{\baselinestretch}{1.5}

\begin{document}

\begin{center}{\LARGE\bf On the laws of the iterated logarithm   under sub-linear expectations without the assumption on the continuity of capacities}
\end{center}

\begin{center} {\sc
Li-Xin Zhang\footnote{This work was supported by grants from the NSF of China (Nos. 11731012,12031005),   Ten Thousands Talents Plan of Zhejiang Province (Grant No. 2018R52042), NSF of Zhejiang Province (No. LZ21A010002) and the Fundamental Research Funds for the Central Universities
}
}\\
{\sl \small School  of Mathematical Sciences, Zhejiang University, Hangzhou 310027} \\
(Email:stazlx@zju.edu.cn)\\
\end{center}

 \renewcommand{\abstractname}{~}
\begin{abstract}
{\bf Abstract:}
In this paper, we establish some general forms of the law of the iterated logarithm for independent random variables   in a sub-linear expectation space, where  the random variables are not necessarily identically distributed. Exponential  inequalities for the maximum sum  of independent random variables and Kolmogorov's converse exponential  inequalities are established as tools for showing the law of the iterated logarithm.  As an application, the sufficient and necessary conditions of the law of the iterated logarithm for independent and identically distributed random variables under the sub-linear expectation are obtained. In the paper, it is also shown that if the sub-linear expectation space is rich and regular enough, it will have no continuous capacity. The laws of the iterated logarithm are established without the assumption on the continuity of  capacities.


{\bf Keywords:}  sub-linear expectation, capacity, Kolmogorov's exponential inequality,
 laws of the iterated logarithm

 {\bf AMS 2020 subject classifications:}  60F15, 60F05

\vspace{-3mm}
\end{abstract}

\renewcommand{\baselinestretch}{1.2}


\section{ Introduction and notations.}\label{sect1}
\setcounter{equation}{0}

Let $\{X_n;n\ge 1\}$ be a sequence of independent random variables on a probability space $(\Omega,\mathcal F, \pr)$. Denote $S_n=\sum_{i=1}^n X_i$, $s_n^2=\sum_{i=1}^n \ep X_i^2$, $\log x=\ln \max(e,x)$, where $\ep$ is the expectation with respect to $\pr$.  The almost sure limit behavior of $\{S_n/\sqrt{2s_n^2\log\log s_n^2};n\ge 1\}$ has been studied extensively. It is known, under some conditions, that
\begin{equation}\label{WittmannLIL0}
\pr\left(\limsup_{n\to \infty}\frac{|S_n|}{\sqrt{2s_n^2\log\log s_n^2}}=1\right)=1.
\end{equation}
 This is the "Law of the  Iterated Logarithm" (LIL). In his   well-known  paper, Wittmann (1985) established a general theorem for LIL which  states that \eqref{WittmannLIL0} holds if the following conditions are fulfilled:
 \begin{equation}\label{WittmannLIL1}
 \ep X_n=0 \; \text{ and } \; \ep X_n^2<\infty, \;\; n\ge 1,
 \end{equation}
 \begin{equation}\label{WittmannLIL2}
 \sum_{n=1}^{\infty} \frac{\ep |X_n|^p}{(2s_n^2\log\log s_n^2)^{p/2}}< \infty \;\; \text{ for some } 2<p\le 3,
 \end{equation}
 [ or $\ep X_n^3=0$, $n\ge 1$, and \eqref{WittmannLIL2} holds for some $3<p\le 4$],
 \begin{equation}\label{WittmannLIL3}
 \lim_{n\to \infty} s_n^2=\infty \;\; \text{ and } \; \limsup_{n\to \infty} \frac{s_{n+1}^2}{s_n^2}<\infty.
 \end{equation}
 According to Wittmann, the classical result of Hartman and Wintner (1941) is just a corollary of his theorem. That is, if $\{X_n;n\ge 1\}$ is a sequence of independent and identically distributed (i.i.d.) random variables, then
 \begin{equation}\label{HWLIL1}
\pr\left(\limsup_{n\to \infty}\frac{|S_n|}{\sqrt{2n\log\log n}} =(\ep X_1^2)^{1/2}\right)=1
\end{equation}
 if
 \begin{equation}\label{HWLIL2}
 \ep X_1=0 \;\;\text{ and }\;\; \ep X_1^2<\infty.
\end{equation}
Wittmann (1987) showed that his theorem also holds when $p>3$. Chen (1993) extended Wittmann's theorem to the case of random variables taking their values  in a Banach space and weakened  the condition \eqref{WittmannLIL2} to    that for every $\epsilon>0$ there exists $p>2$ such that
 \begin{equation}\label{ChenX1}
 \sum_{n=1}^{\infty} \frac{\ep |X_n|^pI\{|X_n|\ge \epsilon \sqrt{s_n^2/\log\log s_n^2}\}}{(2s_n^2\log\log s_n^2)^{p/2}}< \infty.
 \end{equation}

In this paper, we consider the random variables in a   sub-linear expectation space. Let $\{X_n;n\ge 1\}$ be sequence of independent random variables in a sub-linear expectation space $(\Omega, \mathscr{H},\Sbep)$   with a related upper capacity $\Capc$. Chen and Hu (2014) showed that, if $\{X_n;n\ge 1\}$  is a sequence of i.i.d. random variables,
 then
 \begin{equation}\label{CH0}
\Capc\left(\limsup_{n\to \infty}\frac{|S_n|}{\sqrt{2n\log\log n}} >(\Sbep X_1^2)^{1/2}\right)=0
\end{equation}
if
\begin{equation}\label{CH1}
 \Sbep[X_1]=\Sbep[-X_1]=0,
\end{equation}
\begin{equation}\label{CH2}
 X_1,X_2,\ldots, \text{ are bounded random variables}.
\end{equation}
Zhang (2016) showed that \eqref{CH0} holds if \eqref{CH1}, and
\begin{equation}\label{Zhang1}
\lim_{c\to\infty} \Sbep[(X_1^2-c)^+]=0,
\end{equation}
\begin{equation}\label{Zhang2}
\int_0^{\infty}\Capc\left(|X_1|^2/\log\log |X_1|\ge x\right)dx<\infty.
\end{equation}
It is obvious that \eqref{Zhang1} and \eqref{Zhang2} are much weaker than \eqref{CH2}, and together with \eqref{CH1} are very close to Hartman and Wintner's condition \eqref{HWLIL2}. Zhang also showed that \eqref{Zhang2} is necessary for \eqref{CH0} to hold. Nevertheless two important questions remained unanswered:
\begin{enumerate}
  \item Is \eqref{Zhang1} also a necessary  condition? It is known that Hartman and Wintner's condition \eqref{HWLIL2} is also  necessary   for \eqref{HWLIL1} to hold (cf.  Strassen (1966)).   What are the sufficient and necessary conditions for  \eqref{CH0} to hold? It should be noted that,  unlike in the classical probability space, $\Sbep[X_1^2]<\infty$ and \eqref{Zhang1} are not equivalent and neither of them  implies \eqref{Zhang2}.
  \item Does Wittmann's theorem also hold  under the sub-linear expectation?
\end{enumerate}
A big difficulty for showing  the necessity  of a kind of the condition \eqref{Zhang1} for \eqref{CH0} is that the symmetrization argument is not valid under the sub-linear expectation.     As for Wittmann's LIL, beside we have not enough powerful exponential inequalities, a difficulty is that we can not use the truncation argument under the sub-linear expectation  as freely as under the classical expectation because, if a random variable $X$ is partitioned to $X_1+X_2$, the sub-linear expectation $\Sbep[X]$ is no longer $\Sbep[X_1]+\Sbep[X_2]$. The purpose of this paper is to establish   LIL  for independent random   variables under the sub-linear expectation, where the random variables are not necessarily identically-distributed. As a corollary, we obtain the sufficient and necessary conditions of the LIL for i.i.d. random  variables.

This paper is organized as follows. In section \ref{sectBasic}, some notation  under the sub-linear expectation is introduced.  The main tools are established in section~\ref{sectIneq}, including  exponential inequalities for the maximum sum of independent random variables and Kolmogorov's converse exponential inequality.  Some properties of the capacities are given in  section \ref{sectCapc} where, as a corollary, it is showed that the $G$-capacity is not continuous and the converse part of the Borel-Cantelli lemma is usually not valid for a capacity.
  In section~\ref{sectLIL}, the theorems on   the laws of the iterated logarithm are given.    The proofs of the laws of the iterated logarithm are shown in section \ref{sectProof}.

\section{Basic settings}\label{sectBasic}
\setcounter{equation}{0}

We use the framework and notations of Peng (2008, 2009, 2019). Let  $(\Omega,\mathcal F)$
 be a given measurable space  and let $\mathscr{H}$ be a linear space of real measurable functions
defined on $(\Omega,\mathcal F)$ such that if $X_1,\ldots, X_n \in \mathscr{H}$  then $\varphi(X_1,\ldots,X_n)\in \mathscr{H}$ for each
$\varphi\in C_{l,Lip}(\mathbb R^n)$,  where $C_{l,Lip}(\mathbb R^n)$ denotes the linear space of (local Lipschitz)
functions $\varphi$ satisfying
\begin{eqnarray*} & |\varphi(\bm x) - \varphi(\bm y)| \le  C(1 + |\bm x|^m + |\bm y|^m)|\bm x- \bm y|, \;\; \forall \bm x, \bm y \in \mathbb R^n,&\\
& \text {for some }  C > 0, m \in \mathbb  N \text{ depending on } \varphi. &
\end{eqnarray*}
$\mathscr{H}$ is considered as a space of ``random variables''. In this case we denote $X\in \mathscr{H}$.  We also denote $C_{b,Lip}(\mathbb R^n)$ the space of bounded  Lipschitz
functions.

\begin{definition}\label{def1.1} A  sub-linear expectation $\Sbep$ on $\mathscr{H}$  is a function $\Sbep: \mathscr{H}\to \overline{\mathbb R}$ satisfying the following properties: for all $X, Y \in \mathscr H$, we have
\begin{description}
  \item[\rm (a)]  Monotonicity: If $X \ge  Y$ then $\Sbep [X]\ge \Sbep [Y]$;
\item[\rm (b)] Constant preserving: $\Sbep [c] = c$;
\item[\rm (c)] Sub-additivity: $\Sbep[X+Y]\le \Sbep [X] +\Sbep [Y ]$ whenever $\Sbep [X] +\Sbep [Y ]$ is not of the form $+\infty-\infty$ or $-\infty+\infty$;
\item[\rm (d)] Positive homogeneity: $\Sbep [\lambda X] = \lambda \Sbep  [X]$, $\lambda\ge 0$.
 \end{description}
 Here $\overline{\mathbb R}=[-\infty, \infty]$, $0\cdot \infty$ is defined to be $0$. The triple $(\Omega, \mathscr{H}, \Sbep)$ is called a sub-linear expectation space. Given a sub-linear expectation $\Sbep $, let us denote the conjugate expectation $\cSbep$of $\Sbep$ by
$$ \cSbep[X]:=-\Sbep[-X], \;\; \forall X\in \mathscr{H}. $$
\end{definition}

From the definition, it is easily shown that    $\cSbep[X]\le \Sbep[X]$, $\Sbep[X+c]= \Sbep[X]+c$ and $\Sbep[X-Y]\ge \Sbep[X]-\Sbep[Y]$ for all
$X, Y\in \mathscr{H}$ with $\Sbep[Y]$ being finite. Further, if $\Sbep[|X|]$ is finite, then $\cSbep[X]$ and $\Sbep[X]$ are both finite. By Theorem 1.2.1 of Peng (2019), there exists a family of finite additive linear expectations $E_{\theta}: \mathscr{H}\to \overline{\mathbb R}$ indexed by $\theta\in \Theta$, such that
\begin{equation}\label{linearexpression} \Sbep[X]=\max_{\theta\in \Theta} E_{\theta}[X] \; \text{ for } X \in \mathscr{H} \text{ with } \Sbep[X] \text{ being finite}. \end{equation}
Moreover, for each $X\in \mathscr{H}$, there exists $\theta_X\in \Theta$ such that $\Sbep[X]=E_{\theta_X}[X]$ if $\Sbep[X]$ is finite.

\begin{definition} ({\em See Peng (2008, 2019)})

\begin{description}
  \item[ \rm (i)] ({\em Identical distribution}) Let $\bm X_1$ and $\bm X_2$ be two $n$-dimensional random vectors defined
respectively in sub-linear expectation spaces $(\Omega_1, \mathscr{H}_1, \Sbep_1)$
  and $(\Omega_2, \mathscr{H}_2, \Sbep_2)$. They are called identically distributed, denoted by $\bm X_1\overset{d}= \bm X_2$  if
$$ \Sbep_1[\varphi(\bm X_1)]=\Sbep_2[\varphi(\bm X_2)], \;\; \forall \varphi\in C_{b,Lip}(\mathbb R^n). $$
 A sequence $\{X_n;n\ge 1\}$ of random variables is said to be identically distributed if $X_i\overset{d}= X_1$ for each $i\ge 1$.
\item[\rm (ii)] ({\em Independence})   In a sub-linear expectation space  $(\Omega, \mathscr{H}, \Sbep)$, a random vector $\bm Y =
(Y_1, \ldots, Y_n)$, $Y_i \in \mathscr{H}$ is said to be independent to another random vector $\bm X =
(X_1, \ldots, X_m)$ , $X_i \in \mathscr{H}$ under $\Sbep$  if for each test function $\varphi\in C_{l,Lip}(\mathbb R^m \times \mathbb R^n)$
we have
$ \Sbep [\varphi(\bm X, \bm Y )] = \Sbep \big[\Sbep[\varphi(\bm x, \bm Y )]\big|_{\bm x=\bm X}\big],$
whenever $\overline{\varphi}(\bm x):=\Sbep\left[|\varphi(\bm x, \bm Y )|\right]<\infty$ for all $\bm x$ and
 $\Sbep\left[|\overline{\varphi}(\bm X)|\right]<\infty$.
 \item[\rm (iii)] ({\em Independent random variables}) A sequence of random variables $\{X_n; n\ge 1\}$
 is said to be independent, if  $X_{i+1}$ is independent to $(X_{1},\ldots, X_i)$ for each $i\ge 1$.
 \end{description}
\end{definition}
It is easily seen that,  if $\{X_1,\ldots,X_n\}$ are independent, then $\Sbep[\sum_{i=1}^n X_i]=\sum_{i=1}^n\Sbep[X_i]$.

Next, we consider the capacities corresponding to the sub-linear expectations. Let $\mathcal G\subset\mathcal F$. A function $V:\mathcal G\to [0,1]$ is called a capacity if
$$ V(\emptyset)=0, \;V(\Omega)=1 \; \text{ and } V(A)\le V(B)\;\; \forall\; A\subset B, \; A,B\in \mathcal G. $$
It is called to be sub-additive if $V(A\bigcup B)\le V(A)+V(B)$ for all $A,B\in \mathcal G$  with $A\bigcup B\in \mathcal G$.

 Let $(\Omega, \mathscr{H}, \Sbep)$ be a sub-linear expectation space.  We denote   $(\Capc,\cCapc)$ to be a pair of  capacities with the properties that
\begin{equation}\label{eq1.3}
   \Sbep[f]\le \Capc(A)\le \Sbep[g]\;\;
\text{ if } f\le I_A\le g, f,g \in \mathscr{H} \text{ and } A\in \mathcal F,
\end{equation}
 $$ \Capc \text{ is sub-additive}  $$
and $\cCapc(A):= 1-\Capc(A^c)$, $A\in \mathcal F$.
It is obvious that
\begin{equation}\label{eq1.4}
  \cCapc(A\bigcup B)\le \cCapc(A)+\Capc(B).
\end{equation}
We call $\Capc$ and $\cCapc$ the upper and the lower capacity, respectively. In general, we can choose $\Capc$ as
\begin{equation}\label{eq1.5} \Capc(A):=\inf\{\Sbep[\xi]: I_A\le \xi, \xi\in\mathscr{H}\},\;\; \forall A\in \mathcal F.
\end{equation}
To distinguish this   capacity from others, we denote it by $\upCapc$, and $\lowCapc(A)=1-\upCapc(A)$. $\upCapc$ is the largest capacity satisfying \eqref{eq1.3}.

When there exists  a family of probability measure on $(\Omega,\mathscr{F})$ such that
\begin{equation}\label{eq1.7} \Sbep[X]=\sup_{P\in \mathscr{P}}P[X]=:\sup_{P\in \mathscr{P}}\int  XdP ,
\end{equation} $\Capc$ can be defined as
\begin{equation}\label{eq1.6} \Capc(A)=\sup_{P\in \mathscr{P}}P(A).
\end{equation}
We denote this capacity by $\Capc^{\mathscr{P}}$, and $\cCapc^{\mathscr{P}}(A)=1-\Capc^{\mathscr{P}}(A)$.

 Also, we define the  Choquet integrals/expecations $(C_{\Capc},C_{\cCapc})$  by
$$ C_V[X]=\int_0^{\infty} V(X\ge t)dt +\int_{-\infty}^0\left[V(X\ge t)-1\right]dt $$
with $V$ being replaced by $\Capc$ and $\cCapc$ respectively.
If $\Capc_1$ on the sub-linear expectation space $(\Omega_1,\mathscr{H}_1,\Sbep_1)$ and $\Capc_2$ on the sub-linear expectation space $(\Omega_2,\mathscr{H}_2,\Sbep_2)$  are two capacities having the property \eqref{eq1.3}, then for any random variables $X_1\in \mathscr{H}_1$ and $\tilde X_2\in \mathscr{H}_2$ with $X_1\overset{d}=\tilde X_2$, we have
\begin{equation}\label{eqV-V}
\Capc_1(X_1\ge x+\epsilon)\le \Capc_2(\tilde X_2\ge x)\le \Capc_1(X_1\ge x-\epsilon)\;\; \text{ for all } \epsilon>0 \text{ and } x,
\end{equation}
and so
$$
\Capc_1(X_1\ge x)=\Capc_2(\tilde X_2\ge x)\;\;  \text{ for all  except countable many }  x,
$$
\begin{equation}\label{eqV-V3} C_{\Capc_1}(X_1)=C_{\Capc_2}(\tilde X_2).
\end{equation}
In particular, if $\Capc_1$ and $\Capc_2$ are two capacities having the property \eqref{eq1.3}, then for any random variable $X\in \mathscr{H}$,
$$
\Capc_1(X\ge x+\epsilon)\le \Capc_2(X\ge x)\le \Capc_1(X\ge x-\epsilon)\;\; \text{ for all } \epsilon>0 \text{ and } x
$$
and
$$ C_{\Capc_1}(X)=C_{\Capc_2}(X).
$$

Finally, for real numbers $x$ and $y$, denote $x\vee y=\max(x,y)$, $x\wedge y=\min(x,y)$, $x^+=\max(0,x)$, $x^-=\max(0,-x)$ and $\log x=\ln \max(e,x)$. For a random variable $X$, because $XI\{|X|\le c\}$  may be not in $\mathscr{H}$, we will truncate it in the form $(-c)\vee X\wedge c$ denoted by $X^{(c)}$.

 \section{Exponential  inequalities} \label{sectIneq}
\setcounter{equation}{0}
Exponential inequalities    and Kolmogorov's converse exponential inequality are basic tools for establishing the LIL. In this section, we give the   exponential inequalities   under both the upper capacity $\Capc$ and the lower capacity $\cCapc$ and Kolmogorov's converse exponential inequalities under the upper capacity $\Capc$.
   The next lemma gives   the Kolmogorov-type exponential inequalities for maximum sums of independent random variables.
\begin{lemma}\label{lemExpIneq}  Let $\{X_1,\ldots, X_n\}$ be a sequence  of   independent random variables
in the sub-linear expectation space $(\Omega,\mathscr{H},\Sbep)$. Set $S_n=\sum_{i=1}^n  X_i$,  $B_n^2=\sum_{i=1}^n \Sbep[X_i^2]$, $b_n^2=\sum_{i=1}^n \cSbep[X_i^2]$,  and $A_n(p,y)=\sum_{i=1}^n\Sbep[(X_i^+\wedge y )^p]$, $p\ge 2$. Denote
$$ B_{n,y}^2=\sum_{i=1}^n \Sbep[(X_i\wedge y)^2], \;\;  b_{n,y}^2=\sum_{i=1}^n \cSbep[(X_i\wedge y)^2], \;\;  y>0. $$
\begin{description}
  \item[\rm (I)]
  For all $x,y>0$,
\begin{align}\label{eqExpIneq.3}
& \Capc\Big( \max_{k\le n} (S_k-\Sbep[S_k])\ge x\Big)\;\; \left(\text{ resp. } \cCapc\Big( \max_{k\le n} (S_k-\cSbep[S_k])\ge x\Big)\right)\nonumber\\
\le & \Capc\big(\max_{k\le n}X_k> y\big) + \exp\left\{-\frac{x^2}{2(xy+B_{n,y}^2)}\Big(1+\frac{2}{3}\ln \big(1+\frac{xy}{B_{n,y}^2}\big)\Big)\right\};
\end{align}
and for all $p\ge 2$, $x,y>0$ and $0<\delta\le 1$,
\begin{align}\label{eqExpIneq.4}
& \Capc\Big( \max_{k\le n} (S_k-\Sbep[S_k])\ge x\Big)\;\; \left(\text{ resp. } \cCapc\Big( \max_{k\le n} (S_k-\cSbep[S_k])\ge x\Big)\right)\nonumber\\
\le &  \Capc\big(\max_{i\le n} X_i> y \big)
  +2\exp\{p^p\}\Big\{\frac{A_n(p,y)}{y^p} \Big\}^{\frac{\delta x}{10y}}
+\exp\left\{-\frac{x^2}{2(1+\delta) B_{n,y}^2 }\right\}.
\end{align}
 \item[\rm (II)]  For all $x,y>0$,
\begin{align}\label{eqExpIneq.5}
& \cCapc\Big( \max_{k\le n} (S_k-\Sbep[S_k])\ge x\Big)\nonumber\\
\le &  \Capc\big(\max_{k\le n}X_k> y\big)
+ \exp\left\{-\frac{x^2}{2(xy+b_{n,y}^2)}\Big(1+\frac{2}{3}\ln \big(1+\frac{xy}{b_{n,y}^2}\big)\Big)\right\};
\end{align}
and for all $p\ge 2$, $x,y>0$ and $0<\delta\le 1$,
\begin{align}\label{eqExpIneq.6}
& \cCapc\left(\max_{k\le n} (S_k-\Sbep[S_k])\ge x\right) \nonumber\\
\le &\Capc\big(\max_{i\le n} X_i> y \big)
  +2\exp\{p^p\}\Big\{\frac{A_n(p,y)}{y^p} \Big\}^{\frac{\delta x}{10y}}
+\exp\left\{-\frac{x^2}{2 (1+\delta) b_{n,y}^2 }\right\}.
\end{align}
\end{description}
Further, the upper bounds in \eqref{eqExpIneq.4} and \eqref{eqExpIneq.6} can be,  respectively, replaced by
\begin{equation}\label{eqExpIneq.7}
\begin{aligned}
& C_p \delta^{-p} \frac{1}{x^p}\sum_{i=1}^n\Sbep[(X_i^+)^p]+ \exp\left\{-\frac{x^2}{2 (1+\delta) B_n^2 }\right\}, \\
& C_p \delta^{-p} \frac{1}{x^p}\sum_{i=1}^n\Sbep[(X_i^+)^p]+ \exp\left\{-\frac{x^2}{2 (1+\delta) b_n^2 }\right\}.
\end{aligned}
\end{equation}
 \end{lemma}
\begin{remark} \eqref{eqExpIneq.4} and \eqref{eqExpIneq.6} are  Fuk  and  Nagaev (1971)'s type inequalities.
\end{remark}
{\bf Proof}. The upper bound in \eqref{eqExpIneq.3}   for $\Capc(S_n-\Sbep[S_n]\ge x)$ and $\cCapc(S_n-\cSbep[S_n]\ge x)$ are derived by Zhang (2016).  Here, we consider the maximum sums. First, we give the proof  of \eqref{eqExpIneq.5} and \eqref{eqExpIneq.6}.     Let $Y_k=X_k\wedge y$, $T_n=\sum_{i=1}^n (Y_i-\Sbep[X_i])$.  Then $X_k-Y_k=(X_k-y)^+\ge 0$ and
 $\Sbep[Y_k]\le \Sbep[X_k]$.  From the fact  that $\cCapc(A\cup B)\le \cCapc(A)+\Capc(B)$ it follows that
 $$
 \cCapc\left(\max_{k\le n}(S_k-\Sbep[S_k])\ge x\right)\le  \Capc\big(\max_{k\le n}X_k> y\big)+
 \cCapc\left(\max_{k\le n}T_k\ge x\right),
 $$
 and for any $t>0$,  $\varphi(x)=:e^{t(x\wedge y)}$ is a bounded non-decreasing function and belongs to $C_{b,Lip}(\mathbb R)$ since $0\le \varphi^{\prime}(x)\le t e^{ty}$. From
$$e^{tY_k}=1+ tY_k+\frac{e^{tY_k}-1-t Y_k}{Y_k^2}Y_k^2\le 1 +tY_k+\frac{e^{ty}-1-t y}{y^2}Y_k^2   $$
and the facts that $\cSbep[X+Y]\le  \Sbep[X]+\cSbep[Y]$ and $e^y-y-1\ge 0$, we have
$$ \cSbep[e^{t Y_k}]\le 1+t\Sbep[X_k]+\frac{e^{ty}-1-t y}{y^2} \cSbep[Y_k^2]\le\exp\left\{t\Sbep[X_k]+\frac{e^{ty}-1-t y}{y^2}\cSbep[Y_k^2]\right\}, $$
$$\cSbep[e^{t (T_k-T_{k-1})}] \le\exp\left\{\frac{e^{ty}-1-t y}{y^2}\cSbep[Y_k^2]\right\}. $$
Write
$$ U_0=1, \;\; U_k=\exp\Big\{-\frac{e^{ty}-1-t y}{y^2} b_{k,y}^2 \Big\}e^{t T_k},\;\; k=1,\cdots, k_n. $$
Then
\begin{align}\label{eqproofExqI.1} & \cSbep\left[U_k-U_{k-1}|X_1,\ldots,X_{k-1}]\right]\nonumber\\
= & U_{k-1} \cSbep\left[\exp\Big\{-\frac{e^{ty}-1-t y}{y^2} \cSbep[Y_k^2] \Big\}e^{t (T_k-T_{k-1})}-1\right]\le 0.
\end{align}
Here and in the following, for a bounded Lipschitz function $Y=f(X_1,\ldots, X_n)$ of $X_1,\ldots, X_n$, $ \Sbep\left[Y|X_1,\ldots,X_{k-1}]\right]$ and  $\cSbep\left[Y|X_1,\ldots,X_{k-1}]\right]$ are,  respectively, defined  by
$$\Sbep\left[Y|X_1,\ldots,X_{k-1}]\right]=\left(\Sbep\left[f(x_1,\ldots,x_{k-1}, X_k,\ldots,X_n)\right]\right)\Big|_{x_1=X_1,\ldots, x_{d-1}=X_{d-1}}$$
and
$$\cSbep\left[Y|X_1,\ldots,X_{k-1}]\right]=\left(\cSbep\left[f(x_1,\ldots,x_{k-1}, X_k,\ldots,X_n)\right]\right)\Big|_{x_1=X_1,\ldots, x_{d-1}=X_{d-1}}.$$
 For any $\alpha>0$ and   given $\beta\in (0,\alpha)$, let $f(x)$ be a continuous function with bounded derivation such that $I\{x\le \alpha-\beta\}\le f(x) \le I\{x<\alpha\}$. Define $f_0=1$, $f_k=f(U_1)\cdots f(U_k)$. Then
\begin{align*}
&f_0U_0+\sum_{k=1}^n f_{k-1}\big(U_k-U_{k-1}\big)=f_nU_n+\sum_{k=1}^n f_{k-1}\big(1-f(U_k)\big)U_k\\
\ge & f_nU_n+\sum_{k=1}^n f_{k-1}\big(1-f(U_k)\big)(\alpha-\beta)=(\alpha-\beta)(1-f_n)+f_nU_n\\
\ge & (\alpha-\beta)I\{\max_{k\le n}U_k\ge \alpha\}.
\end{align*}
By the independence property,
\begin{align*}
&\cSbep\left[f_0U_0+\sum_{k=1}^n f_{k-1}\big(U_k-U_{k-1}\big)\right]\\
= & \cSbep\left[f_0U_0+\sum_{k=1}^{n-1} f_{k-1}\big(U_k-U_{k-1}\big)+f_{n-1}\cSbep\big[U_n-U_{n-1}\big|X_1,\ldots,X_{n-1}\big]\right]\\
\le & \cSbep\left[f_0U_0+\sum_{k=1}^{n-1} f_{k-1}\big(U_k-U_{k-1}\big)\right]\le \cdots\le \cSbep[f_0U_0].
\end{align*}
It follows that
$$ (\alpha-\beta)\cCapc(\max_{k\le n}U_k\ge \alpha)\le \cSbep[f_0U_0]=\cSbep[U_0]. $$
By letting $\beta\to 0$, we have
\begin{equation}\label{eqproofExqI.2}
\cCapc\left(\max_{k\le n}U_k\ge \alpha\right)\le \frac{\cSbep[U_0]}{\alpha}=\frac{1}{\alpha}.
\end{equation}
Note
$$ \exp\left\{t\max_{k\le n}T_k\right\}\le \max_{k\le n}U_k \exp\Big\{\frac{e^{ty}-1-t y}{y^2} b_{n,y}^2 \Big\}. $$
Hence by \eqref{eqproofExqI.2},
\begin{align*}
\cCapc\left(\max_{k\le n} T_k\ge x\right)  \le & \cCapc\left(\max_{k\le n} U_k\ge \exp\Big\{tx-\frac{e^{ty}-1-t y}{y^2} b_{n,y}^2 \Big\}\right)\\
 \le &  \exp\left\{-tx+\frac{e^{ty}-1-t y}{y^2} b_{n,y}^2\right\}.
 \end{align*}
Choosing $t=\frac{1}{y}\ln \big(1+\frac{xy}{b_{n,y}^2}\big)$ yields
\begin{align}\label{eqproofExqI.3}
\cCapc\left(\max_{k\le n}T_k \ge x\right)  \le  \exp\left\{\frac{x}{y}-\frac{x}{y}\Big(\frac{b_{n,y}^2}{xy}+1\Big)\ln\Big(1+\frac{xy}{b_{n,y}^2}\Big)\right\}.
 \end{align}
 Applying the elementary inequality
 $$ \ln (1+t)\ge \frac{t}{1+t}+\frac{t^2}{2(1+t)^2}\big(1+\frac{2}{3} \ln (1+t)\big)$$
 yields
 $$ \Big(\frac{b_{n,y}^2}{xy}+1\Big)\ln\Big(1+\frac{xy}{b_{n,y}^2}\Big)
 \ge 1+\frac{xy}{2(xy+b_{n,y}^2)}\Big(1+\frac{2}{3}\ln\big(1+\frac{xy}{b_{n,y}^2}\big)\Big). $$
Thus, \eqref{eqExpIneq.5} is proved.

 Next we show \eqref{eqExpIneq.6}. If $xy\le \delta b_{n,y}^2$, then
 $$\frac{x^2}{2(xy+b_{n,y}^2)}\Big(1+\frac{2}{3}\ln\big(1+\frac{xy}{b_{n,y}^2}\big)\Big)\ge \frac{x^2}{2(1+\delta) b_{n,y}^2 }. $$
If  $xy\ge \delta b_{n,y}^2$,
then
 $$\frac{x^2}{2(xy+b_{n,y}^2)}\Big(1+\frac{2}{3}\ln\big(1+\frac{xy}{b_{n,y}^2}\big)\Big)\ge \frac{x}{2(1+1/\delta)y}. $$
 It follows that
\begin{equation}\label{eqproofExqI.4} \cCapc\left(\max_{k\le n}T_k\ge x\right)\le \exp\left\{-\frac{x^2}{2(1+\delta) b_{n,y}^2 }\right\}
 +\exp\left\{-\frac{x}{2 (1+1/\delta)y}\right\}
 \end{equation}
 by \eqref{eqproofExqI.3}. For $z>0$, let
 $$\beta(z)=\beta_p(z)=\frac{1}{z^p}\sum_{k=1}^n\Sbep[(X_k^+\wedge z)^p], $$
 and choose
 $$ \rho=1\wedge \frac{1}{(1+\delta)\ln (1/\beta(z))}, \;\; y= \rho z. $$
 Notice $b_{n,y}^2\le b_{n,z}^2$. Then by \eqref{eqproofExqI.4},
 \begin{align*}
  &\cCapc\big(\max_{k\le n}(S_k-\Sbep[S_k]) \ge (1+\delta)x\big)\le  \cCapc\big(\max_{k\le n}T_k\ge  x\big)+\Capc\big(\sum_{i=1}^n(X_i- y)^+\ge \delta x\big) \\
  \le &\exp\left\{-\frac{x^2}{2(1+\delta) b_{n,z}^2 }\right\}+\Big(\beta(z)\Big)^{\delta x/2z}+\Capc\big(\max_{i\le n} X_i> z \big)+\Capc\big(\sum_{i=1}^n(X_i-\rho z)^+\wedge z \ge \delta x\big).
  \end{align*}

  On the other hand, for $t>0$,
  \begin{align*}
  &\Capc\big(\sum_{i=1}^n(X_i-\rho z)^+\wedge z \ge \delta x\big)=
 \Capc\left(\sum_{i=1}^n\left[\Big(\frac{X_i}{z}- \rho\big)^+\wedge 1\right] \ge  \frac{\delta x}{z} \right)\\
 \le & e^{- t\frac{\delta x}{z} }\Sbep\exp\left\{t\sum_{i=1}^n\left[\Big(\frac{X_i}{z}- \rho\big)^+\wedge 1\right]\right\}
 = e^{- t\frac{\delta x}{z}}\prod_{i=1}^n \Sbep\exp\left\{t\left[\Big(\frac{X_i}{z}- \rho\big)^+\wedge 1\right]\right\}\\
 \le & e^{- t\frac{\delta x}{z}}\prod_{i=1}^n\left[ 1+e^t \Capc(X_i\ge \rho z)\right]
 \le \exp\left\{ - t\frac{\delta x}{z} +e^t \frac{\beta(z)}{\rho^p}\right\}.
 \end{align*}
 By taking the minimum over $t\ge 0$, it follows that
   $$ \Capc\big(\sum_{i=1}^n(X_i-\rho z)^+\wedge z \ge \delta x\big)\le \exp\left\{\frac{\delta x}{z}\Big(1-\ln \frac{\delta x}{z}+\ln \frac{\beta(z)}{\rho^p}\Big)\right\}.
   $$
  Assume  $\beta(z)<1$.   When $\rho=1/[(1+\delta)\ln \frac{1}{\beta(z)}]$, by the fact that $\sqrt{x}(\ln \frac{1}{x})^p\le (2pe^{-1})^p$ we have
   \begin{align*}
   &\frac{\delta x}{z}\Big(1-\ln \frac{\delta x}{z}+\ln \frac{\beta(z)}{\rho^p}\Big) \\
   =& \frac{\delta x}{z}\Big(1-\ln \frac{\delta x}{z}+\ln \beta(z)+p\ln\ln \frac{1}{\beta(z)}+ p\ln  (1+\delta) \Big)\\
   \le & \frac{\delta x}{z}\Big( -\ln \frac{\delta x}{z}+1+p\ln(2(1+\delta)pe^{-1})\Big)+\frac{\delta x}{2z}\ln \beta(z) \\
   \le & (2(1+\delta)pe^{-1})^p+\frac{\delta x}{2z}\ln \beta(z).
   \end{align*}
   The last inequality is due to $\max_{x\ge 0}\{x(C+1-\ln x)\}=e^C$.  When $\rho=1$,
   $$\frac{\delta x}{z}\Big(1-\ln \frac{\delta x}{z}+\ln \frac{\beta(z)}{\rho^p}\Big)\le 1+\frac{\delta x}{z}\ln   \beta(z)
   \le 1+\frac{\delta x}{2z}\ln   \beta(z). $$
  It follows that
 \begin{align*}
   &\cCapc\big(\max_{k\le n}(S_k-\Sbep[S_k])\ge (1+\delta)x\big)\\
  \le &\exp\left\{-\frac{x^2}{2(1+\delta) b_{n,z}^2 }\right\}+\Capc\big(\max_{i\le n} X_i> z \big)\\
  & +2\exp\left\{(2(1+\delta)pe^{-1})^p\right\}\Big(\beta(z)\Big)^{\delta x/2z}.
  \end{align*}
 Therefore,
   \begin{align*}
   &\cCapc\big(\max_{k\le n}(S_k-\Sbep[S_k])\ge  x\big)\\
  \le &\exp\left\{-\frac{x^2}{2 (1+\delta)^3 b_{n,z}^2  }\right\}+\Capc\big(\max_{i\le n} X_i> z \big)\\
  & +2\exp\left\{(2(1+\delta)pe^{-1})^p\right\}\Big(\beta(z)\Big)^{\frac{\delta x}{2(1+\delta)z}}.
  \end{align*}
 For $0<\delta^{\prime}\le 1$, let $\delta=\sqrt[3]{1+\delta^{\prime}}-1$. Then  $\frac{\delta}{1+\delta} \ge \frac{\delta^\prime}{5}$, $2(1+\delta)<e$. It follows that
\begin{align*}
   &\cCapc\big(\max_{k\le n}(S_k-\Sbep[S_k])\ge  x\big)\\
  \le & \exp\left\{-\frac{x^2}{2 (1+\delta^{\prime}) b_{n,z}^2 }\right\} +\Capc\big(\max_{i\le n} X_i> z \big)
  +2\exp\{p^p\}\Big(\beta(z)\Big)^{\frac{\delta^{\prime} x}{10z}}.
 \end{align*}
   If $\beta(z)\ge 1$, then the above inequality is obvious. \eqref{eqExpIneq.6} is proved.

For (I), it is sufficient to notice that
$$ \cSbep[e^{t Y_k}]\le 1+t\cSbep[X_k]+\frac{e^{ty}-1-t y}{y^2} \Sbep[Y_k^2], $$
$$\cSbep[e^{t (Y_k-\cSbep[X_k])}] \le\exp\left\{\frac{e^{ty}-1-t y}{y^2}\Sbep[Y_k^2]\right\}  $$
and
$$ \Sbep[e^{t Y_k}]\le 1+t\Sbep[X_k]+\frac{e^{ty}-1-t y}{y^2} \Sbep[Y_k^2], $$
$$\Sbep[e^{t (Y_k-\Sbep[X_k])}] \le\exp\left\{\frac{e^{ty}-1-t y}{y^2}\Sbep[Y_k^2]\right\}.$$

For \eqref{eqExpIneq.7}, it is sufficient to choose $y=\delta x/10$ and notice that $b_{n,y}^2\le b_n^2$, $B_{n,y}^2\le B_n^2$,
$$ \Capc\big(\max_{i\le n} X_i> y \big)\le \frac{A_n(p,y)}{y^p}, \; y>0. \;\; \Box $$

\bigskip


The following lemma is an analogue of Kolmogorov's converse exponential inequality.
   \begin{lemma} \label{lemExpIneq2}Let   $\{X_{n,i};i=1,\ldots, k_n\}$ be an array  of independent random variables
in the sub-linear expectation space $(\Omega,\mathscr{H},\Sbep)$ with $s_n^2=\sum_{i=1}^{k_n}\Sbep[X_{n,i}^2]$.  Let $\{x_n\}$ be a sequence of positive numbers with   $x_n\to \infty$ and $x_n/s_n\to 0$.
Suppose
\begin{equation}\label{eqlemExpIneq2.1} \frac{ \sum_{i=1}^{k_n}|\Sbep[X_{n,i}]|}{s_nx_n}\to 0, \;\; \frac{ \sum_{i=1}^{k_n}|\cSbep[X_{n,i}]|}{s_nx_n}\to 0
\end{equation}
and there  exists    a positive number  $\alpha$ such that
$$ |X_{n,i}|\le \alpha \frac{s_n}{x_n}, \;\; i=1,\ldots, k_n. $$
Then for any $\gamma>0$, there exists a positive constant $\pi(\gamma)$ (small enough) such that
\begin{equation}\label{eqlemExpIneq2.3} \liminf_{n\to \infty} x_n^{-2}\ln \Capc\left(\sum_{i=1}^{k_n}X_{n,i}\ge z s_n x_n\right)\ge -\frac{z^2}{2}(1+\gamma) \text{ for all } 0<z\alpha\le \pi(\gamma).
\end{equation}
   \end{lemma}

 {\bf Proof}. We use an  argument analogues to Stout (1974) (cf. Petrov (1995, Page 241-243)).  First, it is easily seen that
 $$\frac{1}{s_n^2}\sum_{i=1}^{k_n}\left|\Sbep\big[(X_{n,i}-\Sbep[X_{n,i}])^2]-\Sbep[X_{n,i}^2]\right|\le  \frac{3\alpha}{s_nx_n}
 \sum_{i=1}^{k_n}|\Sbep[X_{n,i}]| \to 0. $$
 Without loss of generality, we can assume that $\Sbep[X_{n,i}]=0$, $i=1,\ldots, k_n$. Otherwise, we instead consider $X_{n,i}-\Sbep[X_{n,i}]$ instead.

 Let $S_n=\sum_{i=1}^{k_n}X_{n,i}$ and $q_n(y)=\Capc\left(S_n\ge y s_nx_n)\right)$. Then by \eqref{eqExpIneq.3},
 \begin{equation}\label{eqPlemExpIneq2.2} q_n(y)\le \exp\left\{-\frac{y^2x_n^2}{2(y\alpha+1)}\right\}, \;\; y>0.
 \end{equation}
 For any $t>0$ with $t \alpha <1/32$, we have
\begin{align*}
\exp\{tX_{n,i}\frac{x_n}{s_n}\}\ge & 1+ tX_{n,i}\frac{x_n}{s_n}+\frac{t^2}{2}X_{n,i}^2\frac{x_n^2}{s_n^2}\left(1-\sum_{i=3}^{\infty}\frac{2(t\alpha)^{i-2}}{i!}\right)\\
\ge & 1  + tX_{n,i}\frac{x_n}{s_n}+\frac{t^2}{2}X_{n,i}^2\frac{x_n^2}{s_n^2}(1-t\alpha/2).
\end{align*}
Then from the fact that $\Sbep[X+Y]\ge \cSbep[X]+\Sbep[Y]$ it follows that
$$
\Sbep\left[\exp\{tX_{n,i}\frac{x_n}{s_n}\}\right]\ge    1+t \frac{x_n}{s_n}\cSbep[X_{n,i}]+\frac{t^2}{2}\frac{x_n^2}{s_n^2}(1-t\alpha/2))\Sbep[X_{n,i}^2].
$$
Applying $\ln (1+x)\ge x-x^2$ ($x\ge -1/4$) yields
\begin{align*}
 \ln \Sbep\left[\exp\{tX_{n,i}\frac{x_n}{s_n}\}\right]\ge
& t \frac{x_n}{s_n}\cSbep[X_{n,i}]+\frac{t^2}{2}\frac{x_n^2}{s_n^2}(1-t\alpha/2)\Sbep[X_{n,i}^2] \\
&-\left(t \frac{x_n}{s_n}|\cSbep[X_{n,i}]|(t\alpha+(t\alpha)^2)+\frac{t^2}{2}\frac{x_n^2}{s_n^2}\Sbep[X_{n,i}^2](1-t\alpha/2)^2\frac{(t\alpha)^2}{2}\right).
\end{align*}
It follows that
\begin{align*}
 x_n^{-2}\ln \Sbep\left[\exp\{tS_n \frac{x_n}{s_n}\}\right]\ge
& -t \frac{\sum_{i=1}^{k_n}|\cSbep[X_{n,i}]|}{s_n x_n }(1+2) +\frac{t^2}{2} (1-t\alpha)  \\
\to & \frac{t^2}{2}   (1-t\alpha).
\end{align*}
Note
$$
\Sbep\left[\exp\{tS_n \frac{x_n}{s_n}\}\right]\le   C_{\Capc}\left(\exp\{tS_n \frac{x_n}{s_n}\}\right)=\int_0^{\infty}\Capc\left(\exp\{tS_n \frac{x_n}{s_n}\}> y\right)dy
$$
It follows that
\begin{equation}\label{eqPlemExpIneq2.7} \liminf_{n\to \infty} x_n^{-2}\ln C_{\Capc}\left(\exp\{tS_n \frac{x_n}{s_n}\}\right)\ge   \frac{t^2}{2}(1-t\alpha)\;\;\text{ for all } 0<t\alpha<1/32.
\end{equation}

Now, for $\delta<1/4$, let $t=z/(1-\delta)$. Then
\begin{align}\label{eqPlemExpIneq2.8}
&C_{\Capc}\left(\exp\{tS_n \frac{x_n}{s_n}\}\right)=\int_{-\infty}^{\infty} t x_n^2e^{tx_n^2y}q_n(y)dy\nonumber\\
=& \left(\int_{-\infty}^0+\int_0^{t(1-\delta)}+\int_{t(1-\delta)}^{t(1+\delta)}+\int_{t(1+\delta)}^{8t}+\int_{8t}^{\infty}\right)t x_n^2e^{tx_n^2y}q_n(y)dy\nonumber\\
=:& I_1+I_2+I_3+I_4+I_5.
\end{align}
It is obvious that
\begin{equation}\label{eqPlemExpIneq2.9} I_3\le 2t^2x_n^2 \exp\{t^2x_n^2(1+\delta)\}q_n(t(1-\delta))=2t^2x_n^2 \exp\{t^2x_n^2(1+\delta)\}q_n(z)
\end{equation}
and
 \begin{equation}\label{eqPlemExpIneq2.10} I_1\le \int_{-\infty}^0t x_n^2e^{tx_n^2y}dy\le 1.
\end{equation}
Assume that $8t\alpha\le 1$.    by \eqref{eqPlemExpIneq2.2}, if   $y\alpha\ge 1$, then
$$ e^{tx_n^2y}q_n(y)\le \exp\left\{tx_n^2y-\frac{yx_n^2}{4\alpha}\right\}\le e^{-tx_n^2y}, $$
and, if $8t\le y\le 1/\alpha$,
$$ e^{tx_n^2y}q_n(y)\le \exp\left\{tx_n^2y-\frac{y^2x_n^2}{4}\right\}\le e^{-tx_n^2y}. $$
It follows that
\begin{equation}\label{eqPlemExpIneq2.12} I_5\le \int_{8t}^{\infty} tx_n^2e^{-tx_n^2y} dy \le 1, \;\; 8t\alpha\le 1.
\end{equation}
Now, consider $I_2$ and $I_4$. Choose a positive constant $\beta$. Then if $y\alpha\le \beta<\delta$, then
$$ q_n(y)\le \exp\left\{-\frac{y^2x_n^2}{2(1+\beta)}\right\},\; \text{ if } y\le 8t \text{ and }  8t\alpha \le \beta. $$
Let $\psi(y)=ty-\frac{y^2}{2(1+\beta)}$. Thus we arrive the inequality
$$ I_2+I_4\le tx_n^2 \int_D e^{\psi(y)x_n^2}dy \; \text{ with } \; D=(0,t(1-\delta))\cup (t(1+\delta),8t).  $$
The function $\psi(y)$ has a maximum at the point $y=t(1+\beta)$ which lies in the interval $(t(1-\delta),t(1+\delta))$. Therefore,
\begin{align*}
 & \sup_{y\in D}\psi(y)=  \max\{\psi(t(1-\delta)),\psi(t(1+\delta))\}\\
=&\frac{t^2}{2}\Big(1-\delta^2+(1+\delta)^2\frac{\beta}{1+\beta}\Big)
\le   \frac{t^2}{2}\Big(1-\frac{\delta^2}{2}-\frac{\delta^2}{2}\frac{\beta}{1+\beta}\Big),
\end{align*}
if $\beta=\delta^2/(2(1+\delta)^2)$. It follows that
\begin{equation}\label{eqPlemExpIneq2.14} I_2+I_4\le 8t^2x_n^2\exp\left\{\max_{y\in D}\psi(y) x_n^2\right\} \le\frac{1}{4} \exp\left\{\frac{t^2x_n^2}{2}\Big(1-\frac{\delta^2}{2}\Big)\right\}
\end{equation}
for $n$ large enough if $8t\alpha\le \beta=\delta^2/(2(1+\delta)^2)$.
On the other hand, if $t\alpha\le \beta=\delta^2/2(1+\delta)^2$, it follows from \eqref{eqPlemExpIneq2.7} that
\begin{equation}\label{eqPlemExpIneq2.15}  C_{\Capc}\left(\exp\{tS_n \frac{x_n}{s_n}\}\right)\ge  \exp\left\{ \frac{t^2x_n^2}{2}\Big(1-\frac{\delta^2}{2}\Big)\right\}\ge 8
\end{equation}
for $n$ large enough. It follows from \eqref{eqPlemExpIneq2.9}-\eqref{eqPlemExpIneq2.15} that, for $0<t\alpha\le \delta^2/(16(1+\delta)^2)$,
\begin{equation}\label{eqPlemExpIneq2.16} I_1+I_2+I_4+I_5\le \frac{1}{2}C_{\Capc}\left(\exp\{tS_n \frac{x_n}{s_n}\}\right),\end{equation}
and therefore,
\begin{equation}\label{eqPlemExpIneq2.17} 2t^2x_n^2 \exp\{t^2x_n^2(1+\delta)\}q_n(z)\ge I_3\ge \frac{1}{2}C_{\Capc}\left(\exp\{tS_n \frac{x_n}{s_n}\}\right)\ge  \frac{1}{2}\exp\left\{ \frac{t^2x_n^2}{2}\Big(1-\frac{\delta^2}{2}\Big)\right\}
\end{equation}
when $n$ is large enough. It follows that
\begin{align*}
 \liminf_{n\to \infty} x_n^{-2}\ln q_n(z)\ge & \frac{t^2}{2} \Big(1-\frac{\delta^2}{2}\Big)-t^2(1+\delta)\\
 > & -\frac{z^2}{2}\frac{(1+\delta)^2}{(1-\delta)^2},
\;\;\text{ if }  0<z\alpha<\delta^2/(16(1+\delta)^2), \;\delta<1/4.
\end{align*}
At last, for every $\gamma>0$, choose $0<\delta<1/4$ such that $\frac{(1+\delta)^2}{(1-\delta)^2}\le 1+\gamma$. Then \eqref{eqlemExpIneq2.3} holds with $\pi(\gamma)=\delta^2/(16(1+\delta)^2)$.
The proof is completed.  $\Box$

We conjecture that   for the lower capacity $\cCapc$, we have an analogue Kolmogorov's converse exponential inequality.

\begin{conjecture}   Let   $\{X_{n,i};i=1,\ldots, k_n\}$ be an array  of independent random variables in
the sub-linear expectation space $(\Omega, \mathscr{H}, \Sbep)$ with $\underline{s}_n^2=\sum_{i=1}^{k_n}\cSbep[X_{n,i}^2]$.  Let $x_n$ be a sequence of positive numbers with   $x_n\to \infty$ and $x_n/\underline{s}_n\to 0$.
Suppose
$$ \frac{ \sum_{i=1}^{k_n}|\Sbep[X_{n,i}]|}{\underline{s}_nx_n}\to 0, \;\; \frac{ \sum_{i=1}^{k_n}|\cSbep[X_{n,i}]|}{\underline{s}_nx_n}\to 0 $$
and there  exists    a positive number  $\alpha$ such that
$$ |X_{n,i}|\le \alpha  \underline{s}_n/ x_n, \;\; i=1,\ldots, k_n. $$
Then for any $\gamma>0$, there exists a positive constant $\pi(\gamma)$ (small enough) such that
\begin{equation}\label{eqConjecture} \liminf_{n\to \infty} x_n^{-2}\ln \cCapc\left(\sum_{i=1}^{k_n}X_{n,i}\ge z \underline{s}_n x_n\right)\ge -\frac{z^2}{2}(1+\gamma) \text{ for all } 0<z\alpha\le \pi(\gamma).
\end{equation}
\end{conjecture}

 It seems that it is not an easy task  to obtain the lower bound of the tail capacity under $\cCapc$. Recently, Peng, Yang and Yao (2020) and  Peng and Zhou (2020) studied the tail behavior of the G-normal distribution by analyzing a nonlinear heat equation. Let $\xi\sim N(0,[\underline{\sigma}^2,\overline{\sigma}])$ in sub-linear expectation space $(\widetilde{\Omega}, \widetilde{\mathscr{H}},\widetilde{\mathbb E})$. By  Corollary 1 of Peng and Zhou (2020) we have
 \begin{align*} \widetilde{\Capc}(\xi>x)=& \frac{2}{\underline{\sigma}+\overline{\sigma}}\int_x^{\infty}\left( \phi(z/\overline{\sigma})I\{z\ge 0\}+\phi(z/\underline{\sigma})I\{z<0\}\right)dz\\
 =& \begin{cases}
 \frac{2\overline{\sigma}}{\underline{\sigma}+\overline{\sigma}}\left[1-\Phi\Big(\frac{x}{\overline{\sigma}}\Big)\right], & x\ge 0, \\
 1-\frac{2\underline{\sigma}}{\underline{\sigma}+\overline{\sigma}}\Phi\Big(\frac{x}{\underline{\sigma}}\Big), &x\le 0,
 \end{cases}
 \end{align*}
 where $\Phi(x)$ and $\phi(x)$, respectively, are the distribution function and density of a standard normal random variables in a probability space.
 Hence, by the fact that $-\xi\overset{d}=\xi$,
 \begin{align}\label{eqnormcapc} \widetilde{\cCapc}(\xi\ge x)=&\widetilde{\cCapc}(\xi\le -x) =1-\widetilde{\Capc}(\xi>-x)\nonumber \\
 =& \begin{cases}
 \frac{2\underline{\sigma}}{\underline{\sigma}+\overline{\sigma}}\left[1-\Phi\Big(\frac{x}{\underline{\sigma}}\Big)\right], & x\ge 0, \\
 1-\frac{2\overline{\sigma}}{\underline{\sigma}+\overline{\sigma}}\Phi\Big(\frac{x}{\overline{\sigma}}\Big), &x\le 0.
 \end{cases}
  \end{align}

From \eqref{eqnormcapc} and the central limit theorem, we can derive   a lower bound of an    exponential inequality  under $\cCapc$ for independent  and identically distributed random variables.
\begin{lemma}\label{lemLower} Suppose that $\{X_{ni}; i=1,\ldots,k_n\}$ is an array  of   independent  and identically distributed random variables in
the sub-linear expectation space $(\Omega, \mathscr{H}, \Sbep)$  with
$$\Sbep[X_{n1}^2]\to \overline{\sigma}^2<\infty,\;\; \cSbep[X_{n1}^2]\to \underline{\sigma}^2>0. $$
Let $\{x_n\}$ be a sequence of positive numbers such that
   $x_n\to \infty$, $x_n/\sqrt{k_n}\to 0$.
Assume
$$ \Sbep[( X_{n1}^2 -\epsilon k_n/x_n^2)^+]\to 0 \text{ for all } \epsilon>0, $$
and
$$ \frac{\sum_{i=1}^{k_n} (|\Sbep[X_{ni}]|+|\cSbep[X_{ni}]|)}{x_n\sqrt{k_n}}=\frac{\sqrt{k_n} (|\Sbep[X_{n1}]|+|\cSbep[X_{n1}]|)}{x_n} \to 0. $$
Denote $S_n=\sum_{i=1}^{k_n} X_{ni}$.
  Then for any $z>0$,
\begin{equation}\label{eqLowerB}
 \liminf_{n\to \infty} x_n^{-2}\ln \cCapc\left(  S_n \ge z \underline{\sigma} x_n \sqrt{k_n}   \right)\ge  -\frac{z^2}{2}.
   \end{equation}
\end{lemma}
  {\bf Proof}. Denote $S_{n,0}=0$, $S_{n,k}=\sum_{i=1}^k X_{ni}$.
  For $t>2$, let
\begin{eqnarray*}
 N=[k_n t^2/x_n^2],\;\; m=[x_n^2/t^2];\;\; r=\sqrt{k_n}x_n/(tm).
\end{eqnarray*}
 Then $mN\le k_n$, $r\sim\sqrt{N}$ and
\begin{align*}
& \left\{  \frac{S_n}{\underline{\sigma} x_n\sqrt{k_n}} \ge z \right\}
\supset   \left\{\frac{S_{n, Nm}}{\underline{\sigma} x_n\sqrt{k_n}} \ge z+  \epsilon/2\right\}
\bigcap \left\{ \left|\frac{S_n-S_{n,Nm}}{\underline{\sigma}x_n\sqrt{k_n} }\right|\le \epsilon/2\right\}\\
& \;\;=    \left\{  \frac{S_{n,Nm}}{r\underline{\sigma} }\ge tm(z+ \epsilon/2)\right\}
\bigcap \left\{ \left|\frac{S_n-S_{n, Nm}}{\underline{\sigma}x_n\sqrt{k_n} }\right|\le \epsilon/2\right\}\\
&\;\; \supset   \bigcap_{i=1}^m  \left\{   \frac{S_{n,Ni}-S_{n,N(i-1)}}{tr\underline{\sigma} }\ge  z + \epsilon  /2\right\}
\bigcap \left\{ \left|\frac{S_n-S_{n,Nm}}{\underline{\sigma}x_n\sqrt{k_n} }\right|\le \epsilon/2\right\}.
\end{align*}
For given $z>0$ and $\epsilon>0$. Let $f,g\in C_{b,Lip}(\mathbb R)$ such that $I\{x\ge z+\epsilon/2\}\ge f(x)\ge I\{x\ge z+\epsilon\}$ and $I\{|x|\le \epsilon/2\}\ge g(x)\ge I\{|x|\le \epsilon/4\}$.
It follows that
$$ I\left\{  \frac{S_n}{\underline{\sigma} x_n\sqrt{k_n}} \ge z \right\}
\ge \prod_{i=1}^m f\left(\frac{S_{n,Ni}-S_{n,N(i-1)}}{tr\underline{\sigma} }\right)
g\left(\frac{S_n-S_{n,Nm}}{\underline{\sigma}x_n\sqrt{k_n} }\right). $$
Note that
 $\{S_{n,Ni}-S_{n,N(i-1)},i=1,\ldots,m, S_n-S_{n,Nm}\}$ are independent under $\Sbep$ (and $\cSbep$). By \eqref{eq1.3},  we have
$$\cCapc\left( \frac{S_n}{\underline{\sigma} x_n\sqrt{k_n}} \ge z\right)
\ge \left(\cSbep\left[ f\left(\frac{S_{n,N}}{tr\underline{\sigma} }\right)\right]\right)^m
 \cSbep\left[g\left(\frac{S_n-S_{n,Nm}}{\underline{\sigma}x_n\sqrt{k_n} }\right)\right].
$$
Note
$$ \frac{\sum_{i=1}^N\Sbep[(X_{ni}^2-\epsilon N)^+]}{N}\to 0, \;\; \frac{\sum_{i=1}^N(|\Sbep[X_{ni}]|+|\cSbep[X_{ni}]|)}{\sqrt{N}}\to 0. $$
By applying the Lindeberg limit theorem of Zhang (2021),  we have
$$\lim_{n\to\infty}\Sbep\left[\varphi\Big(\frac{S_{n,N}}{r}\Big)\right]=\lim_{n\to\infty}\Sbep\left[\varphi\Big(\frac{S_{n,N}}{\sqrt{N}}\Big)\right]= \widetilde{\mathbb E}\left[\varphi (\xi)\right], \;\; \text{ for all } \varphi\in C_{b,Lip}(\mathbb R), $$
where $\xi\sim N(0,[\underline{\sigma}^2,\overline{\sigma}^2])$ under $\widetilde{\mathbb E}$.
It follows that
\begin{align*}
&\lim_{n\to \infty}\cSbep\left[ f\left(\frac{S_{n,N}}{tr\underline{\sigma} }\right)\right]
=\widetilde{\mathcal E}\left[ f\left(\frac{\xi}{t\underline{\sigma} }\right)\right]\\
\ge & \widetilde{\mathcal V}\Big(\xi\ge t(z+  \epsilon)\underline{\sigma}\Big)=\frac{2\underline{\sigma}}{\underline{\sigma}+\overline{\sigma}} \left[1-\Phi\Big( tz(1+\epsilon)  \Big)\right],
\end{align*}
by \eqref{eqnormcapc}.
On the other hand,
$$ 1-  \cSbep\left[g\left(\frac{S_n-S_{n,Nm}}{\underline{\sigma}x_n\sqrt{k_n} }\right)\right]\le
\Capc\left(\frac{|S_{n,k_n-Nm}|}{\underline{\sigma}x_n\sqrt{k_n} }\ge \epsilon/4\right)\le \frac{C\Sbep[X_{n1}^2]}{\epsilon^2}\frac{k_n-Nm}{\underline{\sigma}^2k_n x_n^2}\to 0.$$
 It follows that
\begin{align*}
& \liminf_{n\to \infty} x_n^{-2}\ln \cCapc\left( \frac{S_n}{\underline{\sigma} x_n\sqrt{k_n}} \ge z \right)\\
\ge &  \liminf_{t\to \infty}\liminf_{n\to \infty} t^{-2} m^{-1}\ln \cCapc\left( \frac{S_n}{\underline{\sigma} x_n\sqrt{k_n}} \ge z\right) \\
\ge & \liminf_{t\to \infty}\liminf_{n\to \infty} t^{-2}  \ln \cSbep\left[ f\left(\frac{S_{n,N}}{tr\underline{\sigma} }\right)\right]
=  \liminf_{t\to \infty} t^{-2}\ln \widetilde{\mathcal E}\left[ f\left(\frac{\xi}{t\underline{\sigma} }\right)\right]\\
\ge & \liminf_{t\to \infty} t^{-2}\ln \left[1-\Phi\Big(  tz(1+\epsilon)  \Big)\right]
=-\frac{(z(1+\epsilon))^2}{2}.
\end{align*}
 The proof   is completed. $\Box$

  \section{Properties of the Capacities}\label{sectCapc}
\setcounter{equation}{0}

Before we give the laws of the iterated logarithm,  we need  more notation and the properties of  capacities.

\begin{definition}\label{def3.1}
\begin{description}
\item{\rm (I)}   A function $V:\mathcal F\to [0,1]$ is called to be  countably sub-additive if
$$ V\Big(\bigcup_{n=1}^{\infty} A_n\Big)\le \sum_{n=1}^{\infty}V(A_n) \;\; \forall A_n\in \mathcal F. $$

\item{\rm (II)}  A capacity $V:\mathcal F\to [0,1]$ is called to be continuous  from below if it satisfies that $V(A_n)\uparrow V(A)$ whenever $A_n\uparrow A$, where $A_n, A\in \mathcal F$, and, it is called to be continuous from above  if it satisfies that $V(A_n)\downarrow  V(A)$ whenever $A_n\downarrow A$, where $A_n, A\in \mathcal F$.
\end{description}
\end{definition}
It is obvious that the continuity from above with the sub-additivity  implies the continuity from below, and the continuity from the below with the sub-additivity implies the countable sub-additivity. So, we call a sub-additive capacity to be  continuous   if it is continuous from above. Also, if $V$ is a  capacity continuous from above, then
\begin{equation}\label{CapcC} V(\bigcap_{i=1}^{\infty}A_i)=1 \text{ for events } \{A_n\} \text{ with } A_n\supset A_{n+1} \text{ and } V(A_n)=1,\; n=1,2,\cdots.
\end{equation}
It is obvious that the lower capacity $\cCapc$ has the property \eqref{CapcC} when the upper capacity $\Capc$ is countably sub-additive.

The following lemma is the Borel-Cantelli  Lemma and its converse  under capacities.

\begin{lemma} \label{lemBC}
\begin{description} \item[\rm (i)]
 Let $\{A_n, n\ge 1\}$ be a sequence of events in $\mathcal F$. Suppose that $V$ is a sub-additive capacity and $\sum\limits_{n=1}^{\infty}V\left (A_n\right)<\infty$.
 Then
 $$ \lim_{ n\to \infty}\max_N V\left(\bigcup_{i=n}^NA_i\right)=0. $$
If $V$ is a countably sub-additive capacity, then
\begin{equation}\label{eqlemBC.1} V\left (A_n\;\; i.o.\right)=0, \;\; \text{ where } \{A_n\;\; i.o.\}=\bigcap_{n=1}^{\infty}\bigcup_{i=n}^{\infty}A_i.
\end{equation}
\item[\rm (ii)] Suppose that $\{\xi_n; n\ge 1\}$ is a sequence of independent random variables in $(\Omega,\mathscr{H},\Sbep)$. Suppose $\sum\limits_{n=1}^{\infty}\Capc(\{\xi_n\ge 1+\epsilon\})=\infty$  for some $\epsilon>0$.
   Then
   \begin{equation}\label{eqlemBC.3}
    \Capc\left (\bigcup_{m=n}^{\infty}\{\xi_m\ge 1\}\right)  \ge \Capc\left (\bigcup_{m=n}^N\{\xi_m\ge 1\}\right)\to  1\; \text{ as }N\to \infty,
 \end{equation}
 and
$$ \Capc\left (\{\xi_n\ge 1\}\;\; i.o.\right)=1 \;\; \text{ if }  \;  \Capc  \; \text{ has the  property \eqref{CapcC} }.
$$
 \item[\rm (iii)] Suppose that  $\{\xi_n; n\ge 1\}$ is a sequence of independent random variables in $(\Omega,\mathscr{H},\Sbep)$, and $\Capc_1$ is  a countably sub-additive capacity with $\Capc_1\le \Capc$.
   Then
 \begin{equation}\label{eqlemBC.2} \cCapc_1\left (\{\xi_n\ge 1\}\;\; i.o.\right)=1\;\; \text{ if } \;  \sum_{n=1}^{\infty}\cCapc(\{\xi_n\ge 1+\epsilon\})=\infty \text{ for some }\epsilon>0.
 \end{equation}
  \item[\rm (iv)] Suppose  that $\{\xi_n; n\ge 1\}$ is a sequence of independent random variables in $(\Omega,\mathscr{H},\Sbep)$.
   Suppose $\sum\limits_{n=1}^{\infty}\cCapc(\{\xi_n\ge 1-\epsilon\})<\infty$  for some $\epsilon>0$. Then
  \begin{equation}\label{eqlemBC.4} \lim_{n\to\infty}\max_N\cCapc\left (\bigcup_{m=n}^N \{\xi_m\ge  1\} \right)= 0,
  \end{equation}
   and
  \begin{equation}\label{eqlemBC.5} \cCapc\left (\{\xi_n\ge 1\}\;\; i.o.\right)=0 \;\; \text{ if }  \;  \Capc  \; \text{ is continuous}.
 \end{equation}
 \end{description}
\end{lemma}
Lemma \ref{lemBC} (i) (resp. (iv)) is the direct part  of the Borel-Cantelli  Lemma for $\Capc$ (resp. $\cCapc$).  Parts (ii) or (iii) are the    converse ones.

{\bf Proof.}  (i) is trivial. For (ii), denote $A_n=\{\xi_n\ge 1\}$. Let $g(x)$ be a Lipschitz function with $I\{x\ge 1+\epsilon\}\le g(x)\le I\{x\ge 1\}$. Then
\begin{align}\label{eqprooflemBC1}
  & \cCapc\left(\bigcap_{i=n}^{\infty }A_i^c\right)\le \cCapc\left(\bigcap_{i=n}^{N}A_i^c\right)
\le     \cSbep\left[\prod_{i=n}^N\big(1-g(\xi_i)\big)\right] \nonumber\\
=  &   \prod_{i=n}^N \cSbep\left[\big(1-g(\xi_i)\big)\right]
=    \prod_{i=n}^N  \big(1-\Sbep[g(\xi_i)]\big)\nonumber\\
\le   & \exp\big\{-\sum_{i=n}^N\Sbep[g(\xi_i)]\big\}\le \exp\big\{-\sum_{i=n}^N\Capc(\xi_i\ge 1+\epsilon)\big\}\nonumber \\
  \to & \exp\big\{-\sum_{i=n}^{\infty}\Capc(\xi_i\ge 1+\epsilon)\big\}= 0 \; \text{ if } \;  \sum_{n=1}^{\infty}\Capc(\{\xi_n\ge 1+\epsilon\})=\infty.
\end{align}
That is $\Capc\left(\bigcup_{i=n}^{\infty }A_i\right)=1$ and $\Capc\left(\bigcup_{i=n}^NA_i\right)\to 1$ as $N\to \infty$.

 For (iii), similarly to \eqref{eqprooflemBC1} we have
$$\Capc\left(\bigcap_{i=n}^{\infty }A_i^c\right)\le \exp\big\{-\sum_{i=n}^{\infty}\cCapc(\xi_i\ge 1+\epsilon)\big\}=0  \; \text{ if } \;  \sum_{n=1}^{\infty}\cCapc(\{\xi_n\ge 1+\epsilon\})=\infty.
$$
 It follows from the countable sub-additivity of $\Capc_1$  that
\begin{align*}
\Capc_1\left(\big\{A_n \;\; i.o.\big\}^c\right) \le \sum_{n=1}^{\infty}\Capc_1\left(\bigcap_{i=n}^{\infty }A_i^c\right)\le  \sum_{n=1}^{\infty}\Capc\left(\bigcap_{i=n}^{\infty }A_i^c\right)=0.
\end{align*}
 Therefore,  $\cCapc_1\left(A_n \;\; i.o.\right)=1$.

 For (iv), we let $g(x)$ be a Lipschitz function with $I\{x\ge 1\}\le g(x)\le I\{x\ge 1-\epsilon\}$. Suppose
$\sum_{n=1}^{\infty}\cCapc(\xi_i\ge 1-\epsilon)<\infty$. Then
\begin{align*}
    &\Capc\left(\bigcap_{i=n}^{N}A_i^c\right)
\ge      \Sbep\left[\prod_{i=n}^N\big(1-g(\xi_i)\big)\right]
=      \prod_{i=n}^N \Sbep\left[\big(1-g(\xi_i)\big)\right]
=     \prod_{i=n}^N  \big(1-\cSbep[g(\xi_i)]\big)\\
\ge &  \prod_{i=n}^N \left(1- \cCapc(\xi_i\ge 1-\epsilon)\right)
\ge \exp\left\{-2\sum_{i=n}^N  \cCapc(\xi_i\ge 1-\epsilon)\right\}\ge \exp\left\{-2\sum_{i=n}^{\infty}  \cCapc(\xi_i\ge 1-\epsilon)\right\}
\end{align*}
for $N,n$ large enough. The last inequality is computed from the fact that $1-x\ge e^{-2x}$ for all $x\le 1/2$. Thus, \eqref{eqlemBC.4} is proved.
If  $\Capc$ is continuous, then
\begin{align*}
   \Capc\left(\bigcup_{n=1}^{\infty}\bigcap_{i=n}^{\infty}A_i^c\right)\ge \lim_{n\to \infty} \lim_{N\to \infty} \Capc\left(\bigcap_{i=n}^{N}A_i^c\right)
\ge  &    \lim_{n\to \infty} \exp\left\{-2\sum_{i=n}^{\infty}  \cCapc(\xi_i\ge 1-\epsilon)\right\}=1.
\end{align*}
Therefore, the proof is completed. $\Box$.

\bigskip
When the converse part of the Borel-Cantelli  lemma is applied, it is usually needed to suppose the continuity of the capacity $\Capc$. However, the following proposition shows that the capacities $\Capc$ and $\cCapc$ are usually not continuous.
\begin{proposition}\label{prop1} Let $(\Omega,\mathscr{H},\Sbep)$ be a sub-linear expectation  space with a sequence of independent and identically distributed  random variables $\{X_n; n\ge 1\}$. Consider the subspace
 \begin{equation}\label{eqprop1.1}\widetilde{\mathscr{H}}=\left\{Y=\varphi(X_1,X_2,\ldots,X_n): \varphi\in C_{l,Lip}(\mathbb R^n), \Sbep[(|Y|-c)^+]\to 0, n\ge 1\right\}.
 \end{equation}
If $\Capc$ is continuous on $\sigma(X_1,X_2,\ldots)$, then $\Sbep$ is linear on $\widetilde{\mathscr{H}}$.
\end{proposition}

{\bf Proof.} It is sufficient to show that
\begin{equation} \label{eqprop1}\Sbep[Y]=-\Sbep[-Y]\;\; \text{ for all } Y \in \widetilde{\mathscr{H}}.
\end{equation}
Without loss of generality, assume $Y=\varphi(X_1)$ and $|Y|\le c$. Denote $Y_n=\varphi(X_n)$. Then $\{Y_n;n\ge 1\}$ is a sequence of independent and identically distributed random variables with $|Y_n|\le c$. By \eqref{eqExpIneq.3},
$$ \cCapc\left(\frac{S_m}{m}\le \Sbep[Y]-\epsilon\right)=\cCapc\left(\sum_{i=1}^m(-Y_i+\Sbep[Y_i])\ge \epsilon m\right)
\le \exp\left\{-\frac{\epsilon^2 m^2}{2(\epsilon m c+ c^2m)}\right\}\to 0. $$
Hence
\begin{equation} \label{eqprop2} \Capc\left(\frac{S_m}{m}> \Sbep[Y]-\epsilon\right)\to 1\;\; \text{ for all }\epsilon>0.
\end{equation}
On the other hand, let $I(k)=\{2^k+1,\ldots, 2^{k+1}\}$. By \eqref{eqExpIneq.3}, for any $0<\epsilon<c/2$,
$$ \cCapc\left( \max_{n \in I(k) } \sum_{j\in I(k), j\le n}(Y_j-\cSbep[Y_j])\ge 2^{k+1}\epsilon\right)
\le \exp\left\{-\frac{\epsilon^2 2^{2(k+1)}}{2(\epsilon 2^{k+1} c+ c^22^k)}\right\}\le \exp\left\{-\frac{\epsilon^2}{c^2}2^k\right\}.
$$
Let $T_n=\sum_{j=1}^n(Y_j-\cSbep[Y_j])$.
Note the independence of the random variables. By Lemma \ref{lemBC} (iv), it follows that
$$ \lim_{l\to \infty}\max_L\cCapc\left( \bigcup_{k=l}^L \left\{\max_{n \in I(k) } \frac{T_n-T_{2^k}}{2^{k+1}}\ge  \epsilon \right\}\right)=0  \;\text{ for all } \epsilon>0, $$
 which implies
 $$ \lim_{n\to \infty}\max_N\cCapc\left( \max_{n\le l\le N} \frac{T_l}{l}\ge \epsilon\right)=0  \;\text{ for all } \epsilon>0 $$
(cf. the proof of \eqref{limsup1}).  That is
\begin{equation} \label{eqprop3} \lim_{ n\to \infty}\min_N\Capc\left( \max_{n\le l\le N} \frac{S_l}{l}\le \cSbep[Y]+\epsilon\right)=1  \;\text{ for all } \epsilon>0.
\end{equation}
 Let $f$ and $g$ be two Lipschitz functions with $I\{x\le \epsilon\}\ge f(x)\ge I\{x\le \epsilon/2\}$ and $I\{x\ge -\epsilon\}\ge g(x)\ge I\{x\ge -\epsilon/2\}$.
  By    the independence of the random variables,  it follows from \eqref{eqprop3} that
 \begin{align}\label{eqprop4}
   & \lim_{ n\to \infty}\min_N \Sbep\left[f\left( \max_{n\le l\le N}\frac{S_l-S_m}{l}-\cSbep[Y]\right)\cdot g\left(\frac{S_m}{m}- \Sbep[Y]\right)\right]\nonumber\\
  \ge & \lim_{ n\to \infty}\min_N\Sbep\left[f\left( \max_{n\le l\le N}\frac{S_l-S_m}{l}-\cSbep[Y]\right)\right]\cdot \Sbep\left[g\left(\frac{S_m}{m}- \Sbep[Y]\right)\right]\nonumber\\
  \ge & \lim_{ n\to \infty}\min_N\Capc\left( \max_{n\le l\le N}\frac{S_l-S_m}{l}\le  \cSbep[Y]+\epsilon/2\right) \cdot\Capc\left(  \frac{S_m}{m}\ge  \Sbep[Y]-\epsilon/2\right) \nonumber\\
  \ge & \lim_{ n\to \infty}\min_N\Capc\left( \max_{n\le l\le N}\frac{S_l}{l}\le  \cSbep[Y]+\epsilon/3\right) \cdot\Capc\left(  \frac{S_m}{m}\ge  \Sbep[Y]-\epsilon/2\right) \nonumber\\
  =&\Capc\left( \frac{S_m}{m}\ge  \Sbep[Y]-\epsilon/2\right).
 \end{align}
 Previously, we have not used the continuity or the property \eqref{CapcC} of $\Capc$. Now, notice that
 \begin{align*}
  \left\{\limsup_{n\to \infty} \frac{S_n}{n}\le \cSbep[Y]+2\epsilon\right\}\supset   \bigcup_{n=1}^{\infty} \bigcap_{N=n}^{\infty}\left\{\max_{n\le l\le N}\frac{S_l-S_m}{l}\le \cSbep[Y]+\epsilon\right\}.
\end{align*}
 By the continuity   of $\Capc$,  it follows from \eqref{eqprop4} that
 \begin{align}\label{eqprop5}
 & \Capc\left(\limsup_{n\to \infty} \frac{S_n}{n}\le \cSbep[Y]+2\epsilon \text{ and }  \frac{S_m}{m}\ge  \Sbep[Y]-\epsilon\right)\nonumber\\
 \ge & \lim_{n\to \infty}\lim_{N\to \infty}\Capc\left( \max_{n\le l\le N}\frac{S_n-S_m}{l}\le \cSbep[Y]+\epsilon \text{ and }  \frac{S_m}{m}\ge  \Sbep[Y]-\epsilon\right)\nonumber \\
 \ge &\lim_{ n\to \infty}\min_N  \Sbep\left[f\left( \max_{n\le l\le N}\frac{S_l-S_m}{l}-\cSbep[Y]\right)\cdot g\left(\frac{S_m}{m}- \Sbep[Y]\right)\right]\nonumber\\
 \ge &  \Capc\left( \frac{S_m}{m}\ge  \Sbep[Y]-\epsilon/2\right).
 \end{align}
 By letting $m\to \infty$, it follows from \eqref{eqprop2} that
 \begin{align*}
 & \Capc\left(\limsup_{n\to \infty} \frac{S_n}{n}\le \cSbep[Y]+2\epsilon \text{ and }  \limsup_{m\to\infty} \frac{S_m}{m}\ge  \Sbep[Y]-\epsilon\right)\\
 \ge & \limsup_{m\to \infty}\Capc\left(\limsup_{n\to \infty} \frac{S_n}{n}\le \cSbep[Y]+\epsilon \text{ and }  \frac{S_m}{m}\ge  \Sbep[Y]-\epsilon\right)=1.
 \end{align*}
 Therefore, $\Sbep[Y]-\epsilon<\cSbep[Y]+2\epsilon$ for every $\epsilon>0$. Hence, \eqref{eqprop1} is verified and the proof is completed. $\Box$.

 Let $B(t)$ be a $G$-Brownian motion. Denote $X_n=\sqrt{n(n+1)}\Big(B(1- 1/(n+1))-B(1-1/n)$. Then $X_1,X_2,\ldots$ is a sequence of independent and identically distributed $G$-normal random variables. Applying Proposition \ref{prop1}, we have the following corollary.
\begin{corollary} The $G$-capacity $\hat{c}$ as defined in section 6.3 of Peng (2019) is not continuous unless $B(t)$ is a classical Brownian motion in a probability space.
\end{corollary}

 \bigskip
  According to Proposition  \ref{prop1}, the continuity of a sub-additive capacity is a very stringent condition. It is needed to avoid assuming the continuity of a   capacity. Because the Borel-Cantelli lemma (Lemma \ref{lemBC} (i)) is needed when   the strong limit theorems, e.g.,  the LIL, are considered, we usually assume that the capacity $\Capc$ is countably sub-additive. Such a condition  is      satisfied when $\Sbep$ can be presented in the form of \eqref{eq1.7} (cf. Chapters 3 and 6 of Peng (2019)).
But the capacity $\upCapc$ defined as in \eqref{eq1.5}  may be not countably sub-additive   so that even the direct part of the Borel-Cantelli lemma is not valid.  So we  consider its countably sub-additive extension.
\begin{definition}   A  countably sub-additive extension $\outCapc$  of $\upCapc$   is defined by
\begin{equation}\label{outcapc} \outCapc(A)=\inf\Big\{\sum_{n=1}^{\infty}\upCapc(A_n): A\subset \bigcup_{n=1}^{\infty}A_n\Big\},\;\; \outcCapc(A)=1-\outCapc(A^c),\;\;\; A\in\mathcal F,
\end{equation}
where $\upCapc$ is defined as in \eqref{eq1.5}.
\end{definition}

As shown in Zhang (2016), $\outCapc$ is countably sub-additive, $\outCapc(A)\le \upCapc(A)$ and  $\outCapc=\upCapc$ when $\upCapc$ is countably sub-additive. So, (i) and (iii) of Lemma \ref{lemBC}  are valid for $\outCapc$. It is shown by Zhang (2016) that, if $V$ is also a sub-additive (resp. countably sub-additive) capacity satisfying
\begin{equation}\label{eqVprop} V(A)\le \Sbep[g]  \text{ whenever }  I_A\le g\in \mathscr{H},
\end{equation}
 then $V(A)\le \upCapc$ (resp. $V(A)\le \outCapc(A)$.  Hence, if there exists a countably sub-additive capacity having   the property \eqref{eq1.3}, then $\outCapc$ has the property \eqref{eq1.3}.

 \begin{definition} Another countably sub-additive capacity generated by $\Sbep$ can be defined as follows:
 \begin{equation}\label{tildecapc} \tildeCapc(A)=\inf\Big\{\lim_{n\to\infty}\Sbep[\sum_{i=1}^n g_i]: I_A\le \sum_{n=1}^{\infty}g_n, 0\le g_n\in\mathscr{H}\Big\},\;\;\; A\in\mathcal F.
\end{equation}
\end{definition}

We can show that $\tildeCapc$ is a countably sub-additive capacity having the property \eqref{eqVprop}, and so, $\tildeCapc(A)\le \outCapc(A)$. Further,
 if $\Sbep$ has the form \eqref{eq1.7}, then
$$ \Capc^{ \mathscr{P}}(A)=\sup_{P\in \mathscr{P}}P(A)\le \tildeCapc(A)\le \outCapc(A), \;\; A\in\mathcal{F},  $$
by noting that
$$ P(A)\le P[\sum_{n=1}^{\infty}g_n]=\lim_{n\to \infty}P[\sum_{i=1}^ng_i]\le \lim_{n\to \infty}\Sbep[\sum_{i=1}^ng_i]$$
when $I_A\le \sum_{n=1}^{\infty}g_n$ and $g_n\ge 0$.

The out capacity $c^{\prime}$ defined in Example 6.5.1 of Peng (2019) coincides with  $\tildeCapc$ if $\mathscr{H}$ is chosen as the family of (bounded) continuous functions on a metric space $\Omega$.

Since  $\outCapc$ and $\mathbb C^{\ast}$ are countably sub-additive capacities, the direct part of    the  Borel-Cantelli lemma (Lemma \ref{lemBC} (i)) is valid for them. But they may   not be continuous unless $\Sbep$ is linear. To make the converse part of the  Borel-Cantelli lemma (Lemma \ref{lemBC} (ii)) valid, it would be  reasonable to  assume  \eqref{CapcC} instead of the continuity of $\Capc$.   Unfortunately, the following proposition tells us that \eqref{CapcC} is also a stringent condition.

\begin{proposition}\label{prop2} Let $(\Omega,\mathscr{H},\Sbep)$ be a sub-linear expectation  space with a sequence of independent and identically distributed  random variables $\{X_n; n\ge 1\}$. Consider the subspace $\widetilde{\mathscr{H}}$ defined as in \eqref{eqprop1.1}.
Suppose that  the following condition is satisfied.
\begin{description}
      \item[\rm (CC)]      The sub-linear expectation $\Sbep$  on $\mathscr{H}_b$  satisfies
\begin{equation} \label{eqexpressbyP} \Sbep[X]=\sup_{P\in \mathscr{P}}P[X], \; X\in \mathscr{H}_b
\end{equation}
where $\mathscr{H}_b=\{f\in \mathscr{H}; f \text{ is bounded}\}$, $\mathscr{P}$ is a  countable-dimensionally weakly compact family of probability measures on $(\Omega,\sigma(\mathscr{H}))$ in sense that, for any $Y_1,Y_2,\ldots \in \mathscr{H}_b$ and any  sequence $\{P_n\}\subset \mathscr{P}$ there are a subsequence $\{n_k\}$ and a probability measure $P\in \mathscr{P}$ for which
\begin{equation}\label{eqcompact1} \lim_{k\to \infty} P_{n_k}[\varphi(Y_1,\ldots,Y_d)]= P[\varphi(Y_1,\ldots,Y_d)],\; \varphi\in C_{b,Lip}(\mathbb R^d), d\ge 1.
\end{equation}
\end{description}
Define
\begin{equation} \label{eqdefC-P} \Capc^{\mathscr{P}}(A)=\sup_{P\in \mathscr{P}}P(A),\;\; A\in \sigma(\mathscr{H}).
\end{equation}
  Then, for  $\Capc=\Capc^{\mathscr{P}}$,  $\mathbb C^{\ast}$,  $\outCapc$ or $\upCapc$ we have
   \begin{description}
       \item[\rm (i)] $\Capc$ has the property \eqref{eq1.3};
       \item[\rm (ii)] If on $\sigma(X_1,X_2,\ldots)$,  $\Capc$ has the property that
  \begin{equation}\label{CapcC2} \Capc(\bigcap_{i=1}^{\infty}A_i)>0 \text{ for events } \{A_n\} \text{ with } A_n\supset A_{n+1} \text{ and } \Capc(A_n)=1,\; n=1,2,\cdots,
\end{equation}
then $\Sbep$ is linear on $\widetilde{\mathscr{H}}$.
     \end{description}
 \end{proposition}

Before proving the proposition, we first give examples for which the condition (CC) is satisfied.
\begin{lemma}\label{lem4.2} If one of the following conditions is satisfied, then the condition (CC) is satisfied.
\begin{description}
 \item[\rm (a)]  $\Omega$ is a complete separable metric space, each element $X(\omega)$ in $\mathscr{H}$ is a continuous function on $\Omega$. The sub-linear expectation $\Sbep$ satisfies
$$ \Sbep[X]= \max_{P\in \mathscr{P}}P[X], \;X\in \mathscr{H}_b,$$
where $\mathscr{P}$ is a weakly compact family of probability measures on the metric space $\Omega$.
 \item[\rm (b)]     $\Omega$ is a complete separable metric space, each element $X(\omega)$ in $\mathscr{H}$ is a continuous function on $\Omega$.  There is a capacity $V$ with the property \eqref{eq1.3}  and is tight in sense that for any $\epsilon>0$, there is a compact set $K\subset \Omega$ such that $V(K^c)<\epsilon$.  Let $\mathscr{P}$  be the family of  all probability measures $P$ on $\sigma(\mathscr{H})$ which satisfies $P[f]\le \Sbep[f]$ for all $f\in \mathscr{H}_b$.
 \item[\rm (c)]  $\Sbep$ on $\mathscr{H}_b$ is regular in the sense that $\Sbep[X_n]\downarrow 0$  for any   elements  $\mathscr{H}_b\ni X_n\downarrow0$.  Let $\mathscr{P}$  be the family of  all probability measures $P$ on $\sigma(\mathscr{H})$ which satisfies $P[f]\le \Sbep[f]$ for all $f\in \mathscr{H}_b$.
 \item[\rm (d)] Let $\mathcal T$ be an index set, $\Omega=\mathbb R^{\bigotimes \mathcal T}=\{\bm x=(x_t;t\in \mathcal T);x_t\in \mathbb R, t\in
 \mathcal{T}\}$ be the product of real spaces. Consider the function space on $\Omega$ as
      $$\mathscr{H}=\{\varphi\circ\pi_{t_1,\ldots,t_d}: \varphi\in C_{l,Lip}(\mathbb R^d), t_1,\ldots,t_d\in \mathcal T, d\ge 1\}, $$
where $\pi_{t_1,\ldots,t_d}$ is a project map, $\pi_{t_1,\ldots,t_d}\bm x=(x_{t_1},\ldots, x_{t_d})$. Let $\Sbep$ be a sub-linear expectation on $\mathscr{H}$
with $\upCapc(|\pi_t\bm x|\ge c)\to 0$ as $c\to \infty$ for all $t\in\mathcal T$, and
  $\mathscr{P}$  be the family of  all probability measures $P$ on $\sigma(\mathscr{H})$ which satisfies $P[f]\le \Sbep[f]$ for all $f\in \mathscr{H}_b$.

  Further, in the definition of $\mathscr{H}$, $C_{l,Lip}(\mathbb R^d)$ can be replaced by any a class $\mathscr{C}(\mathbb R^d)$ of continuous functions on $\mathbb R^d$.
\end{description}
\end{lemma}

{\bf Proof of Proposition \ref{prop2}.}  Notice  \eqref{eqexpressbyP} and \eqref{eqdefC-P}. It is obvious that   $\Capc^{\mathscr{P}}(A)\le \mathbb C^{\ast}(A) \le \outCapc(A)\le \upCapc(A)$. (i) is obvious since $\Capc^{\mathscr{P}}$ and  $\upCapc$ satisfy \eqref{eq1.3}.

For proving (ii),    we first show that,  if $\bm X=(X_1,X_2,\ldots)$ is a sequence of random variables in $\mathscr{H}$ for which each $X_n$ is tight in the sense  that $\Capc(|X_n|\ge c)\to 0$ as $c\to\infty$, then
\begin{equation}\label{eqClose}
\begin{aligned}
 & \Capc^{\mathscr{P}}\left(\bm X\in F_n\right)\searrow \Capc^{\mathscr{P}}\left(\bm X\in F\right)  \\
  \text{whenever } & F_n \text{ are   closed subsets of  }\; \mathbb R^{\infty}\; \text{ with } F_n\searrow F.
  \end{aligned}
  \end{equation}
 Notice that the condition (CC) is satisfied. Consider the family of probability measures $\mathscr{P}$ on $\sigma(\bm X)$. For each $\epsilon>0$, by the tightness of $X_i$, there exists a positive constant $C_i$ such that
$$ \sup_{P\in \mathscr{P}}P\left(|X_i|>C_i\right)\le   \Capc(|X_i|\ge C_i/2)<\epsilon/2^i.$$
 Let $K=\bigotimes_{i=1}^{\infty}[-C_i,C_i]$. Then $K$ is a compact set on the metric space $\mathbb R^{\infty}$, and
$$ \sup_{P\in \mathscr{P}}P\left(\bm X \notin K\right)\le  \sum_{i=1}^{\infty}\sup_{P\in \mathscr{P}}P\left(|X_i|>C_i\right)<\epsilon. $$
Hence, $\mathscr{P}\bm X^{-1}=\{\overline{P}: \overline{P}(A)=P( \bm Y\in A), A\in \mathscr{B}(\mathbb R^{\infty}),P\in \mathscr{P}\}$ is tight and so a relatively weakly compact family of probability measures on the metric space $\mathbb R^{\infty}$ by Prohorov's theorem (cf. Billingsley (1999, Page 58)). Next, we show that $\mathscr{P}\bm X^{-1}$ is closed. Suppose that $\{P_n\bm X^{-1};P_n\in \mathscr{P}\}$ is a weakly convergent sequence on $\mathbb R^{\infty}$.  Then there  exists a probability measure $Q$ on $\mathbb R^{\infty}$ such that
$$ Q[f]=\lim_{n\to \infty}P_n[f(\bm X)], f\in C_b(\mathbb R^{\infty}).
$$
 We must   show that the limit $Q$ is determined by a probability measure $P\in \mathscr{P}$.
  Consider the sequence $\bm Y=\{X_n^{(l)}; n=1,2,\ldots, l=1,2,\ldots\}$ in $\mathscr{H}_b$.
 By the conditions assumed, for the sequence $\{P_n\}$ there exists a subsequence $\{n_k\}$ and a probability measure $P\in \mathscr{P}$ such that \eqref{eqcompact1} holds.  Hence
 $$ Q[f(x_1^{(l)},\ldots,x_d^{(l)})]=P[f(X_1^{(l)},\ldots,X_d^{(l)})], \;\; \forall f\in C_{b,lip}(\mathbb R^d),\; d\ge 1, l\ge 1, $$
 which, by letting $l\to \infty$ and noting the continuity of $f$, $Q$ and $P$, implies
 $$ Q[f(x_1,\ldots,x_d)]=P[f(X_1,\ldots,X_d)], \;\; \forall f\in C_{b,lip}(\mathbb R^d),\; d\ge 1. $$
 Notice that $Q$ and $P\bm X^{-1}$ are both determined by their common finite-dimensional distributions. It follows that
 $$ P[f(\bm X)]=Q[f(\bm x)]=\lim_{n\to\infty}P_n[f(\bm X)], f\in C_b(\mathbb R^{\infty}). $$
We conclude that $\mathscr{P}\bm X^{-1}$ is closed and so weakly compact. If let
$$\widetilde{V}(A)=\Capc^{\mathscr{P}}(\bm X\in A)=\sup_{P\in \mathscr{P}}P(\bm X\in A), \; A\in \mathscr{B}(\mathbb R^{\infty}), $$
then by Lemma 6.1.12 of Peng (2019), for any closed sets $F_n$s in the metric space $\mathbb R^{\infty}$ with $F_n\searrow F$ we have $\widetilde{V}(F_n)\searrow\widetilde{V}(F)$. \eqref{eqClose} is proved.

Now, we prove (ii). Let $Y,Y_1,Y_2,\ldots,$ be independent and identically distributed bounded random variables being defined as in the proof of Proposition \ref{prop1}. It is sufficient to show that $\Sbep[Y]=\cSbep[Y]$.   Write $\bm Y=(Y_1,Y_2,\ldots)$.   With the same arguments as in the proof of Proposition \ref{prop1}, both \eqref{eqprop2}  and \eqref{eqprop4}  also remain   true. Write
 \begin{equation}\label{eqkeyevent} A_l=\left\{\limsup_{n\to \infty} \frac{S_n}{n}\le \cSbep[Y]+2\epsilon \text{ and } \bigcup_{m=l}^{\infty} \left\{\frac{S_m}{m}\ge  \Sbep[Y]-\epsilon\right\}\right\}.
\end{equation}
We will show that $\Capc^{\mathscr{P}}(A_l)=1$ and so $\Capc(A_l)=1$.

 Note, $\big\{ \max_{n\le l\le N}\frac{S_l-S_m}{l}\le \cSbep[Y]+\epsilon \text{ and }  \frac{S_m}{m}\ge  \Sbep[Y]-\epsilon\big\}_{N=n}^{\infty}$ is a decreasing sequence of closed sets of $(Y_1,Y_2,\ldots)$. By \eqref{eqClose} we have
\begin{align*}
 &\Capc^{\mathscr{P}}\left(\max_{l\ge n}\frac{S_l-S_m}{l}\le \cSbep[Y]+\epsilon \text{ and }  \frac{S_m}{m}\ge  \Sbep[Y]-\epsilon\right)\\
 =& \lim_{N\to \infty} \Capc^{\mathscr{P}}\left(\max_{n\le l\le N}\frac{S_l-S_m}{l}\le \cSbep[Y]+\epsilon \text{ and }  \frac{S_m}{m}\ge  \Sbep[Y]-\epsilon\right)\\
 \ge &   \lim_{N\to \infty} \Sbep\left[f\left( \max_{n\le l\le N}\frac{S_l-S_m}{l}-\cSbep[Y]\right)\cdot g\left(\frac{S_m}{m}- \Sbep[Y]\right)\right],
\end{align*}
where the  inequality is due to \eqref{eq1.3}. Notice that
 \begin{align*}
  \left\{\limsup_{n\to \infty} \frac{S_n}{n}\le \cSbep[Y]+2\epsilon\right\}\supset   \bigcup_{n=1}^{\infty} \left\{\max_{l\ge n}\frac{S_l-S_m}{l}\le \cSbep[Y]+\epsilon\right\}.
\end{align*}
It follows from \eqref{eqprop4}  that
 \begin{align*}
 &\Capc^{\mathscr{P}}\left(\limsup_{n\to \infty} \frac{S_n}{n}\le \cSbep[Y]+2\epsilon \text{ and }  \frac{S_m}{m}\ge  \Sbep[Y]-\epsilon\right)\nonumber\\
 \ge &   \lim_{  n\to \infty}\Capc^{\mathscr{P}}\left(\max_{l\ge n}\frac{S_l-S_m}{l}\le \cSbep[Y]+\epsilon \text{ and }  \frac{S_m}{m}\ge  \Sbep[Y]-\epsilon\right) \nonumber\\
 \ge &   \lim_{  n\to \infty}\lim_{N\to \infty} \Sbep\left[f\left( \max_{n\le l\le N}\frac{S_l-S_m}{l}-\cSbep[Y]\right)\cdot g\left(\frac{S_m}{m}- \Sbep[Y]\right)\right]\nonumber \\
 \ge &   \lim_{ n\to \infty}\min_N  \Sbep\left[f\left( \max_{n\le l\le N}\frac{S_l-S_m}{l}-\cSbep[Y]\right)\cdot g\left(\frac{S_m}{m}- \Sbep[Y]\right)\right]\nonumber \\
 \ge &  \Capc\left( \frac{S_m}{m}\ge  \Sbep[Y]-\epsilon/2\right)\to 1 \text{ as } m\to \infty
 \end{align*}
by  \eqref{eqprop2}.
 Therefore,
$$
 \Capc^{\mathscr{P}}\left(\limsup_{n\to \infty} \frac{S_n}{n}\le \cSbep[Y]+2\epsilon \text{ and } \bigcup_{m=l}^{\infty} \left\{\frac{S_m}{m}\ge  \Sbep[Y]-\epsilon\right\}\right)=1 \; \text{ for all }  l\ge 1.
$$
Hence
$$\Capc(A_l)=1\; \text{ for all }l.
$$
Therefore, by the property \eqref{CapcC2}  we have
$$
 \Capc\left(\limsup_{n\to \infty} \frac{S_n}{n}\le \cSbep[Y]+2\epsilon \text{ and } \limsup_{m \to\infty}  \frac{S_m}{m}\ge  \Sbep[Y]-\epsilon \right)=
 \Capc\left(\bigcap_{l=1}^{\infty}A_l\right)>0.
$$
It follows that $\Sbep[Y]-\epsilon\le \cSbep[Y]+2\epsilon$ for all $\epsilon>0$. Therefore, $\Sbep[Y]=\cSbep[Y]$.
The proof is completed.
 $\Box$

{\bf Proof of Lemma \ref{lem4.2}}.  It is obvious that  the condition (a) implies the condition (CC), since $\varphi(Y_1,\ldots,Y_d)$ is a continuous function on $\Omega$.

 For the case that (b) is satisfied,    it is sufficient to show that $\Sbep$ is regular on $\mathscr{H}_b$  and so that the condition (c) is satisfied.
  Suppose that $\mathscr{H}_b\ni f_n\downarrow 0$ and $f_n\le L$. By the tightness of $V$, for any $\epsilon>0$, there is a compact $K$ such that $V(K^c)<\epsilon/L$. Notice that $f_n$ is continuous and so the sequence of functions $\{f_n\}$  is uniformly convergent  on the compact set $K$. Then
$$ \delta_n=: \sup_{\omega\in K} f_n(\omega)\downarrow 0 \text{ and } f_n\le \delta_n+LI_{K^c}. $$
It is obvious that
$$0\le \Sbep[f_n]\le \delta_n+LV(K^c) \le \delta_n+\epsilon$$
by \eqref{eq1.3}. That is,  $\Sbep[f_n]\downarrow 0$.  Hence $\Sbep$ is regular on $\mathscr{H}_b$.

 Suppose that the condition (c) is satisfied. Notice the expression \eqref{linearexpression}.
Consider the linear expectation $E_{\theta}$ on $\mathscr{H}_b$. If $\mathscr{H}_b\ni f_n\downarrow 0$, then $0\le E_{\theta}[f_n]\le \Sbep[f_n]\to 0$.
Hence, similar to Lemmas 1.3.5 and   6.2.2 of Peng (2019), by the Daniell-Stone  theorem, there is a unique probability $P_{\theta}$ on $\sigma(\mathscr{H}_b)=\sigma(\mathscr{H})$ such that
 $$ E_{\theta}[f]=   P_{\theta}[f]  \text{ for all bounded } f\in \mathscr{H}. $$
 Hence
 $$\Sbep[f]=\sup_{\theta\in \Theta}E_{\theta}[f]=  \sup_{\theta\in \Theta} P_{\theta}[f]  \text{ for all bounded } f\in \mathscr{H}. $$
Recall that $\mathscr{P}$ the family of all probability measures $P$ on $\sigma(\mathscr{H})$ which satisfies $P[f]\le \Sbep[f]$ for all $f\in \mathscr{H}_b$.
 Then
 $$\Sbep[X]=\sup_{\theta\in \Theta} P_{\theta}[X]\le \sup_{P\in \mathscr{P}}P[X]\le \Sbep[X],\;\;  X\in \mathscr{H}_b. $$
 \eqref{eqexpressbyP} holds.

Consider the family of probability measures $\mathscr{P}$ on $\sigma(\bm Y)$. Notice that each $Y_i$ is bounded and so is tight, which implies that
$\mathscr{P}\bm Y^{-1}$ is a relatively weakly compact family of probability measures on $\mathbb R^{\infty}$ as shown in the proof of Proposition \ref{prop2}. Next, we show that $\mathscr{P}\bm Y^{-1}$ is closed. Suppose that $\{P_n\bm Y^{-1};P_n\in \mathscr{P}\}$ is a weakly convergent sequence.  Let a   linear expectation $E$ be defined as
$$E[f(\bm Y)]=\lim_{n\to \infty}P_n[f(\bm Y)], f\in C_b(\mathbb R^{\infty}).$$
Then, $E$ is a linear expectation on the subspace $\mathscr{L}=\{f(Y_1,\ldots,Y_d):f\in C_{b,Lip}(\mathbb R^d), d\ge 1\}$ with $E\le \Sbep$. So,  by the Hahn-Banach theorem, there exists  a finite additive linear expectation  $E^e$ defined on $\mathscr{H}$ such that, $E^e=E$ on $\mathscr{L}$ and, $E^e\le \Sbep$ on $\mathscr{H}$. For $E^e$, as shown before, there is probability measure $P^e$ on $\sigma(\mathscr{H})$ such that $P^e[f]=E^e[f]$ for all $f\in \mathscr{H}_b\supset \mathscr{L}$. Hence $P^e\in \mathscr{P}$ and
$$ P^e[f(\bm Y)]=E[f(\bm Y)]=\lim_{n\to \infty}P_n[f(\bm Y)], f\in C_b(\mathbb R^{\infty}). $$
It follows that $\mathscr{P}\bm Y^{-1}$ is closed and so weakly compact.  So, the condition (CC) is satisfied.

Suppose that the condition (d) is satisfied.  It is sufficient to show that $\Sbep$  is regular on $\mathscr{H}_b$.
Suppose $\mathscr{H}_b\ni X_n\searrow 0$ and $0\le X_n\le L$. For each $X_n$, there are $t_{n,1},t_{n,2},\ldots,t_{n,d_n}\in \mathcal T$ such that
$X_n(\bm x)=\varphi_n(x_{t_{t,1}},\ldots,x_{t_{n, d_n}})$, $\varphi_n\in C_{l,Lip}(\mathbb R^{d_n})$. So, we can choose an index set $\mathcal{S}=:\{t_1,t_2,\ldots\}\subset \mathscr{T}$ such that $X_n(\bm x)=\psi_n(x_{t_1},\ldots, x_{t_{p_n}})$ with $\psi_n\in C_{l,Lip}(\mathbb R^{p_n})$. Hence $X_n(\bm x)$ is a continuous function of $\pi_{\mathcal S}\bm x=(x_{t_1},x_{t_2},\ldots)$ on $\mathbb R^{\infty}$ and can be written by $f_n(\pi_{\mathcal S}\bm x)$. By the condition assumed, for any $\epsilon>0$, there exists $C_i$ such that $\upCapc(|\pi_{t_i}\bm x|>C_i)<\epsilon/(2^iL)$. Let $K=\bigotimes_{i=1}^{\infty}[-C_i,C_i]$ and $K_n=\bigotimes_{i=1}^{p_n}[-C_i,C_i]$. Then $K$ is a compact set on $\mathbb R^{\infty}$. Hence
$$ \delta_n=:\sup_{\pi_{\mathcal S}\bm x\in K}|f_n(\pi_{\mathcal S}\bm x)|\to 0. $$
Notice that
\begin{align*}
& |X_n(\bm x)|\le   \delta_n+\sup_{\pi_{\mathcal S}\bm x\not \in K}|\psi_n(x_{t_1},\ldots,t_{p_n})| \\
=& \delta_n+\sup_{(x_{t_1},\ldots, x_{t_{p_n}})\not \in K^n}|\psi_n(x_{t_1},\ldots,t_{p_n})|
\le \delta_n+LI\{(x_{t_1},\ldots, x_{t_{p_n}})\not \in K^n\}.
\end{align*}
By the (finite) sub-additivity of $\upCapc$, it follows that
\begin{align*}
&\Sbep[|X_n|]\le \delta_n +\upCapc\left(\bm x:(x_{t_1},\ldots, x_{t_{p_n}})\not \in K^n\right) \\
\le  & \delta_n+L\sum_{i=1}^{p_n}\upCapc(|\pi_{t_i}\bm x|>C_i)
\le  \delta_n +\sum_{i=1}^{\infty}\epsilon/2^i\le \delta_n+\epsilon.
\end{align*}
Letting $n\to \infty$ and then $\epsilon\to 0$ yields $\Sbep[X_n]\to 0$.  Hence, $\Sbep$ is regular on $\mathscr{H}_b$.
The proof is completed. $\Box$

\begin{remark} \label{remarkCon} Actually, the conditions (CC), (c) and the following statement are  are equivalent:
 \begin{description}
   \item[\rm (e)]  there is a capacity $V$ with the property \eqref{eq1.3} such that any sequence $\{X_n; n\ge 1\}$ of tight random variables satisfies \eqref{eqClose}.
 \end{description}
 When $\Omega$ is a complete separable metric space and $\mathscr{H}=C(\Omega)$ or $C_b(\Omega)$, they are also equivalent to (a) and (b).
\end{remark}
 In fact, (c) $ \implies $ (CC)$\implies$ (e) is proved above. For (e) $ \implies $ (c), suppose $\mathscr{H}_b\ni X_n\searrow 0$ and consider $\bm X=(X_1,X_2,\ldots)$. Notice $\{X_n\ge \epsilon\}$ is a closed set of $\bm X$ and $\{X_n\ge \epsilon\}\searrow\emptyset$. By \eqref{eqClose},
 $$ 0\le\Sbep[X_n]\le \epsilon+V(X_n\ge \epsilon)\searrow \epsilon. $$
 Hence, $\Sbep[X_n]\to 0$. (c) holds. It is obvious that (a) implies (b) with $V=\Capc^{\mathscr{P}}$, and (a) or (b) $\implies$ (CC) is proved in the proof of Lemma \ref{lem4.2}.
 At last, suppose that $\Omega$ is a complete separable metric space,  $\mathscr{H}=C(\Omega)$ or $C_b(\Omega)$, and the condition (CC) is satisfied.
 We want to prove that $\mathscr{P}$ is a weakly compact family of probability measures on $\Omega$, and so the condition (a) is satisfied.  Since  $C(\Omega)$ is a separable topological space,  there exists a countable family $\{h_1,h_2,\ldots\}\subset C(\Omega)$ which is dense in $C(\Omega)$, i.e., for any $h\in C(\Omega)$ there exists a sequence $i_j$ such that $h(\omega)=\lim_{j\to \infty}h_{i_j}(\omega)$ for all $\omega\in \Omega$. Now, let $\{P_n\}\subset \mathscr{P}$. Applying the condition (CC) to $\{P_n\}$ and the sequence $\{h_i^{(l)};i,l=1,2,\ldots\}$ yields that there is a subsequence $\{P_{n_k}\}$ and   $P\in \mathscr{P}$ for which
\begin{equation}\label{eqremarkCon.1} P_{n_k}[\varphi(h_1^{(l)},\ldots,h_d^{(l)})]\to P[\varphi(h_1^{(l)},\ldots,h_d^{(l)})], \;\; \varphi\in C_{b,Lip}(\mathbb R^d), \; d,l\ge 1.
\end{equation}
Let $h\in C(\Omega)$ with $|h(\omega)|\le L$. For (a), it is sufficient to show that
\begin{equation}\label{eqremarkCon.2} P_{n_k}[h]\to P[h].
\end{equation}
It is sufficient to show that, for any subsequence of $\{n_k\}$ (without loss of generality, we assume that it is $\{n_k\}$ itself), there is a further subsequence $\{n_k^{\prime}\}\subset \{n_k\}$ such that
$$ P_{n_k^{\prime}}[h]\to P[h]. $$
Consider $\{P_{n_k}\}$ and random variables $\{h, h_i^{(l)};i,l=1,2,\ldots\}$. By the condition (CC) again, there exists a subsequence $\{n_k^{\prime}\}$ and  $Q\in \mathscr{P}$ such that
\begin{equation}\label{eqremarkCon.3} P_{n_k^{\prime}}[\varphi(h, h_1^{(l)},\ldots,h_d^{(l)})]\to Q[\varphi(h,h_1^{(l)},\ldots,h_d^{(l)})], \;\; \varphi\in C_{b,Lip}(\mathbb R^{d+1}), \; d,l\ge 1.
\end{equation}
Combining \eqref{eqremarkCon.1} and \eqref{eqremarkCon.3} yields
$$ Q[h_d^{(l)}]=P[h_d^{(l)}], \;\; d,l\ge 1. $$
Since there exists a sequence $i_j$ such that $h=\lim_{j\to \infty}h_{i_j}=\lim_{j\to \infty}h_{i_j}^{(2L)}$, we have  that
$$ Q[h]=\lim_{j\to \infty}Q[h_{i_j}^{(2L)}]=\lim_{j\to \infty}P[h_{i_j}^{(2L)}]=P[h], $$
which, together with \eqref{eqremarkCon.3}, implies that
$$ P_{n_k^{\prime}}[h]\to Q[h]=P[h]. $$
Hence, \eqref{eqremarkCon.2} holds.
\bigskip

For understanding the capacity $\outCapc$, we give the last lemma of this section which shows that it coincides with a classical probability measure $P$ if $\Sbep$ is introduced by $P$.

\begin{lemma}\label{lemlinear} Let $(\Omega, \mathscr{H},\Sbep)$ be a sub-linear space with
\begin{equation}\label{eqlemlinear.1}\Sbep[f]=P[f], \;\; f\in \mathscr{H},
\end{equation}
 where $P$ is a probability measure on $(\Omega, \mathcal{F})$. Then
\begin{equation}\label{eqlemlinear.2}
\outCapc(A)=\tildeCapc(A)=P(A), \;\; \forall \; A\in \sigma(\mathscr{H}).
\end{equation}
Further, if $\Capc$ is a countably sub-additive capacity with the property \eqref{eq1.3}, then $\Capc(A)=P(A)$ for all $A\in \sigma(\mathscr{H})$.
\end{lemma}

{\bf Proof. }  First, it is obvious that $P(A)\le \tildeCapc(A)$.

  Fix  $X_1,\ldots,X_d\in \mathscr{H}$ and let $\bm X=(X_1,\ldots,X_d)$. Let $F$ be a closed set in $\mathbb R^d$, let $F^{\delta}=\{\bm x: \|\bm x-\bm y\|<\delta  \text{ for some } \bm y\in F\}$ be the $\delta$-neighborhood of $F$. Then there exists a   Lipschitz function $f_{F,\delta}$ such that
$$ I_F\le f_{F,\delta}\le I_{F^{\delta}} $$
(c.f., page 9 of Billingsley (1999)).
By \eqref{eq1.3}  and \eqref{eqlemlinear.1}, it follows that
$$ \upCapc(\bm X \in F)\le \Sbep[f_{F,\delta}(\bm X)]=P[f_{F,\delta}(\bm X)]\le P(\bm X\in F^{\delta}). $$
Letting $\delta\to 0$ yields
$$ \upCapc(\bm X \in F)\le  P(\bm X \in F)  $$
by the continuity of $P$.
Suppose that $O$ is an open set. Then $F=O^c$ is a closed set. So
$$ I_O=1-I_F\ge 1-f_{F,\delta}\ge 1- I_{F^{\delta}}. $$
It follows that
$$ \upCapc(\bm X \in O)\ge \Sbep[1-f_{F,\delta}(\bm X)]=P[1-f_{F,\delta}(\bm X)]\ge 1-P(\bm X\in F^{\delta}), $$
by \eqref{eq1.3}  and \eqref{eqlemlinear.1}  again.
Letting $\delta\to 0$ yields
$$ \upCapc(\bm X\in O)\ge 1-P(\bm X\in F)=P(\bm X\in O). $$
Now, for a set $A\in \mathcal{B}(\mathbb R^d)$, let $A^o$ and $A^-$ be the interior and closure of $A$, and $\partial A=A^-\setminus A^o$ be the boundary of $A$. Then
$$ P(\bm X\in A^o)\le \upCapc(\bm X\in A^o)\le \upCapc(\bm X\in A)\le \upCapc(\bm X\in A^-)\le P(\bm X\in A^-). $$
Hence
\begin{equation}\label{eqlemlinear.3}  \upCapc(\bm X\in A)=P(\bm X\in A)  \;\; \text{ whenever } P(\bm X\in \partial A)=0.
\end{equation}
Let $\mathcal F_0$ be the family of sets of the form $A=\{(X_1,\ldots, X_d)\in B\}$ for some group of  $d$, $B\in \mathcal{B}(R^d)$, and elements $X_1,\ldots,X_d\in \mathscr{H}$ satisfying  $P((X_1,\ldots, X_d)\in \partial B)=0$. Then $\mathscr{F}_0\subset \sigma(\mathscr{H})$ is a field. Define
$$P^{\ast}(A)=\inf\{\sum_{n=1}^{\infty}P(A_n): \forall A_n\in \mathcal{F}_0 \text{ with } A\subset \bigcup_{n=1}^{\infty}A_n\Big\}.$$
Then   $P^{\ast}$ is an outer measure extension of $P\big|_{\mathcal{F}_0}$. By the measure extension theorem,
$$ P^{\ast}(A)=P(A), \;\; \forall A\in \sigma(\mathcal{F}_0), $$
cf. section 4.1 of Lo\`eve (1977) or Theorems 2 and 4 in \S 5.2 of Yan, Wang and Liu (1997).

On the other hand,  it is obvious that
\begin{align*}
\outCapc(A)=&\inf\Big\{\sum_{n=1}^{\infty}\upCapc(A_n): A\subset \bigcup_{n=1}^{\infty}A_n\Big\}\\
\le & \inf\Big\{\sum_{n=1}^{\infty}\upCapc(A_n): \forall A_n\in \mathcal{F}_0 \text{ with } A\subset \bigcup_{n=1}^{\infty}A_n\Big\}\\
=&\inf\Big\{\sum_{n=1}^{\infty}P(A_n): \forall A_n\in \mathcal{F}_0 \text{ with } A\subset \bigcup_{n=1}^{\infty}A_n\Big\} \;\; (by\;\; \eqref{eqlemlinear.3}) \\
 =&P^{\ast}(A).
\end{align*}
Hence, for \eqref{eqlemlinear.2}  it is sufficient to show that $\sigma(\mathscr{H})=\sigma(\mathcal{F}_0)$. Let $X_1,\ldots,X_d\in \mathscr{H}$ and $O\subset \mathbb R^d$ be an open set. For each $\bm x\in O$, there exists an open ball $B(\bm x,\delta_x)=\{\bm y; |\bm x-\bm y\|<\delta_x)$ such that $B(\bm x,\delta_x)\subset O$. Notice that
$\partial B(\bm x,\delta)\subset \{\bm y;\|\bm x-\bm y\|=\delta\}$, and so $\{\partial B(\bm x,\delta); 0<\delta\le \delta_x\}$ is a family of disjoint sets. There are at most countable number of them satisfying $P(\bm X\in \partial B(\bm x,\delta))\ne 0$. Therefore,  there exists a  $0<\delta_x^{\prime}\le \delta_x$ such that $\{\bm X\in B(\bm x,\delta_x^{\prime})\}\in \mathcal{F}_0$. Now, $\{B(\bm x,\delta_x^{\prime}); \bm x\in O\}$ is a cover of $O$. So, there exists a countable subset $U$ of $O$ such that
$ O=\bigcup_{\bm x\in U} B(\bm x,\delta_x^{\prime}). $
Then
$$\{\bm X\in O\}=\bigcup_{\bm x\in U} \{\bm  X\in B(\bm x,\delta_x^{\prime})\}\in \sigma(\mathcal F_0).  $$
Therefore,
$$ \sigma(X_1,\ldots,X_d)=\big\{\{\bm X\in B\}; B\in \mathcal{B}(\mathbb R^d)\big\} \subset \sigma(\mathcal F_0).  $$
By the arbitrariness of $X_1,\ldots,X_d\in \mathscr{H}$, we have $\sigma(\mathscr{H}) \subset \sigma(\mathcal F_0)$. The proof of  \eqref{eqlemlinear.2}is completed.

Now, suppose that $\Capc$ is a countably sub-additive capacity with the property \eqref{eq1.3}, then $\Capc(A)\le \outCapc(A)=P(A)$ for all $A\in \sigma(\mathscr{H})$.
With the same argument for $\upCapc$, we can show that $\Capc(A)=P(A)$ for all $A\in \mathcal F_0$. Let $\mathcal{M}=\{A\in \sigma(\mathscr{H}); \Capc(A)=P(A)\}$. It is sufficient to show that $\mathcal{M}$ is a monotone class, and then $\sigma(\mathscr{H})=\sigma(\mathcal F_0)\subset \mathcal{M}$. Firstly, if $A\in \mathcal M$, then $P(A^c)\ge \Capc(A^c)\ge 1-\Capc(A)=1-P(A)=P(A^c)$, and so $A^c\in \mathcal M$. Secondly,  for $\mathcal M\ni A_n\uparrow A$, we have $P(A)\ge \Capc(A)\ge \lim_{n\to \infty}\Capc(A_n)=\lim_{n\to\infty}P(A_n)=P(A)$, and then $A\in \mathcal M$. Hence, $\mathcal M$ is a monotone class and the proof is completed.
$\Box$.

   \section{The law of the iterated logarithm}\label{sectLIL}
\setcounter{equation}{0}

\subsection{General results}
 We state the results for the general LIL. The first two theorems describe Wittmann's LIL for independent random variables which are   not necessarily identically distributed. Let $\{X_n; n\ge 1\}$ be a sequence of   independent random variables in a sub-linear expectation space $(\Omega,\mathscr{H},\Sbep)$ with a sub-additive capacity $\Capc$ satisfying \eqref{eq1.3}. Denote  $s_n^2=\sum_{k=1}^n \Sbep[X_k^2]$, $t_n=\sqrt{2\log\log s_n^2}$, $a_n=s_nt_n$,
  $$\Gamma_n(p, \alpha)=\Sbep\Big[\Big(\big(|X_n|-\alpha  s_n/t_n \big)^+\Big)^p\Big], \overline{\Gamma}_n(p, \alpha)=\Sbep\Big[\Big(\big(|X_n|\wedge a_n-\alpha  s_n/t_n \big)^+\Big)^p\Big] $$
and
$$\Lambda_n(p,\alpha)=\sum_{j=1}^n \Sbep\Big[\Big(\big(|X_j|-\alpha  s_n/t_n \big)^+\Big)^p\Big], \;\; \overline{\Lambda}_n(p,\alpha)=\sum_{j=1}^n \Sbep\Big[\Big(\big(|X_j|\wedge a_n-\alpha  s_n/t_n \big)^+\Big)^p\Big]. $$

 \begin{theorem} \label{thLILnonid1}Let  $\{X_n; n\ge 1\}$ be a sequence of  independent random variables in the sub-linear expectation space $(\Omega,\mathscr{H},\Sbep)$ with $s_n^2\to \infty$.
    Suppose that
     \begin{equation}
     \label{eqthLILnonid1.1} \sum_{n=1}^{\infty} \Capc\left(|X_n|\ge \epsilon a_n\right)<\infty \; \text{ for all } \epsilon>0,
     \end{equation}
    and,  for every $\alpha>0$ there exist   some $p\ge 2$ and $d\ge 0$ such that
   \begin{equation}\label{eqthLILnonid1.2} \sum_{n=1}^{\infty}\frac{\overline{\Gamma}_n(p, \alpha)}{a_n^p}\left(\frac{ \overline{\Lambda}_n(p, \alpha)}{a_n^p}\right)^d<\infty.
    \end{equation}
   Then, for every $\epsilon>0$,
    \begin{equation}\label{eqthLILnonid1.3}
     \lim_{n\to \infty} \max_N\Capc\Big(\inf_{ n\le m\le N}  \frac{ S_m-\cSbep[S_m] }{a_m}<   -1 -\epsilon   \text{ or }   \sup_{ n\le m\le N}  \frac{ S_m-\Sbep[S_m] }{a_m}>   1 +\epsilon  \Big)=0.
     \end{equation}
     If $\Capc$ is countably sub-additive, then
     \begin{equation}\label{eqthLILnonid1.3ad}
     \lim_{n\to \infty}  \Capc\Big(\inf_{ m\ge n}  \frac{ S_m-\cSbep[S_m] }{a_m}<   -1 -\epsilon   \text{ or }   \sup_{m\ge n}  \frac{ S_m-\Sbep[S_m] }{a_m}>   1 +\epsilon  \Big)=0.
     \end{equation}
If $\outCapc$ is defined as in \eqref{outcapc},    then
 \begin{equation}\label{eqthLILnonid1.4}
  \outCapc\left( \liminf_{n\to \infty}\frac{S_n-\cSbep[S_n]}{a_n}< -1 \;\text{ or }\; \limsup_{n\to \infty}\frac{S_n-\Sbep[S_n]}{a_n}> 1 \right)=0.
 \end{equation}
\end{theorem}

 \begin{theorem} \label{thLILnonid2} Let  $\{X_n; n\ge 1\}$  be a sequence of independent random variables in the sub-linear expectation space $(\Omega,\mathscr{H},\Sbep)$.  Suppose \eqref{eqthLILnonid1.1}
    and that,  for every $\alpha>0$ there exist   some $p\ge 2$ and $d\ge 0$ such that
   \begin{equation}\label{eqthLILnonid2.1} \sum_{n=1}^{\infty}\frac{\Gamma_n(p, \alpha)}{a_n^p}\left(\frac{  \Lambda_n(p, \alpha)}{a_n^p}\right)^d<\infty,
    \end{equation}
    and further,
    \begin{equation}\label{eqthLILnonid2.2} \sum_{n=1}^{\infty} s_n^{-2} (\log s_n^2)^{\delta-1}\Sbep[X_n^2]= \infty\; \text{ for all } \delta>0,
  \end{equation}
 \begin{equation}\label{eqthLILnonid2.3}  \frac{\sum\limits_{j=1}^n |\Sbep[X_j]|+\sum\limits_{j=1}^n |\cSbep[X_j]|}{a_n}\to 0.
    \end{equation}
 Then, for every $\epsilon>0$,
\begin{equation} \label{eqthLILnonid2.4}
\lim_{n\to\infty}\max_N\Capc\left(\max_{n\le m\le N}\frac{|S_m|}{a_m}> 1+\epsilon\right)=0,
\end{equation}
\begin{align}\label{eqthLILnonid2.5ad}   & \lim_{ n\to \infty}\max_N\Capc\left( \max_{ n\le m\le N} \frac{-S_m}{a_m}\ge 1-\epsilon\right)
=  \lim_{  n\to \infty}\max_N\Capc\left(  \max_{ n\le m\le N} \frac{S_m}{a_m}\ge 1-\epsilon\right)=1.
\end{align}
Further, if $\outCapc$ is defined as in \eqref{outcapc},    then
\begin{equation}\label{eqthLILnonid2.6}
\outCapc\left(\limsup_{n\to \infty}\frac{|S_n|}{a_n}> 1\right)=0.
\end{equation}
\end{theorem}

The following are some remarks on the conditions.

\begin{remark}\label{remark3}  When $d=0$, \eqref{eqthLILnonid2.1}  and \eqref{eqthLILnonid1.2}  are
 \begin{equation}\label{eqremkLIL3.1} \sum_{n=1}^{\infty}\frac{\Sbep\Big[\Big(\big(|X_n|-\alpha  s_n/t_n \big)^+\Big)^p\Big]}{a_n^p}<\infty
    \end{equation}
    and
\begin{equation}\label{eqremkLIL3.2} \sum_{n=1}^{\infty}\frac{\Sbep\Big[\Big(\big(|X_n|\wedge a_n-\alpha  s_n/t_n \big)^+\Big)^p\Big]}{a_n^p}<\infty,
    \end{equation}
respectively. Further, if  $n$ is large enough such that $\epsilon a_n/2>\alpha s_n/t_n$, then
$$ \frac{\epsilon}{2} I\{|X_n|\ge \epsilon a_n\}\le  \frac{\big(|X_n|\wedge a_n-\alpha s_n/t_n \big)^+}{a_n}\le \frac{\big(|X_n|-\alpha s_n/t_n \big)^+}{a_n}.
$$
Therefore, \eqref{eqremkLIL3.1}  implies \eqref{eqremkLIL3.2}, and \eqref{eqremkLIL3.2}  implies \eqref{eqthLILnonid1.1}. It follows that, if \eqref{eqremkLIL3.1}  or \eqref{eqremkLIL3.2}  is satisfied, then \eqref{eqthLILnonid1.1}  can be removed.
\end{remark}

By Theorem \ref{thLILnonid2} and Remark \ref{remark3},   we have the following corollary.
\begin{corollary} Let    $\{X_n; n\ge 1\}$  be a sequence of  independent random variables in the sub-linear expectation space $(\Omega,\mathscr{H},\Sbep)$. Suppose that \eqref{eqthLILnonid2.2}  and \eqref{eqthLILnonid2.3}  hold, and for every $\alpha>0$, there exists $p\ge 2$ such that \eqref{eqremkLIL3.1}  holds. Then \eqref{eqthLILnonid2.4} -\eqref{eqthLILnonid2.6}  hold.
\end{corollary}

\begin{remark} \label{remark1}    As shown by Wittmann, \eqref{eqthLILnonid2.2}  is implied by
 \begin{equation}\label{eqthLILnonid2.0} s_n^2=\sum_{i=1}^n \Sbep[X_i^2]\to \infty \;\;\text{ and }\;\; \limsup\limits_{n\to\infty} \frac{s_{n+1}}{s_n}<\infty.
  \end{equation}
\end{remark}
In fact,
\begin{align*}
  & \sum_{n=3}^{\infty} s_n^{-2} (\log s_n^2)^{\delta-1}\Sbep[X_n^2]\\
  =&\sum_{n=3}^{\infty}\int_{s_{n-1}^2}^{s_n^2} s_n^{-2} (\log s_n^2)^{\delta-1}dx
\ge   c\sum_{n=3}^{\infty}\int_{s_{n-1}^2}^{s_n^2} s_{n-1}^{-2} (\log s_{n-1}^2)^{\delta-1}dx \\
\ge & c\sum_{n=3}^{\infty}\int_{s_{n-1}^2}^{s_n^2} x (\log x)^{\delta-1}dx=c\int_{s_2^2}^{\infty} x (\log x)^{\delta-1}dx=\infty.
\end{align*}

\begin{remark}\label{remark4} If $\{X_n;n\ge 1\}$ satisfies Kolmogorov's (1929) condition as
\begin{equation}\label{eqremkLIL4.1}s_n^2   \to \infty,\;\; |X_n|\le \alpha_n\frac{s_n}{t_n}, \;n=1,2,\ldots  \text{ and }  \alpha_n\to 0,
\end{equation}
then the conditions \eqref{eqremkLIL3.1}  (thus \eqref{eqremkLIL3.2} ) and \eqref{eqthLILnonid2.0}  are satisfied.
\end{remark}
Hence, we have the following corollary.
\begin{corollary} ({\em Kolmogorov's LIL}) Let    $\{X_n; n\ge 1\}$  be a sequence of  independent random variables in the sub-linear expectation space $(\Omega,\mathscr{H},\Sbep)$. Suppose that \eqref{eqthLILnonid2.3}  and \eqref{eqremkLIL4.1}  hold. Then \eqref{eqthLILnonid2.4} -\eqref{eqthLILnonid2.6}  hold.
\end{corollary}

\begin{remark} \label{remark2} It is obvious that \eqref{eqthLILnonid2.1}  implies \eqref{eqthLILnonid1.2}. If we have an additional condition $s_n^{-2}\sum_{j=1}^n \Sbep[X_j^2\wedge a_j^2]\to 1$, then \eqref{eqthLILnonid2.1}   in Theorem \ref{thLILnonid2} can be weakened to \eqref{eqthLILnonid1.2}.
Also, it is obvious that
$$ \frac{\overline{\Lambda}_n(p,\alpha)}{a_n^p}\le \frac{\overline{\Lambda}_n(2,\alpha)}{a_n^2}\le \frac{\sum_{j=1}^n\Sbep[X_j^2]}{a_n^2}\le \frac{1}{2\log\log s_n^2}, $$
and therefore, \eqref{eqthLILnonid1.2}  is satisfied if
\begin{equation}\label{eqremkLIL2.1} \sum_{n=1}^{\infty}\frac{\overline{\Gamma}_n(p, \alpha)}{s_n^{p/2} (\log\log s_n^2)^{d^{\prime}}}<\infty \; \text{ for some } d^{\prime}>0.
    \end{equation}

 As for the condition \eqref{eqthLILnonid2.1},   when $p=2$ it is just  Theorem 4.1 (i) of Wittmann (1987). Hence, Theorem \ref{thLILnonid2} has extended Wittmann's LIL   in even the classical case.

Notice $\Gamma_n(2, 2\alpha)/a_n^2\le \frac{\Gamma_n(p, \alpha)}{s_n^{p/2}}\alpha^{2-p} t_n^{2p-2}$, $\frac{\Lambda_n(2,2\alpha)}{a_n^2}\le t_n^{-2}$. Hence, if for every $\alpha>0$, there exist constants $p\ge 2$ and $d^{\prime}>0$ such that
\begin{equation}\label{eqremkLIL2.2} \sum_{n=1}^{\infty}\frac{\Gamma(p, \alpha)}{s_n^{p/2} (\log\log s_n^2)^{d^{\prime}}}<\infty,
    \end{equation}
then \eqref{eqthLILnonid1.2}  is satisfied with $p=2$ and every $\alpha>0$.
\end{remark}

Next, we consider the i.i.d. case. For a random variable $X$, we  denote
$$\breve{\mathbb E}[X]=\lim_{c\to \infty} \Sbep[(-c)\vee (X\wedge c)]$$  if the limit exists.
It can be verified that $\breve{\mathbb E}[X]$ exists if $C_{\Capc}(|X|)<\infty$ or $\breve{\mathbb E}[|X|^{1+\epsilon}]<\infty$,
$\breve{\mathbb E}[|X|]\le C_{\Capc}(|X|)$. Further, $\breve{\mathbb E}[X]=\Sbep[X]$ if $\Sbep[(|X|-c)^+]\to 0$ as $c\to \infty$.

The following two theorems on the LIL for a sequence of independent and identically distributed random variables are corollaries of  Theorems \ref{thLILnonid1} and \ref{thLILnonid2}.

   \begin{theorem} \label{thLILiid1}
   Let     $\{Y_n; n\ge 1\}$  be  independent and identically distributed random variables in the sub-linear expectation space $(\Omega,\mathscr{H},\Sbep)$, and $\outCapc$ be defined as in \eqref{outcapc}.
   Suppose
   \begin{equation}\label{eqLILmomentcondition1} C_{\Capc}\left[\frac{Y_1^2}{\log\log|Y_1|}\right]<\infty.
    \end{equation}
      Denote
 $\overline{\sigma}_2^2=\breve{\mathbb E}[(Y_1-\breve{\mathbb E}[Y_1])^2]$  and  $\overline{\sigma}_1^2=\breve{\mathbb E}[(Y_1+\breve{\mathbb E}[-Y_1])^2]$ (finite or infinite).
     Then
\begin{equation}\label{eqthLILiid1.3}
\outCapc\left(\liminf_{n\to \infty}\frac{\sum_{i=1}^n (Y_i+\breve{\mathbb E}[-Y_i])}{\sqrt{2n\log\log n}}< -\overline{\sigma}_1\;\text{ or }\;   \limsup_{n\to \infty}\frac{\sum_{i=1}^n (Y_i-\breve{\mathbb E}[Y_i])}{\sqrt{2n\log\log n}}> \overline{\sigma}_2 \right)=0.
\end{equation}
\end{theorem}

   \begin{theorem} \label{thLILiid2}
   Let   $\{Y_n; n\ge 1\}$  be independent and identically distributed random variables in the sub-linear expectation space $(\Omega,\mathscr{H},\Sbep)$, and $\outCapc$ be defined as in \eqref{outcapc}. Denote
 $\overline{\sigma}^2=\breve{\mathbb E}[Y_1^2]$ (finite or infinite).
  \begin{description}
  \item[\rm (a)]
 Suppose \eqref{eqLILmomentcondition1}  and
    \begin{equation}\label{eqLILmeanzerocondition}\breve{\mathbb E}[Y_1]=\breve{\mathbb E}[-Y_1]=0. \end{equation}
    Then
\begin{equation}\label{eqthLILiid2.3}
\outCapc\left( \limsup_{n\to \infty}\frac{|\sum_{i=1}^n Y_i|}{\sqrt{2n\log\log n}}> \overline{\sigma} \right)=0,
\end{equation}
\begin{equation}\label{eqthLILiid2.4}
\lim_{n\to \infty}\max_N\Capc\left(\sigma_1\le \sup_{ n\le m\le N}\frac{\sum_{i=1}^m Y_i}{\sqrt{2m\log\log m}}\le \sigma_2\right)=1  \; \text{ for all } \sigma_1<\overline{\sigma}<\sigma_2.
\end{equation}
\item[\rm (b)]   Suppose that there exist a $n_0$ and $M$ such that
\begin{equation}\label{eqthLILiid2.6}
 \lim_{N\to \infty}\Capc\left(   \sup_{n\le m\le N}\frac{|\sum_{i=1}^m Y_i|}{\sqrt{2m\log\log m}}\ge M \right)<1 \; \text{ for all } n\ge n_0.
\end{equation}
Then  \eqref{eqLILmomentcondition1}  and \eqref{eqLILmeanzerocondition}  hold, and
\begin{equation}\label{eqthLILiid2.8}\overline{\sigma}^2=\lim_{c\to \infty}\Sbep[Y_1^2\wedge c]<\infty. \end{equation}
\end{description}
\end{theorem}

\begin{remark} Theorem \ref{thLILiid2}   indicates  us that \eqref{eqLILmomentcondition1}, \eqref{eqLILmeanzerocondition}  and \eqref{eqthLILiid2.8}  are the sufficient and necessary conditions for a Hartman and Wintner  type LIL under  sub-linear expectations.

Compared with the LIL, the sufficient and necessary conditions of the central limit theorem for independent and identically distributed random variables are shown to be \eqref{eqLILmeanzerocondition},\eqref{eqthLILiid2.8}  and
$\Capc(|X_1|\ge x)=o(x^{-2})$ as $x\to\infty$  by Zhang (2020).
\end{remark}

The following theorem gives the result under the lower capacity $\cCapc$.
\begin{theorem}\label{thLIL5}
   Let  $\{Y_n; n\ge 1\}$  be independent and identically distributed random variables in the sub-linear expectation space $(\Omega,\mathscr{H},\Sbep)$, and $\outCapc$ be defined as in \eqref{outcapc}.     Recall
 $\overline{\sigma}^2=\breve{\mathbb E}[Y_1^2]$. Denote
 $\underline{\sigma}^2=\lim\limits_{c\to \infty}\cSbep[Y_1^2\wedge c]$, $T_n=\sum_{i=1}^nY_i$ and $d_n=\sqrt{2n\log\log n}$.   Suppose that \eqref{eqLILmomentcondition1}, \eqref{eqLILmeanzerocondition} and \eqref{eqthLILiid2.8}  are satisfied.
 Then
 $$ \outcCapc\left(\underline{\sigma}\le  \limsup_{n\to \infty}\frac{T_n}{d_n}\le \overline{\sigma}\right)=1, $$
    $$ \outcCapc\left(-\overline{\sigma}\le  \liminf_{n\to \infty}\frac{T_n}{d_n}\le -\underline{\sigma}\right)=1 $$
and
$$\outcCapc\left([-\overline{\sigma}, \;\;\overline{\sigma}]\supset C\Big\{\frac{T_n}{d_n}\Big\}=\Big[\liminf_{n\to \infty}\frac{T_n}{d_n},\;\;\limsup_{n\to \infty}\frac{T_n}{d_n}\Big]\supset
[-\underline{\sigma}, \;\;\underline{\sigma}]\right)=1, $$
where $C \{x_n\} $ denotes the cluster set of a sequence of $\{x_n\}$ in $\mathbb R$.
\end{theorem}
\begin{remark}
Theorem \ref{thLIL5} removes the continuity of $\Capc$ in Corollary 3.13 of Zhang (2016)  so that it is consistent with Theorem 1 of Chen and Hu (2014) where the random variables are assumed to be bounded.

The condition that $\lim\limits_{c\to \infty} \Sbep[(|Y_1|^2-c)^+]= 0$ in Theorem 3.11 and Corollary 3.13 of Zhang (2016) is now weakened to   \eqref{eqthLILiid2.8}. It should   be noted that \eqref{eqthLILiid2.8}, $  \Sbep[|Y_1|^2]< \infty$ and $\lim\limits_{c\to \infty} \Sbep[(|Y_1|^2-c)^+]= 0$ are not equivalent under the sub-linear expectation, and \eqref{eqthLILiid2.8}  is the weakest among them. Also, they do not imply \eqref{eqLILmomentcondition1}.
\end{remark}

For Theorem \ref{thLIL5}, we conjecture that it is also true when $\underline{\sigma}=\infty$.
\begin{conjecture}  If \eqref{eqLILmomentcondition1}, \eqref{eqLILmeanzerocondition}  and $ \lim\limits_{c\to \infty}\cSbep[Y_1^2\wedge c]=\infty$, then
$$ \outcCapc\left(   \liminf_{n\to \infty}\frac{T_n}{d_n}=-\infty \text{ and } \limsup_{n\to \infty}\frac{T_n}{d_n}=\infty\right)=1. $$
\end{conjecture}

\subsection{The exact lower bound}

If $\Capc$ is a continuous capacity, it is obvious that \eqref{eqthLILnonid2.4} and \eqref{eqthLILnonid2.5ad} imply
\begin{equation}\label{eqLILforProb1}     \Capc\left(\liminf_{n\to \infty} \frac{S_n}{a_n}= -1\right)
=  \Capc\left(\limsup_{n\to \infty} \frac{S_n}{a_n}= 1\right)=1,
\end{equation}
\eqref{eqthLILiid2.4} implies
\begin{equation}\label{eqLILforProb2} \Capc\left(\limsup_{n\to \infty}\frac{\sum_{i=1}^n Y_i}{\sqrt{2n\log\log n}}= \overline{\sigma}\right)=1,
\end{equation}
and \eqref{eqthLILiid2.6} is equivalent to
 $$ \Capc\left(  \limsup_{n\to \infty}\frac{|\sum_{i=1}^n Y_i|}{\sqrt{2n\log\log n}}= \infty \right)<1. $$
So Theorem \ref{thLILnonid2} and \ref{thLILiid2} include the LILs for independent random variables in a probability space. As we have shown,  a general capacity   is neither    continuous nor has the property \eqref{CapcC}. So, the converse part of the Borel-Cantelli lemma is not valid in general, and the lower bound of the LIL   becomes complex. In general, we  have no   result on the capacity of
$$ \limsup_{n\to \infty} \frac{S_n}{a_n}= 1 \; \text{ or } \; \limsup_{n\to\infty}\frac{\sum_{i=1}^n Y_i}{\sqrt{2n\log\log n}}=\overline{\sigma}^2. $$
Next, we consider two special   cases.  The first one is that the sub-linear expectation satisfies the condition (CC) in Proposition \ref{prop2}.

\begin{theorem} \label{thLILnew1} Let $(\Omega,\mathscr{H},\Sbep)$  be a sub-linear expectation space satisfying   the condition (CC) in Proposition \ref{prop2}.  And let  $\{X_n; n\ge 1\}$  be a sequence of independent random variables in  $(\Omega,\mathscr{H},\Sbep)$ satisfying the conditions in Theorem \ref{thLILnonid2}.     Then for $V=\Capc^{\mathscr{P}}$, $\mathbb C^{\ast}$ or $\outCapc$,
\begin{equation}\label{eqthLILnew1.1}
V\left(\limsup_{n\to \infty} \frac{|  S_n|}{a_n}>1\right)=0
\end{equation}
and
\begin{equation}\label{eqthLILnew1.2}
V\left(\liminf_{n\to \infty} \frac{  S_n}{a_n}=-1\; \text{ and } \;\limsup_{n\to \infty} \frac{  S_n }{a_n}=1 \right)=1,
\end{equation}
where $\Capc^{\mathscr{P}}$ is defined in Proposition \ref{prop2}.
\end{theorem}

\begin{theorem}\label{thLILnew2} Let $(\Omega,\mathscr{H},\Sbep)$  be a sub-linear expectation space satisfying  the condition (CC) in Proposition \ref{prop2} and having  a capacity $\Capc$  with the property \eqref{eq1.3}. Suppose that  $\{Y_n; n\ge 1\}$  is a sequence  of independent and identically distributed random variables under $\Sbep$.
\begin{description}
\item[\rm (a)]
 If $\eqref{eqLILmomentcondition1}$ and  $\eqref{eqLILmeanzerocondition}$, then for $V=\Capc^{\mathscr{P}}$, $\mathbb C^{\ast}$ or $\outCapc$,
 \begin{equation}\label{eqthLILnew2.2}
 V\left(  \limsup_{n\to\infty}\frac{|\sum_{i=1}^n  Y_i|}{\sqrt{2n\log\log n}}>\overline{\sigma} \right)=0,
\end{equation}
 \begin{equation}\label{eqthLILnew2.3}
 V\left(  C \left\{\frac{\sum_{i=1}^n  Y_i}{\sqrt{2n\log\log n}}\right\}=[-\overline{\sigma},\overline{\sigma}]\right)=1.
\end{equation}
\item[\rm (b)]    If $V=\Capc^{\mathscr{P}}$, $\mathbb C^{\ast}$ or $\outCapc$,
\begin{equation}\label{eqthLILnew2.5}
V\left(   \limsup_{n\to \infty}\frac{|\sum_{i=1}^n  Y_i|}{\sqrt{2n\log\log n}}=+\infty \right)<1.
\end{equation}
then \eqref{eqLILmomentcondition1}, \eqref{eqLILmeanzerocondition}  and \eqref{eqthLILiid2.8} hold.
\end{description}
\end{theorem}

In the second special case,  we consider the copy of the random sequences. We show that    we can redefine the space and random variables on $\mathbb R^{\infty}$ such that  \eqref{eqLILforProb1} and \eqref{eqLILforProb2} hold.

We consider the real space $\widetilde{\Omega}=\mathbb R^{\infty}=\{\bm x=(x_1,x_2,\ldots); x_i\in \mathbb R\}$ with the Borel $\sigma$-field $\widetilde{\mathcal F}=\mathscr{B}(\mathbb R^{\infty})$. Define the function space $\widetilde{\mathscr{H}}=\{\varphi; \varphi(\bm x)=\varphi_1(x_1,\ldots,x_d); \varphi_1\in C_{l,Lip}(\mathbb R^d), d\ge 1\}$. For a sequence $\{X_n;n\ge 1\}$ of random variables on $(\Omega,\mathscr{H},\Sbep)$, we define a copy $\{\tilde X_n;n\ge 1\}$ on $(\widetilde{\Omega}, \widetilde{\mathscr{H}})$ as follows.
First, define the sub-linear expectation on $(\widetilde{\Omega},\widetilde{\mathscr{H}})$ by
\begin{equation}\label{eqcopyE} \widetilde{\mathbb E}[\varphi]=\Sbep[\varphi_1(X_1,\ldots,X_d)]=\Sbep[\varphi(\bm X)], \;\; \varphi=\varphi_1\circ\pi_d, \varphi_1\in C_{l,Lip}(\mathbb R^d),
\end{equation}
where $\pi_d: \mathbb R^{\infty}\to \mathbb R^d$ is the projection map,  $\pi_d\bm x=(x_1,\ldots, x_d)$.
On the space $(\widetilde{\Omega},\widetilde{\mathscr{H}},\widetilde{\mathbb E})$,
let
\begin{equation} \label{eqcopyPs}
\begin{aligned}
 \widetilde{\mathscr{P}}=  \Big\{P: P & \text{ is a probability measure on } (\widetilde{\Omega},\widetilde{\mathcal F}) \text{ satisfying} \\
& \;\;  P[\varphi]\le \widetilde{\mathbb E}[\varphi]  \;\; \text{ for all bounded } \varphi \in \widetilde{\mathscr{H}}\Big\}.
\end{aligned}
\end{equation}
  Define a set function on $\widetilde{\mathcal F}$ by
\begin{equation}\label{eqcopyV} \widetilde{\Capc}^{\widetilde{\mathscr{P}}}(A)=\sup_{P\in \widetilde{\mathscr{P}}}P(A), \;\; A\in \widetilde{\mathcal F},
\end{equation}
where $\widetilde{\Capc}^{\widetilde{\mathscr{P}}}\equiv 0$ if $\widetilde{\mathscr{P}}$ is empty. It is obvious that $\widetilde{\Capc}^{\mathscr{P}}$ is countably sub-additive and  $\widetilde{\Capc}^{\widetilde{\mathscr{P}}}(A)\le \widetilde{\mathbb C}^{\ast}(A)\le \widetilde{\Capc}^{\ast}(A)\le \widetilde{\Capc}(A)$ for any $A\in \widetilde{\mathcal F}$, where $\widetilde{\mathbb C}^{\ast}$, $\widetilde{\Capc}^{\ast}$ and $\widetilde{\Capc}$ are defined on $(\widetilde{\Omega},\widetilde{\mathscr{H}},\widetilde{\mathbb E})$ the same as $\mathbb C^{\ast}$,  $\outCapc$ and $\upCapc$ in \eqref{tildecapc}, \eqref{outcapc} and \eqref{eq1.5}, respectively. It can shown that
$\mathbb C^{\ast}(\bm X\in A)\le \widetilde{\mathbb C}^{\ast}(A)$,  $\outCapc(\bm X\in A)\le \widetilde{\Capc}^{\ast}(A)$ and $\upCapc(\bm X\in A)\le \widetilde{\Capc}(A)$ for all $A\in \widetilde{\mathcal F}$.

Now, define the random variable $\tilde X_n$ by $\tilde X_n(\widetilde{\omega})=x_n$ for $\widetilde{\omega}=\bm x$. Then
\begin{equation}\label{eqcopy} \widetilde{\mathbb E}[\varphi(\tilde X_1,\ldots,\tilde X_d)]= \widetilde{\mathbb E}[\varphi\circ\pi_d]=\Sbep[\varphi(X_1,\ldots,X_d)], \;\; \varphi\in C_{l,Lip}(\mathbb R^d).
\end{equation}
Hence $(\tilde X_1,\ldots,\tilde X_d)\overset{d}=(X_1,\ldots,X_d)$, $d\ge 1$. It follows that $\{\tilde X_n;n\ge 1\}$ is a copy of $\{X_n;n\ge 1\}$. We call such a copy the copy of $\{X_n;n\ge 1\}$ on $\mathbb R^{\infty}$.

\begin{theorem} \label{thLILnew3}  Let  $\{X_n; n\ge 1\}$  be a sequence of independent random variables in the sub-linear expectation space $(\Omega,\mathscr{H},\Sbep)$ satisfying the conditions in Theorem \ref{thLILnonid2}. Denote its copy on $\mathbb R^{\infty}$ defined as above by   $\{\tilde X_n; n\ge 1\}$, and $\tilde S_n=\sum_{i=1}^n \tilde X_i$. Then for $\widetilde V=\widetilde{\Capc}^{\widetilde{\mathscr{P}}}$, $\widetilde{\mathbb C}^{\ast}$ or $\widetilde{\Capc}^{\ast}$,
\begin{equation}\label{eqthLILnew3.1}
\widetilde V\left(\limsup_{n\to \infty} \frac{|\tilde S_n|}{a_n}>1\right)=0
\end{equation}
and
\begin{equation}\label{eqthLILnew3.2}
\widetilde V\left(\liminf_{n\to \infty} \frac{\tilde S_n}{a_n}=-1\; \text{ and } \;\limsup_{n\to \infty} \frac{\tilde S_n}{a_n}=1 \right)=1.
\end{equation}
\end{theorem}

The following theorem is the result about the i.i.d. random variables.
\begin{theorem}\label{thLILnew4}  Let $(\Omega,\mathscr{H},\Sbep)$ be a sub-linear expectation space with a capacity $\Capc$ satisfying \eqref{eq1.3},   and   $\{Y_n; n\ge 1\}$  be  a sequence  of independent and identically distributed random variables under $\Sbep$. Let   $\{\tilde{Y}_n; n\ge 1\}$ be a copy on $\mathbb R^{\infty}$ of $\{Y_n; n\ge 1\}$.
\begin{description}
\item[\rm (a)]
 If $\eqref{eqLILmomentcondition1}$ and  $\eqref{eqLILmeanzerocondition}$, then for $\widetilde V=\widetilde{\Capc}^{\widetilde{\mathscr{P}}}$, $\widetilde{\mathbb C}^{\ast}$ or $\widetilde{\Capc}^{\ast}$,
 \begin{equation}\label{eqthLILnew4.2}
\widetilde  V\left(  \limsup_{n\to\infty}\frac{|\sum_{i=1}^n \tilde Y_i|}{\sqrt{2n\log\log n}}>\overline{\sigma} \right)=0,
\end{equation}
 \begin{equation}\label{eqthLILnew4.3}
\widetilde  V\left(  C \left\{\frac{\sum_{i=1}^n \tilde Y_i}{\sqrt{2n\log\log n}}\right\}=[-\overline{\sigma},\overline{\sigma}]\right)=1.
\end{equation}
Furthermore,
\begin{equation} \label{eqthLILnew4.4}
\widetilde V\left(   C \left\{\frac{\sum_{i=1}^n \tilde Y_i}{\sqrt{2n\log\log n}}\right\}=[-\sigma, \sigma]\right)=\begin{cases}
  1, & \text{ when } \sigma\in [\underline{\sigma},\overline{\sigma}],\\
  0, & \text{ when } \sigma\not\in [\underline{\sigma},\overline{\sigma}].
  \end{cases}
\end{equation}
\item[\rm (b)] Suppose that each $Y_n$ is tight in the sense that $\lim_{c\to \infty}\Capc(|Y_n|\ge c)=0$.  If for $\widetilde V=\widetilde{\Capc}^{\widetilde{\mathscr{P}}}$, $\widetilde{\mathbb C}^{\ast}$ or $\widetilde{\Capc}^{\ast}$,
\begin{equation}\label{eqthLILnew4.5}
\widetilde V\left(   \limsup_{n\to \infty}\frac{|\sum_{i=1}^n \tilde Y_i|}{\sqrt{2n\log\log n}}=+\infty \right)<1.
\end{equation}
Then \eqref{eqLILmomentcondition1}, \eqref{eqLILmeanzerocondition}   and \eqref{eqthLILiid2.8} hold.
\item[\rm (c)] Suppose that each $Y_n$ is tight.  If there exists a constant $b$ such that for $\widetilde V=\widetilde{\Capc}^{\widetilde{\mathscr{P}}}$, $\widetilde{\mathbb C}^{\ast}$ or $\widetilde{\Capc}^{\ast}$,
\begin{equation}\label{eqthLILnew4.6}
\widetilde v\left(   \limsup_{n\to \infty}\frac{ |\sum_{i=1}^n \tilde Y_i| }{\sqrt{2n\log\log n}}=b \right)>0,
\end{equation}
where $\widetilde v=1-\widetilde V$, then \eqref{eqLILmomentcondition1}, \eqref{eqLILmeanzerocondition}   and \eqref{eqthLILiid2.8} hold, and $\underline{\sigma}=\overline{\sigma}=b$.
\end{description}
\end{theorem}

\begin{remark} By \eqref{eqthLILnew4.4},
$$ \widetilde V\left(   \limsup_{n\to \infty}\frac{ |\sum_{i=1}^n \tilde Y_i| }{\sqrt{2n\log\log n}}=\sigma \right)=1 \;\; \text{ for all } \sigma\in [\underline{\sigma},\overline{\sigma}]. $$
However, by \eqref{eqthLILnew4.6}, under $\widetilde v$ the $\limsup$ can not be a constant unless $\underline{\sigma}=\overline{\sigma}$.
\end{remark}

{\bf Open problems:} (i) \eqref{eqthLILnew4.4} is shown for the copy  $\{\tilde Y_n;n\ge 1\}$. We conjecture that it holds for the original sequence $\{ Y_n;n\ge 1\}$.

(ii) Proving \eqref{eqthLILnew1.2}, \eqref{eqthLILnew2.3} and (b) of Theorem \ref{thLILnew2} under $\Capc=\mathbb C^{\ast}$ or $\outCapc$ without the condition (CC) in Proposition \ref{prop2}, or finding counterexamples, is an open problem.

(iii) For a sequence of independent and identically distributed random variables  in a classical probability space,  Martikainen (1980), Rosalsky (1980) and Pruitt (1981) proved that
if $\limsup\limits_{n\to\infty} \sum_{i=1}^n Y_i/\sqrt{2n\log\log  n}=1$ a.s., then $\ep Y_1=0$ and $\ep X_1^2=1$.  We conjecture that Theorem \ref{thLILnew2} (b) remains true when $|\sum_{i=1}^n Y_i|$ in \eqref{eqthLILnew2.5} is replaced by $\sum_{i=1}^n Y_i$.

  \section{Proofs of the laws of the iterated logarithm}\label{sectProof}
  \setcounter{equation}{0}

 In this section, we give the proofs of the theorems in Section \ref{sectLIL}.
\subsection{Proof of the general results}
We first we prove Theorems \ref{thLILnonid1}-\ref{thLIL5}.

{\bf Proof of Theorem \ref{thLILnonid1}}.
    By Wittmann (1985, Lemma 3.3), for any $\lambda>1$, there exists a sequence $\{n_k\}\subset \mathbb N$ with
 \begin{equation}\label{eqproofLILnonid1.1}\lambda a_{n_k}\le a_{n_{k+1}}\le \lambda^3 a_{n_k+1}.
 \end{equation}
 It can be checked that
  \begin{equation}\label{eqproofLILnonid1.2} \lambda  s_{n_k}^2\le s_{n_{k+1}}^2\le \lambda^6 s_{n_k+1}^2, \;\; \lambda^{1/2} \frac{s_{n_k}}{t_{n_k}}\le \frac{s_{n_{k+1}}}{t_{n_{k+1}}}\le \lambda^3
 \frac{s_{n_k+1}}{t_{n_k+1}}
  \end{equation}
 and $\log s_{n_{k+1}}^2\sim \frac{1}{2}\log a_{n_{k+1}}\ge c k$. Hence
  \begin{equation}\label{eqproofLILnonid1.4} \sum_{k=1}^{\infty}\exp\big\{- \frac{(1+\epsilon)t_{n_{k+1}}^2}{2} \big\}=\sum_{k=1}^{\infty}\exp\left\{-(1+\epsilon)\log\log s_{n_{k+1}}^2\right\}<\infty\;\;\text{ for all } \epsilon>0.
  \end{equation}
We write  $I(k)$ to denote the set $\{n_k+1,\ldots, n_{k+1}\}$. Denote $b_j=\alpha \lambda^3 s_j/t_j$, where $0<\alpha<1/10$ is a constant and to be specified.
Denote
$$\overline{\Lambda}_{n_k,n_{k+1}}(p,\alpha)=\sum_{j\in I(k)}\Sbep[\big((|X_j|\wedge a_{n_{k+1}}-\alpha s_{n_{k+1}}/t_{n_{k+1}})^+\big)^p]. $$
It follows from \eqref{eqthLILnonid1.2}  and \eqref{eqproofLILnonid1.1}  that (cf. the arguments of Wittmann (1987, page 526))
\begin{equation}\label{eqproofLILnonid1.5} \sum_{k=1}^{\infty}\left(\frac{\overline{\Lambda}_{n_k,n_{k+1}}(p,\lambda^3\alpha)}{a_{n_{k+1}}^p}\right)^{d+1}<\infty.
\end{equation}
 Let
\begin{equation} \label{eqproofLILnonid1.6} \mathbb N_1=\left\{k\in \mathbb N; \frac{\overline{\Lambda}_{n_k,n_{k+1}}(p,\lambda^3\alpha)}{a_{n_{k+1}}^p}\le t_{n_{k+1}}^{-2p}\right\}.
\end{equation}
It follows from \eqref{eqproofLILnonid1.5}  that
\begin{equation} \label{eqproofLILnonid1.7} \sum_{k\in \mathbb N\setminus \mathbb N_1}t_{n_{k+1}}^{-2p(d+1)}<\infty.
 \end{equation}

 We consider the sequences $\{X_j; j\in \mathbb N\setminus \mathbb N_1\}$ and $\{X_j; j\in \mathbb N_1\}$, respectively.
 Let $\overline{X}_j$ be $X_j$ if $j\in I(k)$ and $k\in \mathbb N\setminus\mathbb N_1$, and $0$ for otherwise. Denote $\widehat{X}_j=X_j-\overline{X}_j$.

 First, we consider $\{X_j; j\in \mathbb N\setminus \mathbb N_1\}$.  Denote $\overline{S}_n=\sum_{j=1}^n (\overline{X}_j-\Sbep[\overline{X}_j])$. Then
 $\overline{S}_n=\sum_{j\in I(k), k\in \mathbb N\setminus \mathbb N_1,j\le n}(X_j-\Sbep[X_j])$.
Let $x=\epsilon a_{n_{k+1}}$ and $y=\epsilon^{\prime}a_{n_{k+1}}$, where  $\epsilon^{\prime}>0$ is chosen  such that  $x/(10 y)\ge   d+1 $. By the inequality \eqref{eqExpIneq.4}  (with $\delta=1$) in Lemma \ref{lemExpIneq},
\begin{align*}
&\Capc\left(\max_{n\in I(k)}\frac{\sum_{j=n_k+1}^n(X_j- \Sbep[X_j])}{a_n}\ge \epsilon\lambda^3\right)\\
\le &\Capc\left(\max_{n\in I(k)}\sum_{j=n_k+1}^n(X_j- \Sbep[X_j])\ge \epsilon\lambda^3 a_{n_k+1}\right) \\
\le & \Capc\left(\max_{n\in I(k)}\sum_{j=n_k+1}^n(X_j - \Sbep[X_j])\ge \epsilon a_{n_{k+1}}\right)  \\
\le & \exp\left\{-\frac{\epsilon^2}{2(1+1)} t_{n_{k+1}}^2\right\}+\sum_{j\in I(k)} \Capc\big(|X_j|\ge \epsilon^{\prime} a_{n_{k+1}}\big) \\
&+C\left(\frac{\sum_{j\in I(k)}\Sbep[|X_j|^{2p}\wedge a_{n_{k+1}}^{2p}]}{a_{n_{k+1}}^{2p}}\right)^{d+1}.
\end{align*}
For $j\in I(k)$, we have
\begin{align*}
\frac{\Sbep[|X_j|^{2p}\wedge a_{n_{k+1}}^{2p}]}{a_{n_{k+1}}^{2p}}\le & \frac{\Sbep[|X_j|^{2p}\wedge b_{n_{k+1}}^{2p}]}{a_{n_{k+1}}^{2p}}
+\frac{\Sbep[\big((|X_j|\wedge a_{n_{k+1}}-b_{n_{k+1}})^+\big)^{2p}]}{a_{n_{k+1}}^{2p}} \\
\le & C\frac{\Sbep[ X_j^2]}{s_{n_{k+1}}^2t_{n_k+1}^{4p-2}}
+\frac{\Sbep[\big((|X_j|\wedge a_{n_{k+1}}-b_{n_{k+1}})^+\big)^{p}]}{a_{n_{k+1}}^{p}}.
\end{align*}
 It follows that
$$
\sum_{j\in I(k)} \frac{\Sbep[|X_j|^{2p}\wedge a_{n_{k+1}}^{2p}]}{a_{n_{k+1}}^{2p}}
\le   C t_{n_{k+1}}^{-2p}
+  \frac{  \overline{\Lambda}_{n_k,n_{k+1}}(p,\lambda^3\alpha)}{a_{n_{k+1}}^{2p}}.
$$
Therefore, for $k\not\in \mathbb N_1$,
\begin{align*}
&\Capc\left(\max_{n\in I(k)}\frac{\sum_{j=n_k+1}^n(X_j- \Sbep[X_j])}{a_n}\ge \epsilon\lambda^3\right)\\
\le & C t_{n_{k+1}}^{-2p(d+1)}  +C\left(\frac{\overline{\Lambda}_{n_k,n_{k+1}}(p,\lambda^3\alpha)}{a_{n_{k+1}}^{p}}\right)^{d+1}
+\sum_{j\in I(k)} \Capc\big(|X_j|\ge \epsilon^{\prime} a_j\big).
\end{align*}
Hence, it follows from \eqref{eqthLILnonid1.1}, \eqref{eqproofLILnonid1.5}  and  \eqref{eqproofLILnonid1.7}  that
\begin{equation} \label{eqproofLILnonid1.9}  \sum_{k\in \mathbb N\setminus \mathbb N_1}\Capc\left(\max_{n\in I(k)}\frac{\sum_{j=n_k+1}^n(X_j- \Sbep[X_j])}{a_n}\ge \epsilon \right)<\infty
\; \text{ for all } \; \epsilon>0.
\end{equation}
  That is,
$$  \sum_{k=1}^{\infty} \Capc\left(\max_{n\in I(k)}\frac{\overline{S}_n-\overline{S}_{n_k}}{a_n}\ge \epsilon \right)<\infty
\; \text{ for all } \; \epsilon>0.
$$
It follows from \eqref{eqproofLILnonid1.9}  that
\begin{align*}
&  \max_M\Capc\left(\sup_{K\le k\le M}\max_{n\in I(k)}\frac{\overline{S}_n-\overline{S}_{n_k}}{a_n}>\epsilon\right)\\
\le & \sum_{k=K}^{\infty} \Capc\left(\max_{n\in I(k)}\frac{\overline{S}_n-\overline{S}_{n_k}}{a_n}\ge \epsilon \right)\to 0 \text{ as }    K\to\infty\; \text{ for all } \; \epsilon>0,
\end{align*}
which implies that for any $\epsilon>0$ and $\delta>0$, there exists a $K_0\ge 1$ such that $\Capc(A_{K,M})\le \delta$ for all $M\ge K\ge K_0$, where
$$A_{K,M}=\left\{\sup_{K\le k\le M}\max_{n\in I(k)}\frac{\overline{S}_n-\overline{S}_{n_k}}{a_n}>\epsilon\right\}. $$
Choose a $K_1\ge K_0$ such that
$$ \Capc\left(\frac{\overline{S}_{n_{K_0}}}{a_{n_{K_1}}}\ge \epsilon\right)\le \frac{\Sbep[|\overline{S}_{n_{K_0}}|}{\epsilon a_{n_{K_1}}}<\delta. $$
Let $n_0=n_{K_1}+1$. Then for $n_0\le n\le m\le N$, there is a $K_1\le k\le M=:N$ such that $m\in I(k)$. Hence, on $A_{K_0,M}^c$ and $\left\{\frac{\overline{S}_{n_{K_0}}}{a_{n_{K_1}}}< \epsilon\right\}$,
\begin{align*}
\overline{S}_m=&(\overline{S}_m-\overline{S}_{n_k})+(\overline{S}_{n_k}-\overline{S}_{n_{k-1}})+\ldots+(\overline{S}_{n_{K_0+1}}-\overline{S}_{n_{K_0}})
+\overline{S}_{n_{K_0}}\\
< & \epsilon a_m+\epsilon a_{n_k}+\ldots+\epsilon a_{n_{K_0+1}}+\epsilon a_{n_{K_1}} \\
\le & \epsilon a_m \left(1+1+\frac{1}{\lambda}+\frac{1}{\lambda^2}+\ldots\right)+\epsilon a_m\le \epsilon a_m \frac{3\lambda}{\lambda-1},
\end{align*}
by \eqref{eqproofLILnonid1.1}.
It follows that
$$ \Capc\left(\sup_{n\le m\le N}  \frac{\overline{S}_m}{a_m}>\epsilon\frac{3\lambda}{\lambda-1}\right)\le \Capc(A_{K,M})+ \Capc\left(\frac{\overline{S}_{n_{K_0}}}{a_{n_{K_1}}}\ge \epsilon\right)\le 2\delta, \;\; N\ge n\ge n_0. $$
Hence
\begin{equation}\label{limsup1}\lim_{ n\to \infty} \max_N\Capc\left(\sup_{n\le m\le N}  \frac{\overline{S}_m}{a_m}>\epsilon\right)=0 \; \text{ for all } \; \epsilon>0.
\end{equation}
 Next, consider $\{X_j; j\in   \mathbb N_1\}$. We use the truncation method. Denote
$$Z_j=  \widehat{X}_j\wedge (2b_{n_{k+1}}),  j\in I(k), k\ge 0, $$
 and $\widehat{S}_n=\sum_{j=1}^n (Z_j-\Sbep[\widehat{X}_j])$.
Then
\begin{equation} \label{eqproofLILnonid1.10} S_n-\Sbep[S_n]=\overline{S}_n+\widehat{S}_n+\sum_{j=1}^n (\widehat{X}_j-Z_j).
\end{equation}
Notice that $\widehat{X}_j-Z_j= 0$ when $j\in I(k)$ and $k\in \mathbb N_1$, $=(X_j-2b_{n_{k+1}})^+\ge 0$ for otherwise. It is easily seen that
\begin{align*}
  \Sbep[(X_j\wedge a_{n_{k+1}}-2b_{n_{k+1}})^+]\le  &\Sbep\big[\big((X_{n_{k+1}}\wedge a_{n_{k+1}}-b_{n_{k+1}})^+\big)^p\big] b_{n_{k+1}}^{1-p} \\
  \Sbep[\big((X_j\wedge a_{n_{k+1}}-2b_{n_{k+1}})^+\big)^2] \le   &\Sbep\big[\big((X_j\wedge a_{n_{k+1}}-b_{n_{k+1}})^+\big)^p]b_{n_{k+1}}^{2-p}.
\end{align*}
It follows that for $k\in \mathbb N_1$,
\begin{align}
 &\frac{\sum_{j\in I(k)}
\Sbep[\widehat{X}_j\wedge a_{n_{k+1}}-Z_j]}{a_{n_{k+1}}}\le \frac{\overline{A}_{n_k,n_{k+1}}(p,\alpha\lambda^3)}{a_{n_{k+1}}^p} \alpha ^{1-p}t_{n_{k+1}}^{2p-2} \le \alpha^{1-p}t_{n_{k+1}}^{-2}\to 0,\label{eqproofLILnonid1.11}\\
 &\frac{\sum_{j\in I(k)}
  \Sbep[(\widehat{X}_j\wedge a_{n_{k+1}}-Z_j)^2] }{s_{n_{k+1}}^2}  \le \frac{\overline{A}_{n_k,n_{k+1}}(p,\alpha\lambda^3)}{a_{n_{k+1}}^p} \alpha ^{2-p}t_{n_{k+1}}^{2p-2} \le    \alpha^{2-p}t_{n_{k+1}}^{-2}\to 0.\label{eqproofLILnonid1.12}
\end{align}
Let $x=\epsilon a_{n_{k+1}}/2$ and $y=\epsilon^{\prime} a_{n_{k+1}}$, where $0<\epsilon^{\prime}<1$  is chosen such that $x/(10y)\ge d+1$. By the inequality \eqref{eqExpIneq.4}  (with $\delta=1$) in Lemma \ref{lemExpIneq}, we have that for $k\in \mathbb N_1$ large enough,
\begin{align}\label{eqproofLILnonid1.13}
&\Capc\Big(\sum_{j\in I(k)}|\widehat{X}_j-Z_j|\ge \epsilon a_{n_{k+1}}\Big)\nonumber\\
\le & \Capc\Big(\sum_{j\in I(k)}|\widehat{X}_j\wedge a_{n_{k+1}}-Z_j|\ge \epsilon a_{n_{k+1}}\Big)+\sum_{j\in I(k)}\Capc(X_j>a_{n_{k+1}}) \nonumber\\
\le  &\Capc\Big(\sum_{j\in I(k)}(|\widehat{X}_j\wedge a_{n_{k+1}}-Z_j|-\Sbep[|\widehat{X}_j\wedge a_{n_{k+1}}-Z_j|])\ge \epsilon a_{n_{k+1}}/2\Big)+\sum_{j\in I(k)}\Capc(X_j>a_j) \nonumber\\
\le  &  \exp\left\{-\frac{\epsilon^2/4}{2(1+1)} \frac{a_{n_{k+1}}^2}{\sum_{j\in I(k)}\Sbep[(\widehat{X}_j\wedge a_j-Z_j)^2] }\right\}\nonumber\\
& +C\left(\sum_{j\in I(k)}\frac{\Sbep\big[\big((|X_j|\wedge a_{n_{k+1}}-b_{n_{k+1}})^+\big)^p\big]}{a_{n_{k+1}}^p}\right)^{d+1}+2\sum_{j\in I(k)}\Capc(X_j>\epsilon^{\prime} a_{n_{k+1}}) \nonumber \\
\le &  \exp\left\{-2 t_{n_{k+1}}^2\right\}
    +C\left(\frac{\overline{\Lambda}_{n_k,n_{k+1}}(p,\lambda^3\alpha)]}{a_{n_{k+1}}^p}\right)^{d+1}+2\sum_{j\in I(k)}\Capc(X_j>\epsilon^{\prime} a_j).
\end{align}
Notice $\widehat{X}_j-Z_j=0$ when $j\in I(k)$ and $k\in \mathbb N\setminus \mathbb N_1$. It follows from \eqref{eqproofLILnonid1.13}, \eqref{eqproofLILnonid1.4}  and \eqref{eqproofLILnonid1.5}  and \eqref{eqthLILnonid1.1}  that
$$\sum_{k=1}^{\infty} \Capc\Big(\sum_{j\in I(k)}|\widehat{X}_j-Z_j|\ge \epsilon a_{n_{k+1}}\Big)=\sum_{k\in \mathbb N_1} \Capc\Big(\sum_{j\in I(k)}|\widehat{X}_j-Z_j|\ge \epsilon a_{n_{k+1}}\Big)<\infty. $$
Thus,
$$ \lim_{K\to \infty}\max_M\Capc\Big(\sup_{ K\le k\le M}\frac{\sum_{j\in I(k)}|\widehat{X}_j-Z_j|}{a_{n_{k+1}}}>\epsilon\Big)=0\; \text{ for all } \; \epsilon>0, $$
which implies
\begin{equation}\label{limsup2} \lim_{n\to\infty}\max_N\Capc\Big(\sup_{n\le m\le N}\frac{\sum_{j=1}^m|\widehat{X}_j-Z_j|}{a_m}>\epsilon\Big)=0 \; \text{ for all } \; \epsilon>0,
\end{equation}
similar to \eqref{limsup1}.
 At last, we consider $Z_j$. For $0<\epsilon<1/2$, choose $\alpha>0$ such that $8\alpha\lambda^3<\epsilon$.
  Notice that $\Sbep[Z_j]\le \Sbep[\widehat{X}_j]$ and $\sum_{j=1}^n \Sbep[Z_j^2]\le s_n^2$. Let $y_k= 2b_{n_{k+1}}$, $x_k=(1 +\epsilon)   a_{n_{k+1}}$. Applying \eqref{eqExpIneq.3}  in Lemma \ref{lemExpIneq} yields
 \begin{align*}
 &\Capc\Big(\max_{n\in I(k)}\frac{\widehat{S}_n}{a_n}\ge  (1 +\epsilon)\lambda^3    \Big)
 \le \Capc\Big(\max_{n\in I(k)} \widehat{S}_n \ge  (1 +\epsilon) a_{n_{k+1}}    \Big) \\
 \le & \exp\left\{-\frac{x_k^2}{2(x_ky_k+\sum_{j=1}^{n_{k+1}}\Sbep[Z_j^2])}\right\} \le \exp\left\{-\frac{x_k^2}{2(x_ky_k+s_{n_{k+1}}^2)}\right\} \\
 = & \exp\left\{-\frac{(1+\epsilon)^2t_{n_{k+1}^2}}{2(1+8\alpha\lambda^3)}\right\}
 \le  \exp\{-(1+\epsilon)\log\log s_{n_{k+1}}^2\},
 \end{align*}
 when $k$ is large enough. It follows by\eqref{eqproofLILnonid1.4}  that
 $$ \sum_{k=1}^{\infty} \Capc\Big(\max_{n\in I(k)}\frac{\widehat{S}_n}{a_n}\ge  (1 +\epsilon )\lambda^3    \Big)<\infty\;\; \forall \epsilon>0, $$
 which implies
 \begin{equation}\label{limsup3}
\lim_{n\to \infty} \max_N\Capc\Big(\sup_{n\le m\le N} \frac{\widehat{S}_m}{a_m}>   (1 +\epsilon )\lambda^3    \Big)=0.
 \end{equation}
By combining \eqref{eqproofLILnonid1.10}, \eqref{limsup1}, \eqref{limsup2}  and \eqref{limsup3}, it follows that
$$ \lim_{n\to \infty} \max_N\Capc\Big(\sup_{n\le m\le N}  \frac{S_m-\Sbep[S_m]}{a_m}>   (1 +\epsilon )\lambda^3  +\epsilon  \Big)=0\; \text{ for all } \;  \epsilon>0 \;\text{ and } \;\lambda>1. $$
 Therefore,
\begin{align*}
  \lim_{n\to \infty} \max_N\Capc\Big(\sup_{n\le m\le N}  \frac{S_m-\Sbep[S_m]}{a_m}>  1+\epsilon    \Big) =0\; \text{ for all } \;   \epsilon>0.
\end{align*}
When $\Capc$ is countably sub-additive, the $\sup_{K\le k\le M}$ in the brackets can be replaced by $\sup_{K\le k<\infty}$, and so
\begin{align*}
  \lim_{n\to \infty} \Capc\Big(\sup_{m\ge n}  \frac{S_m-\Sbep[S_m]}{a_m}>  1+\epsilon    \Big) =0\; \text{ for all } \;   \epsilon>0.
\end{align*}
For $-X_j$s, we have the same result. The proofs of \eqref{eqthLILnonid1.3}  and \eqref{eqthLILnonid1.3ad}  are now completed.

For \eqref{eqthLILnonid1.4}, let $\epsilon_k=1/2^k$. Notice that $\Capc$ and $\upCapc$ satisfy \eqref{eq1.3}. By \eqref{eqV-V}, \eqref{eqthLILnonid1.1}  is equivalent to that as it holds for $\upCapc$. So, \eqref{eqthLILnonid1.3}  holds for $\upCapc$. Then, there is a sequence of $n_k\uparrow\infty$ such that $\outCapc(B_k)\le \upCapc(B_k)\le \epsilon_k$ with
$$  B_k=\left\{\inf_{n_k\le m\le n_{k+1}}  \frac{S_m-\cSbep[S_m]}{a_m}< -1-\epsilon_k  \; \text{ or }\; \sup_{n_k\le m\le n_{k+1}}  \frac{S_m-\Sbep[S_m]}{a_m}>  1+\epsilon_k \right\}.$$
Notice $\sum_{k=1}^{\infty}\outCapc(B_k)<\infty$. By the countable sub-additivity of $\outCapc$, we have
 $$ \outCapc\left(B_k \; i.o. \right)=0, $$
which implies \eqref{eqthLILnonid1.4}.   $\Box$

 \bigskip
 {\bf Proof of Theorem \ref{thLILnonid2}}. \eqref{eqthLILnonid2.4}   follows from Theorem \ref{thLILnonid1}.
Now, we consider \eqref{eqthLILnonid2.5ad}.  Let $\lambda > 1$ be large enough. Let  $\{n_k=n_k(\lambda)\}\subset \mathbb N$ satisfy \eqref{eqproofLILnonid1.1}, and denote $I(k)=I(k,\lambda)=\{n_k+1,\ldots,n_{k+1}\}$. Then
\eqref{eqproofLILnonid1.2}  is satisfied.
Denote $b_j=\alpha \lambda^ 3 s_j/t_j$, where $0<\alpha  <1/(10\lambda^3)$ is a constant and to be specified. Redefine
$$Z_j=\big(-2b_{n_{k+1}}\big)\vee X_j\wedge \big(2b_{n_{k+1}}\big), \; j\in I(k), \; k\ge 0. $$
 Let $B_k^2=\sum_{j\in I(k)}\Sbep[X_j^2]$. Then $B_k^2=s_{n_{k+1}}^2-s_{n_k}^2\ge (1-1/\lambda)s_{n_{k+1}}^2$, $\log s_{n_k}^2\ge ck$. From the condition \eqref{eqthLILnonid2.2}, it follows that
\begin{equation} \label{eqproofLILnonid2.3} \sum_{k=1}^{\infty} (\log s_{n_k}^2)^{\delta-1}=\infty \text{ for all }\delta>0
\end{equation}
 by Lemma 2.3 of Wittmann (1987).   Similarly to \eqref{eqproofLILnonid1.5},  it follows from \eqref{eqthLILnonid2.1}  and \eqref{eqproofLILnonid1.2}  that
\begin{equation}\label{eqproofLILnonid2.5} \sum_{k=1}^{\infty}\left(\frac{   \Lambda_{n_k,n_{k+1}}(p,\lambda^3\alpha)}{a_{n_{k+1}}^p}\right)^{d+1}<\infty,
\end{equation}
where
$$\Lambda_{n_k,n_{k+1}}(p,\alpha)=\sum_{j\in I(k)}\Sbep[\big((|X_j|-\alpha s_{n_{k+1}}/t_{n_{k+1}})^+\big)^p]. $$
 Let
\begin{equation} \label{eqproofLILnonid2.6} \mathbb N_1=\left\{k\in \mathbb N; \frac{  \Lambda_{n_k, n_{k+1}}(p,\lambda^3\alpha)}{a_{n_{k+1}}^p}\le t_{n_{k+1}}^{-2p}\right\}.
\end{equation}
It follows from \eqref{eqproofLILnonid2.5}  that \eqref{eqproofLILnonid1.7}  holds.
By \eqref{eqthLILnonid2.3}, we have
\begin{equation} \label{eqproofLILnonid2.8} \frac{\sum_{j\in I(k)}|\Sbep[X_j]|}{a_{n_{k+1}}}+\frac{\sum_{j\in I(k)}|\cSbep[X_j]|}{a_{n_{k+1}}}\to 0.
\end{equation}
Note
\begin{align*}
& \Sbep|X_j-Z_j|\le \Sbep\big[\big((|X_j|-b_{n_{k+1}})^+\big)^p\big]b_{n_{k+1}}^{1-p},\\
& \Sbep|X_j-Z_j|^2+\Sbep|X_j^2-Z_j^2|\le 2\Sbep\big[\big((|X_j|-b_{n_{k+1}})^+\big)^p\big]b_{n_{k+1}}^{2-p}.
\end{align*}
Similarly to \eqref{eqproofLILnonid1.11}  and   \eqref{eqproofLILnonid1.12}, it follows that for $k\in \mathbb N_1$,
\begin{align}
 &\frac{\sum_{j\in I(k)}
\Sbep|X_j-Z_j|}{a_{n_{k+1}}}\le  \alpha^{1-p}t_{n_{k+1}}^{-2}\to 0,\label{eqproofLILnonid2.11}\\
 &\frac{\sum_{j\in I(k)}
 \big(\Sbep|X_j-Z_j|^2+\Sbep|X_j^2-Z_j^2|\big)}{s_{n_{k+1}}^2}   \le  2 \alpha^{2-p}t_{n_{k+1}}^{-2}\to 0.\label{eqproofLILnonid2.12}
\end{align}
Thus, similarly to \eqref{eqproofLILnonid1.13}, by Lemma \ref{lemExpIneq}  we have that for $k\in \mathbb N_1$ large enough,
\begin{align*}
&\Capc\Big(\sum_{j\in I(k)}|X_j-Z_j|\ge \epsilon a_{n_{k+1}}\Big)  \\
\le &   \exp\left\{-2 t_{n_{k+1}}^2\right\}+C\left( \frac{\Lambda_{n_k,n_{k+1}}(p,\lambda^3\alpha)\big]}{a_{n_{k+1}}^p}\right)^{d+1}+\sum_{j\in I(k)}\Capc\left(|X_j|\ge \epsilon^{\prime}a_j\right).
\end{align*}
It follows that
\begin{equation}\label{eqproofLILnonid2.14}
\sum_{k\in \mathbb N_1} \Capc\Big(\sum_{j\in I(k)}|X_j-Z_j|\ge \epsilon a_{n_{k+1}}\Big) <\infty\;\;\text{ for all } \epsilon>0.
\end{equation}
Next, we consider $Z_j$. Let $\widetilde{B}_k^2=\sum_{j\in I(k)}\Sbep[Z_j^2]$. It follows from \eqref{eqproofLILnonid1.2}  and \eqref{eqproofLILnonid2.12}  that
$$ \widetilde{B}_k^2\sim B_k^2 \ge (1-1/\lambda) s_{n_{k+1}}^2, \;\; k\in \mathbb N_1. $$
Without loss of generality, we assume that
$$  \frac{s_{n_{k+1}}}{\widetilde{B}_k}\le \frac{\lambda}{\lambda-1}, \;\; k\in \mathbb N_1. $$
It follows from \eqref{eqproofLILnonid2.8}  and \eqref{eqproofLILnonid2.11}  that
$$ \frac{\sum_{j\in I(k)}|\Sbep[Z_j]|}{a_{n_{k+1}}}+\frac{\sum_{j\in I(k)}|\cSbep[Z_j]|}{a_{n_{k+1}}}\to 0, \;\; \mathbb N_1\ni k\to \infty.
$$
Further,
$$ |Z_j|\le 2b_{n_{k+1}}=2\alpha\lambda^3\frac{s_{n_{k+1}}}{t_{n_{k+1}}}\le 2\alpha\lambda^3 \frac{\lambda}{\lambda-1}\frac{\widetilde{B}_k}{t_{n_{k+1}}}, \; j\in I(k) $$
for $k$ large enough. For every $\epsilon>0$, let $\gamma=\epsilon/2$ and $\pi(\gamma)$ be the constant defined as in  Lemma \eqref{lemExpIneq2}. Choose
$\alpha$ such that $2\alpha\frac{\lambda^4}{\lambda-1}<\pi(\gamma)$.
By Lemma \ref{lemExpIneq2}, we have that for $k\in \mathbb N_1$ large enough,
\begin{align}\label{eqproofLILnonid2.16}
& \Capc\Big(\sum_{j\in I(k)}Z_j\ge (1-\epsilon) (1-1/\lambda)  a_{n_{k+1}}\Big)\ge \Capc\Big(\sum_{j\in I(k)}Z_j\ge (1-2\epsilon)    \widetilde{B}_k t_{n_{k+1}}\Big)\nonumber\\
\ge &\exp\left\{ -\frac{(1-\epsilon)^2   t_{n_{k+1}}^2}{2}(1+\epsilon)\right\}
\ge   \exp\left\{ -\frac{(1-\epsilon^2) t_{n_{k+1}}^2}{2}\right\}.
\end{align}
It follows from \eqref{eqproofLILnonid2.14}, \eqref{eqproofLILnonid2.16}, \eqref{eqproofLILnonid2.3}  and \eqref{eqproofLILnonid1.7}   that
\begin{align*}
&\sum_{k\in \mathbb N_1}\Capc\Big(\sum_{j\in I(k)}X_j\ge (1-2\epsilon) (1-1/\lambda)  a_{n_{k+1}}\Big) \\
\ge &\sum_{k\in \mathbb N_1} \Capc\Big(\sum_{j\in I(k)}Z_j\ge (1-\epsilon) (1-1/\lambda)  a_{n_{k+1}}\Big)-C\\
\ge & \sum_k\exp\left\{ -\frac{(1-\epsilon^2) t_{n_{k+1}}^2}{2}\right\}-\sum_{k\in \mathbb N\setminus \mathbb N_1}
\exp\left\{ -\frac{(1-\epsilon^2) t_{n_{k+1}}^2}{2}\right\}-C\\
\ge& \sum_k  (\log s_{n_{k+1}}^2)^{\epsilon^2-1}-c\sum_{k\in \mathbb N\setminus \mathbb N_1}t_{n_{k+1}}^{-2p(d+1)}-C=\infty.
\end{align*}
Hence
\begin{equation}\label{eqproofLILnonid2.17}\sum_{k=1}^{\infty}\Capc\Big(\sum_{j\in I(k)}X_j\ge (1-2\epsilon) (1-1/\lambda)  a_{n_{k+1}}\Big)=\infty.
\end{equation}
Noting the  independence of $\sum_{j\in I(k)}X_j$, $k=1,2, \ldots,$  by   Lemma \ref{lemBC} (ii)  it follows that
\begin{equation}\label{eqproofLILnonid2.22}\Capc\Big(\max_{ K\le k\le M} \frac{\sum_{j\in I(k)}X_j}{a_{n_{k+1}}}\ge (1-3\epsilon) (1-1/\lambda)  \Big)
\to 1\;\text{ as } M\to \infty  \; \text{ for all } K.
\end{equation}
On the other hand,
\begin{equation}\label{eqproofLILnonid2.23} \frac{|S_{n_k}|}{a_{n_{k+1}}}\le  \frac{|S_{n_k}|}{a_{n_k}} \cdot  \frac{a_{n_k}}{a_{n_{k+1}}}\le \frac{|S_{n_k}|}{a_{n_k}}\frac{1}{\lambda}
\end{equation}
for $k$ large enough.
It follows that for $K$ large enough,
\begin{align*}
&\max_N\Capc\Big(\max_{K\le   k\le N}\max_{n\in I_k} \frac{S_m}{a_m}\ge (1-3\epsilon) (1-1/\lambda)-(1+\epsilon)/\lambda  \Big) \\
\ge &\Capc\Big(\max_{K\le k\le M} \frac{S_{n_{k+1}}}{a_{n_{k+1}}}\ge (1-3\epsilon) (1-1/\lambda)-(1+\epsilon)/\lambda  \Big) \\
\ge & \Capc\Big(\max_{K\le k\le M} \frac{\sum_{j\in I(k)}X_j}{a_{n_{k+1}}}\ge (1-3\epsilon) (1-1/\lambda)  \Big)
-\Capc\Big(\max_{K\le k\le M} \frac{|S_{n_k}|}{a_{n_k}}\ge  1+\epsilon   \Big)\\
&\to    1 \text{ as } M\to \infty \text{ and then } K\to \infty,
\end{align*}
by \eqref{eqthLILnonid2.4}   and \eqref{eqproofLILnonid2.22}. By the arbitrariness of    $\epsilon> 0$  being small enough and $\lambda>1$ being large enough,  we obtain
$$ \lim_{n\to \infty}\max_{N} \Capc\Big(\max_{n\le m\le N} \frac{S_m}{a_m}\ge  1- \epsilon   \Big)=1\;\text{ for all } \epsilon>0. $$
 For $-X_j$, we have the same conclusion. \eqref{eqthLILnonid2.5ad}  is proved.

At last, as in Theorem \ref{thLILnonid1}, by the countable sub-additivity of $\outCapc$, \eqref{eqthLILnonid2.4}  implies \eqref{eqthLILnonid2.6}.
   $\Box$

\bigskip
For proving Theorems \ref{thLILiid1} and \ref{thLILiid2} for independent and identically distributed random variables, we need more two  lemmas.
\begin{lemma} \label{lem2} Suppose $X\in \mathscr{H}$.
\begin{description}
  \item[\rm (i)]
 For any $\delta>0$,
$$ \sum_{n=1}^{\infty} \Capc\big(|X|\ge \delta \sqrt{n\log\log n} \big)<\infty \;\; \Longleftrightarrow C_{\Capc}\left[\frac{X^2}{\log\log|X|}\right]<\infty.
$$
 \item[\rm (ii)]
   If $C_{\Capc}\left[\frac{X^2}{\log\log|X|}\right]<\infty$, then for any $\delta>0$ and $p>2$,
$$ \sum_{n=1}^{\infty} \frac{\Sbep\big[\big(|X|\wedge (\delta \sqrt{n\log\log n})\big)^p\big]}{(n\log\log n)^{p/2}}<\infty. $$
 \item[\rm (iii)] $C_{\Capc}\left[\frac{X^2}{\log\log|X|}\right]<\infty$, then for any $\delta>0$,
 $$ \Sbep[X^2\wedge (2\delta n\log\log n)]=o(\log \log n) $$
 and $$ \breve{\mathbb E}[(|X|-\delta \sqrt{2  n\log\log n})^+]=o(\sqrt{\log\log n/n}).\;\; \Box $$
\end{description}
\end{lemma}
 {\bf Proof.}  The proof of (i) and (ii)
  can be found in Zhang (2016). For (iii), we denote $d_n =\sqrt{2n\log\log n}$. Let $ f(x)$ be the inverse function of $g(x)=\sqrt{2x\log\log x}$ $(x>0)$.   Then $ c x^2/\log\log |x|\le f(|x|)\le Cx^2/\log\log |x|$. It follows that
 $\int_0^{\infty}\Capc\big(f(|X|\ge y)dy<\infty$. Hence
 \begin{align*}
 & \frac{\Sbep\big[ (|X|\wedge \delta d_n)^2 \big]}{\log\log n} \le   \frac{C_{\Capc}\big( (|X|\wedge \delta d_n)^2 \big)}{\log\log n}
 = \frac{1}{\log\log n} \int_0^{(\delta d_n)^2}\Capc\left(|X|^2>x\right)dx \\
= &\frac{1}{\log\log n}2\int_0^{\delta d_n}x\Capc\left(|X|>x\right)dx  = \frac{2}{\log\log n} \int_0^{g(\delta d_n)}g(y)\Capc\left(|X|>g(y)\right)dg(y) \\
\le  & 4\int_0^{C_{\delta} n}  \frac{\log\log y}{\log\log n}  \Capc\left(f(|X|)>y\right)dy\to 0
\end{align*}
 and
 \begin{align*}
 &\Sbep\big[ (|X|- \delta d_n)^+ \big]\le C_{\Capc}\big((|X|- \delta d_n)^+\big)\le \int_{\delta d_n}^{\infty}\Capc(|X|\ge x)dx \\
 =& \int_{g(\delta d_n)}^{\infty}\Capc(|X|\ge g(y))d g(y)
 \le   2\sqrt{2}\int_{c_{\delta}n}^{\infty} \sqrt{\log \log y/y} \Capc(f(|X|)\ge y)d y\\
 \le & 2\sqrt{2}\sqrt{\log\log n/n}\int_{c_{\delta}n}^{\infty}\Capc(f(|X|)\ge y)d y=o(\sqrt{\log\log n/n}). \; \Box
 \end{align*}

\begin{lemma} \label{lem3}
   Let  $\{Y_n; n\ge 1\}$ be a sequence of   independent and identically distributed random variables in the sub-linear expectation space $(\Omega,\mathscr{H},\Sbep)$ with  $C_{\Capc}\left[\frac{Y_1^2}{\log\log|Y_1|}\right]<\infty$. Then
   \begin{align}
   & \Capc\left(\frac{\sum_{i=1}^n(Y_i-\breve{\mathbb E}[Y_1])}{\sqrt{2n\log\log n}}\ge \epsilon\right)\to 0 \; \text{ for all } \epsilon>0, \label{eqlem3.1}\\
   &\cCapc\left(\frac{\sum_{i=1}^n(-Y_i+\breve{\mathbb E}[Y_1])}{\sqrt{2n\log\log n}}\ge \epsilon\right)\to 0 \; \text{ for all } \epsilon>0. \label{eqlem3.2}
   \end{align}
\end{lemma}
{\bf Proof}. For a random variable $Y$, we denote $Y^{(c)}=(-c)\vee Y\wedge c$. Denote $d_n=\sqrt{2\log\log n}$. Then by applying \eqref{eqExpIneq.4}  and \eqref{eqExpIneq.7}  with $p=2$, we obtain
\begin{align*}
\Capc\left(\sum_{i=1}^n(Y_i^{(d_n)}-n\Sbep[Y_1^{(d_n)}])\ge \epsilon d_n\right)\le C \frac{n \Sbep[(|Y_1|\wedge d_n)^2]}{\epsilon^2d_n^2}\to 0
\end{align*}
and
$$ \frac{n\big|\breve{\mathbb E}[Y_1]-\Sbep[Y_1^{(d_n)}]\big|}{d_n}=\frac{n \breve{\mathbb E}\big[(|Y_1|-d_n)^+\big]}{d_n}\to 0 $$
by Lemma \ref{lem2} (iii).  On the other hand,
$$ \Capc\left( Y_i^{(d_n)} \ne Y_i\text{ for some } i=1,\ldots, n\right)\le n\Capc\left(|Y_1|\ge d_n\right)\to 0 $$
by Lemma \ref{lem2} (i). Therefore, \eqref{eqlem3.1}  holds. The proof of \eqref{eqlem3.2}  is similar. $\Box$

\bigskip
{\bf Proof of Theorems \ref{thLILiid1} and \ref{thLILiid2}}.   If $\overline{\sigma}^2=0$, then $|\breve{\mathbb E}[\pm Y_1]|
\le \breve{\mathbb E}[|Y_1|]\le (\breve{\mathbb E}[Y_1^2])^{1/2}=0$, $\overline{\sigma}_1^2=\overline{\sigma}_2^2=\overline{\sigma}^2=0$, and for any $\epsilon>0$, $\Capc(|Y_1|\ge \epsilon)
\le \Sbep[Y_1^2\wedge \epsilon^2]/\epsilon^2=0$. By the countable sub-additivity of $\outCapc$, it follows that $\outCapc(|Y_n|>0)\le \sum_{j=1}^{\infty}\outCapc(|Y_n|>1/j)\le \sum_{j=1}^{\infty} l^2\Sbep[Y_1^2\wedge j^2]=0$. Hence
$\outCapc(|Y_n|\ne 0\text{ for some } n)=0$. And then, \eqref{eqLILmomentcondition1}, \eqref{eqthLILiid1.3}  in Theorem \ref{thLILiid1} and \eqref{eqLILmeanzerocondition} -\eqref{eqthLILiid2.4}  in Theorem \ref{thLILiid2} hold automatically. Therefore, without loss of generality, we assume $0<\overline{\sigma}\le \infty$.

We first  suppose that \eqref{eqLILmomentcondition1}   is  satisfied.
 Let $d_n=\sqrt{2n\log\log n}$ and $X_n=(-d_n)\vee Y_n\wedge d_n$.
Denote $S_n=\sum_{i=1}^n X_i$, $s_n^2=\sum_{i=1}^n \Sbep[X_i^2]$,  $t_n=\sqrt{2\log\log s_n}$ and $a_n=t_ns_n$.  Then $n\le Cs_n^2$, $d_n\le Ca_n$. Notice that
\begin{align*}
&\sum_{n=1}^{\infty} \outCapc(Y_n\ne X_n)=\sum_{n=1}^{\infty}\outCapc(|Y_n|>d_n)\le \sum_{n=1}^{\infty}\Capc(|Y_1|>d_n/2)<\infty,  \\
   & |\breve{\mathbb E}[Y_n]-\Sbep[X_n]|+ |\breve{\mathbb E}[-Y_n]-\Sbep[-X_n]|\\
   \le &  2 \breve{\mathbb E}[|Y_n-X_n|]|\le 2\breve{\mathbb E}[(|Y_1|-d_n)^+]
=o\left(\sqrt{\log\log n}/\sqrt{n}\right)
\end{align*}
by Lemma \ref{lem2}.
It follows that
\begin{align}
& \frac{\sum_{i=1}^n \big|\breve{\mathbb E}[Y_i]-\Sbep[X_i]\big|+ \big|\breve{\mathbb E}[-Y_i]-\Sbep[-X_i]\big|}{d_n}\to 0,\label{eqproofLILiid2.5}\\
  \sum_{n=1}^{\infty} & \frac{\Sbep[|X_n|^p]}{a_n^p} \le c\sum_{n=1}^{\infty}\frac{\Sbep[|Y_1|^p\wedge d_n^p ]}{d_n^p}<\infty
 \;\;(\text{by Lemma \ref{lem2}}),\label{eqproofLILiid2.6}\\
 &\outCapc\left(\lim_{n\to\infty}\frac{\sum_{i=1}^n(X_i-Y_i)}{d_n}\ne 0\right)
 \le  \outCapc\left(Y_n\ne X_n \;\; i.o.\right)=0.
 \label{eqproofLILiid2.8}
\end{align}
Moreover,
$$\frac{\Sbep[X_{n+1}^2]}{s_n^2}\le \frac{C_0d_nC_{\Capc}(|Y_1|)}{n}\to 0, $$
which implies
\begin{equation}\label{eqproofLILiid2.9} s_n^2\to \infty\; \text{ and } \; \frac{s_{n+1}^2}{s_n^2}\to 1, \;\; \frac{a_{n+1}}{a_n}\to 1 .
\end{equation}

We first show that
\begin{equation}\label{eqproofLILiid2.10}
\outCapc\left(\limsup_{n\to \infty}\frac{\sum_{i=1}^n (Y_i-\breve{\mathbb E}[Y_i])}{\sqrt{2n\log\log n}}> \overline{\sigma}_2 \right)=0.
\end{equation}
Without loss of generality, assume $\breve{\mathbb E}[Y_1]=0$. It follows from  \eqref{eqproofLILiid2.6}  and \eqref{eqproofLILiid2.9}  that the conditions  \eqref{eqremkLIL3.1}  in Remark \ref{remark3}  and \eqref{eqthLILnonid2.0}  in Remark \ref{remark1}    are satisfied. By Theorem \ref{thLILnonid1} and Remarks \ref{remark1} and \ref{remark3}, we have
$$
\outCapc\left(\limsup_{n\to \infty}\frac{\sum_{i=1}^n (X_i-\breve{\mathbb E}[X_i])}{s_nt_n}> 1 \right)=0,
$$
which, together with \eqref{eqproofLILiid2.5}  and \eqref{eqproofLILiid2.8}, implies
$$
\outCapc\left(\limsup_{n\to \infty}\frac{\sum_{i=1}^n (Y_i-\breve{\mathbb E}[Y_i])}{s_nt_n}> 1 \right)=0.
$$
Note
\begin{equation}\label{eqlimvar} \lim_{n\to \infty}\frac{s_n^2}{n}=\overline{\sigma}^2\;\; (\text{finite and infinite}).
\end{equation}
\eqref{eqproofLILiid2.10}  is proved. For $-X_j$ we reach a similar conclusion, and therefore, Theorem \ref{thLILiid1} is proved.

 Next, we turn to the proof of Theorem \ref{thLILiid2}. For the part (a), besides \eqref{eqLILmomentcondition1}  we  further assume \eqref{eqLILmeanzerocondition}, i.e., $\breve{\mathbb E}[Y_n]=\breve{\mathbb E}[-Y_n]=0$.  It follows from \eqref{eqproofLILiid2.5}   that the condition  \eqref{eqthLILnonid2.3}  in Theorem \ref{thLILnonid2} is also satisfied.
Then \eqref{eqthLILiid2.3}     and \eqref{eqthLILiid2.4}  are implied by \eqref{eqthLILnonid2.6}  and \eqref{eqthLILnonid2.5ad}.

Now, we consider the part (b).   Suppose
$$C_{\Capc}\left[\frac{Y_1^2}{\log\log|Y_1|}\right]=\infty. $$
By Lemma  \ref{lem2},
\begin{align*}
  \sum_{n=1}^{\infty} \Capc\big(|Y_n|\ge 2Md_n \big)
\ge  \sum_{n=1}^{\infty} \Capc\big(|Y_1|\ge 3Md_n \big)=\infty\;\; \text{ for all } M>0.
\end{align*}
It follows that there exists a sequence $M_n\nearrow \infty$ such that
\begin{align*}
  \sum_{n=1}^{\infty} \Capc\big(|Y_n|\ge 2M_nd_n \big) =\infty.
\end{align*}
By Lemma \ref{lemBC} (ii), it follows that
$$ \lim_{N\to \infty} \Capc\Big(\max_{n\le m\le N } \frac{|Y_m|}{d_m}\ge M\Big)=1\; \text{ for all } M>0,n\ge 1. $$
Notice $|Y_m|\le |\sum_{i=1}^m Y_i|+|\sum_{i=1}^{m-1}Y_i|$. It follows that
$$   \lim_{N\to \infty} \Capc\Big(\max_{n\le m\le N}\frac{|\sum_{i=1}^mY_i|}{d_m}\ge M\Big)=1\; \text{ for all }M>0, n\ge 1, $$
which contradicts \eqref{eqthLILiid2.6}. It follows that \eqref{eqLILmomentcondition1}  holds, and then there exist $0<\tau<1$, $M>1$ and $n_0\ge 1$ such that
\begin{equation}\label{eqproofLILiid2.14}  \lim_{N\to \infty}\Capc\Big(\max_{n\le m\le N} \frac{|\sum_{i=1}^mY_i|}{d_m}\ge M\Big)<\tau<1\; \text{ for all } n\ge n_0.
\end{equation}
Under \eqref{eqLILmomentcondition1}, $\breve{\mathbb E}[Y_1]$ and  $\breve{\mathbb E}[-Y_1]$ exist and are finite.     On the other hand, by Lemma \ref{lem3},
$$ \Capc\left(\frac{\sum_{i=1}^n Y_i-n \breve{\mathbb E}[Y_1]}{d_n}\ge -\epsilon \right)=1-\cCapc\left(\frac{\sum_{i=1}^n (-Y_i+ \breve{\mathbb E}[Y_1])}{d_n}>  \epsilon \right)\to 1 \;\;\text{ for all } \epsilon>0. $$
It follows that
$$\liminf_{n\to \infty}\lim_{N\to \infty}\Capc\left(\frac{\sum_{i=1}^n Y_i-n \breve{\mathbb E}[Y_1]}{d_n}> -\epsilon,\;\; \max_{N\ge m\ge n} \frac{|\sum_{i=1}^mY_i|}{d_m}< M\right)
\ge 1-\tau>0. $$
Therefore,
$$ \breve{\mathbb E}[Y_1]\le \frac{(M+\epsilon)d_n}{n}\to 0. $$
Similarly, $\breve{\mathbb E}[-Y_1]\le 0$. From the fact that $\breve{\mathbb E}[Y_1]+\breve{\mathbb E}[-Y_1]\ge 0$, it follows that \eqref{eqLILmeanzerocondition}  holds.

 Under \eqref{eqLILmomentcondition1}  and \eqref{eqLILmeanzerocondition}, we  still have  \eqref{eqthLILiid2.4}  which contradicts \eqref{eqproofLILiid2.14}   if $\overline{\sigma}=\infty$.
 Hence \eqref{eqthLILiid2.8}  holds.
  $\Box$

\bigskip
{\bf Proof of Theorem \ref{thLIL5}}.
By \eqref{eqthLILiid2.3}, it is sufficient to show that
$$
  \outcCapc\left(\limsup_{n\to \infty}\frac{T_n}{d_n}\ge \underline{\sigma}\right)=1.
$$
Notice that $\outcCapc$ has the property \eqref{CapcC}  since $\outCapc$ is countably sub-additive. It is sufficient to show that
$$
 \outcCapc\left(\limsup_{n\to \infty}\frac{T_n}{d_n}\ge \underline{\sigma}-\epsilon\right)=1,\;\; \forall \epsilon>0,
$$
because
$$\outcCapc\left(\limsup_{n\to \infty}\frac{T_n}{d_n}\ge \underline{\sigma}\right)
=\outcCapc\left(\bigcap_{l=1}^{\infty}\Big\{\limsup_{n\to \infty}\frac{T_n}{d_n}\ge \underline{\sigma}-\frac{1}{l}\Big\}\right). $$
When $\underline{\sigma}=0$, the conclusion  is trivial because   $\outcCapc$ has the property \eqref{CapcC}  and
$$  \outcCapc\left(\bigcup_{n=m}^{\infty} \frac{T_n}{d_n}\ge  -\epsilon\right)
\ge \outcCapc\left(  \frac{T_n}{d_n}\ge  -\epsilon\right)\ge 1-\upCapc\left(  \frac{-T_n}{d_n}> \epsilon\right)
\to 1, \text{ for all } m, $$
by Lemma \ref{lem3}.

Suppose $\underline{\sigma}>0$.
Let $\lambda > 1$. Denote $n_k=[\lambda^k]$ and $I(k)=\{n_k+1,\ldots,n_{k+1}\}$. Then
 $n_k/n_{k+1}\to  {1}/{\lambda} $, $d_{n_k}/d_{n_{k+1}}\to 1/\sqrt{\lambda}$. Notice \eqref{eqthLILiid2.3},
$$ \frac{T_{n_{k+1}}}{d_{n_{k+1}}}=\frac{T_{n_{k+1}}-T_{n_k}}{\sqrt{2(n_{k+1}-n_k)\log\log n_{k+1}}}\sqrt{1-\frac{n_k}{n_{k+1}}}+
 \frac{T_{n_k}}{d_{n_k}}\frac{d_{n_k}}{d_{n_{k+1}}},
$$
  and that $\lambda>1$ can be chosen large enough. It is sufficient to show that,
\begin{equation}\label{eqproofthLIL5.6}   \outcCapc\left( \limsup_{k\to \infty} \frac{T_{n_{k+1}}-T_{n_k}}{\sqrt{2(n_{k+1}-n_k)\log\log n_{k+1}}}\ge \underline{\sigma}  -\epsilon \right)=1, \;\; \forall \epsilon>0.
\end{equation}

Denote $t_j=\sqrt{2\log\log j}$ and $b_j=\alpha_j   \sqrt{j}/\sqrt{2\log\log j}$, where $\alpha_j \to 0$  is  specified such that $\alpha_j\to 0$  and $\alpha_j^{1-p}t_j^{-2}\to 0$.  Define
$$Z_j=\big(-2b_{n_{k+1}}\big)\vee Y_j\wedge  2b_{n_{k+1}}, \; j\in I(k), \; k\ge 0. $$
By Lemma \ref{lem2} (ii), we have
\begin{equation} \label{eqproofthLIL5.7}  \sum_{k=1}^{\infty} \frac{   \Lambda_{n_k,n_{k+1}}(p)}{d_{n_{k+1}}^p}<\infty,
\end{equation}
where
$$\Lambda_{n_k,n_{k+1}}(p)=\sum_{j\in I(k)}\Sbep[\big((|Y_j|\wedge d_{n_{k+1}}\big)^p]. $$
 Let
\begin{equation} \label{eqproofthLIL5.8}   \mathbb N_1=\left\{k\in \mathbb N; \frac{  \Lambda_{n_k, n_{k+1}}(p)}{d_{n_{k+1}}^p}\le t_{n_{k+1}}^{-2p}\right\}.
\end{equation}
Similar to \eqref{eqproofLILnonid1.11}  and \eqref{eqproofLILnonid1.12}, we have for $k\in \mathbb N_1$,
\begin{align}
 &\frac{\sum_{j\in I(k)}
\Sbep[ |(-d_{n_{k+1}})\vee Y_j\wedge d_{n_{k+1}}-Z_j|]}{d_{n_{k+1}}}\nonumber \\
\le & \frac{\Lambda_{n_k,n_{k+1}}(p)}{d_{n_{k+1}}^p}  \alpha_{n_{k+1}}^{1-p}t_{n_{k+1}}^{2p-2} \le  \alpha_{n_{k+1}}^{1-p}t_{n_{k+1}}^{-2}\to 0,\label{eqproofthLIL5.11}\\
 &\frac{\sum_{j\in I(k)}
  \Sbep[((-d_{n_{k+1}})\vee Y_j \wedge a_{n_{k+1}}-Z_j)^2] }{ n_{k+1} } \nonumber \\
   \le & \frac{\Lambda_{n_k,n_{k+1}}(p)}{d_{n_{k+1}}^p}  \alpha_{n_{k+1}}^{2-p}t_{n_{k+1}}^{2p-2} \le    \alpha_{n_{k+1}}^{2-p}t_{n_{k+1}}^{-2}\to 0,\label{eqproofthLIL5.12}
\end{align}
by noting $\alpha_j\to 0$ such that $\alpha_j^{1-p}t_j^{-2}\to 0$. Similar to \eqref{eqproofLILnonid1.13},
we have that for $k\in \mathbb N_1$ large enough,
\begin{align*}
&\Capc\Big(\sum_{j\in I(k)}|Y_j-Z_j|\ge \epsilon d_{n_{k+1}}\Big) \\
\le & \Capc\Big(\sum_{j\in I(k)}|(-d_{n_{k+1}})\vee Y_j\wedge d_{n_{k+1}}-Z_j|\ge \epsilon d_{n_{k+1}}\Big)+\sum_{j\in I(k)}\Capc(|X_j|>d_{n_{k+1}})  \\
\le &  \exp\left\{-2 t_{n_{k+1}}^2\right\}
    +C \frac{\Lambda_{n_k,n_{k+1}}(p)}{d_{n_{k+1}}^p} +2\sum_{j\in I(k)}\Capc(X_j>\epsilon^{\prime} d_j).
\end{align*}
It follows that
\begin{equation}\label{eqproofthLIL5.13}
\sum_{k\in \mathbb N_1 } \Capc\Big(\sum_{j\in I(k)}|Y_j-Z_j|\ge \epsilon d_{n_{k+1}}\Big)<\infty.
\end{equation}

Next, we apply Lemma \ref{lemLower} to the array $\{Z_j;j\in I(k)\}$ of independent and identically random variables, $k\in \mathbb N_1$.
By \eqref{eqproofthLIL5.12}, we have $\Sbep[Z_j^2]\sim \Sbep[Y_1^2\wedge d_{n_{k+1}}^2]\to \overline{\sigma}^2$ and $\cSbep[Z_j^2]\sim \cSbep[Y_1^2\wedge d_{n_{k+1}}^2]\to \underline{\sigma}^2$.  By \eqref{eqproofthLIL5.11}, Lemma \ref{lem2} (iii) and the fact that $\breve{\mathbb E}[Y_1]=\breve{\mathbb E}[-Y_1]=0$, we have
$$ \frac{\sum_{j\in I(k)}
(|\Sbep[Z_j]|+|\cSbep[Z_j]|)}{\sqrt{n_{k+1}-n_k}t_{n_{k+1}}} \to 0. $$
Note
$$|Z_j|\le 2 b_{n_{k+1}}=o\Big(\sqrt{n_{k+1}-n_k}/t_{n_{k+1}}\Big),\;\; j\in I(k). $$
Applying Lemma \ref{lemLower} with $k_n=n_{k+1}-n_k$ and $x_n=t_{n_{k+1}}$ yields
\begin{align}\label{eqproofthLIL5.14}
& \cCapc\left( \frac{\sum_{j\in I(k)}Z_j}{\sqrt{2(n_{k+1}-n_k)\log\log n_{k+1}}}\ge \underline{\sigma}(1-\epsilon)\right)\nonumber \\
\ge &
\exp\left\{- (1-\epsilon)\log\log n_{k+1}\right\}\ge c k^{-(1-\epsilon)},
\end{align}
for $k\in \mathbb N_1$ large enough. Notice that
\begin{equation}\label{eqproofthLIL5.15}\sum_{k\not\in \mathbb N_1}\exp\left\{- (1-\epsilon)\log\log n_{k+1}\right\}\le C\sum_{k\not\in \mathbb N_1} t_{n_{k+1}}^{-2p}<\infty,
\end{equation}
by \eqref{eqproofthLIL5.7}  and \eqref{eqproofthLIL5.8}.
From \eqref{eqproofthLIL5.13}, \eqref{eqproofthLIL5.14}  and \eqref{eqproofthLIL5.15}, we conclude that
$$\sum_{k\in \mathbb N_1} \cCapc\left( \frac{\sum_{j\in I(k)}Y_j}{\sqrt{2(n_{k+1}-n_k)\log\log n_{k+1}}}\ge \underline{\sigma}(1-\epsilon)\right)=\infty. $$
  Hence,
$$\sum_{k=1}^{\infty} \cCapc\left( \frac{T_{n_{k+1}}-T_{n_k}}{\sqrt{2(n_{k+1}-n_k)\log\log n_{k+1}}}\ge \underline{\sigma}(1-\epsilon)\right)=\infty, \;\; \forall \epsilon>0. $$
The above equation also holds for $\lowCapc$ by \eqref{eqV-V}.
Notice the independence of random variables and that $\outCapc$ is countably sub-additive with $\outCapc\le \upCapc$. By  Lemma \ref{lemBC} (iii), we have that
$$\outcCapc\left( \frac{T_{n_k}-T_{n_{k-1}}}{\sqrt{2(n_k-n_{k-1})\log\log n_k}}\ge \underline{\sigma}(1-\epsilon)-\epsilon  \right)=1. $$
\eqref{eqproofthLIL5.6}  is proved.
$\Box$

\subsection{Proof of the exact lower bounds}

Now, we consider the results on the exact lower bounds. We first prove Theorems \ref{thLILnew1} and \ref{thLILnew2} and then Theorems \ref{thLILnew3} and \ref{thLILnew4}.

{\bf Proof of  Theorem \ref{thLILnew1}}.   We will  show that  \eqref{eqthLILnonid2.4} and \eqref{eqthLILnonid2.5ad} imply \eqref{eqthLILnew1.1} and  \eqref{eqthLILnew1.2} although $V$ may be   not continuous.  That \eqref{eqthLILnonid2.4} implies \eqref{eqthLILnew1.1} is shown in the proof of Theorem \ref{thLILnonid1}. Now, consider \eqref{eqthLILnew1.2}.  Let  $\epsilon_k=1/2^k$.
By \eqref{eqthLILnonid2.4} and \eqref{eqthLILnonid2.5ad}, there exist sequences $n_k\nearrow \infty$ and $m_k\nearrow \infty$ with $n_k<m_k<n_{k+1}$ such that
\begin{equation}\label{eqproofthLILnew1.3}
\max_N \upCapc\left(\max_{n_k\le m\le N}\frac{|\tilde S_m|}{a_m}\ge 1+\epsilon_k/2\right)\le \epsilon_k/2
 \end{equation}
and
\begin{equation}\label{eqproofthLILnew1.4}
 \begin{aligned}
 & \upCapc\left(\max_{n_k\le m\le m_k}\frac{  S_m}{a_m}\ge 1-\epsilon_k/2\right)\ge 1-\epsilon_k/2 \text{ and } \\
 & \upCapc\left(\max_{n_k\le m\le m_k}\frac{-  S_m}{a_m}\ge 1-\epsilon_k/2\right)\ge 1-\epsilon_k/2.
 \end{aligned}
 \end{equation}
Without loss of generality, we assume $a_{m_{k-1}}/a_{n_k}\le \epsilon_k/4$. From \eqref{eqproofthLILnew1.3} and \eqref{eqproofthLILnew1.4}, it follows that
\begin{align}\label{eqproofthLILnew1.9}
&\upCapc\left(\max_{n_k\le m\le m_k}\frac{ S_m- S_{m_{k-1}}}{a_m}\ge 1-\epsilon_k\right)\nonumber\\
\ge &\upCapc\left(\max_{n_k\le m\le m_k}\frac{\tilde S_m}{a_m}\ge 1-\epsilon_k/2\right)
-\upCapc \left(\frac{|\tilde  S_{m_{k-1}}|}{a_{n_k}}\ge \epsilon_k/2\right)\nonumber \\
\ge &1-\epsilon_k/2-\max_N\upCapc\left(\max_{n_{k-1}\le m\le N}\frac{|\tilde S_m|}{a_m}\ge 2\right)\nonumber \\
\ge &1-\epsilon_k/2-\epsilon_{k-1}/2\ge 1-\epsilon_{k-1}.
\end{align}
Similarly,
$$
\upCapc\left(\max_{n_k\le m\le m_k}\frac{ -(S_m- S_{m_{k-1}})}{a_m}\ge 1-\epsilon_k\right)\ge 1-\epsilon_{k-1}.
$$
Let
$$\widetilde{B}_{k,+}(\epsilon)=\left\{\bm x: \max_{n_{2k}\le m\le m_{2k}}\frac{ s_m- s_{m_{2k-1}}}{a_m}\ge 1-\epsilon\right\},$$
$$\widetilde{B}_{k,-}(\epsilon)= \left\{\bm x:\max_{n_{2k+1}\le m\le m_{2k+1}}\frac{ -(s_m- s_{m_{2k}})}{a_m}\ge 1-\epsilon\right\},  $$
where $s_m=\sum_{i=1}^m x_i. $
 Choose Lipschitz functions $f_{k,+}$ and $f_{k,-}$  such that
$$ I_{\widetilde{B}_{k,+}(2\epsilon_{2k})}\ge f_{k,+}\left(\max_{n_{2k}\le m\le m_{2k}}\frac{ s_m- s_{m_{2k-1}}}{a_m}\right)\ge  I_{\widetilde{B}_{k,+}(\epsilon_{2k})},$$
$$ I_{\widetilde{B}_{k,-}(2\epsilon_{2k+1})}\ge f_{k,-}\left(\max_{n_{2k+1}\le m\le m_{2k+1}}\frac{ -(s_m- s_{m_{2k}})}{a_m}\right)\ge  I_{\widetilde{B}_{k,-}(\epsilon_{2k+1})}.$$
Recall  that $\mathscr{P}$ is the family defined as in Proposition \ref{prop2}. Then
 by noting the independence, we have
\begin{align*} & \Capc^{\mathscr{P}}\left(\bm X\in \bigcap_{k=\ell}^{N} \Big(\widetilde{B}_{k,+}(2\epsilon_{2k})\bigcap \widetilde{B}_{k,-}(2\epsilon_{2k+1})\Big)\right)\\
\ge & \Sbep\left[  \prod_{k=\ell}^{N}f_{k,+}\left(\max_{n_{2k}\le m\le m_{2k}}\frac{ S_m- S_{m_{2k-1}}}{a_m}\right) f_{k,-}\left(\max_{n_{2k+1}\le m\le m_{2k+1}}\frac{ -(S_m- S_{m_{2k}})}{a_m}\right)\right]\\
=&  \prod_{k=\ell}^{N}\Sbep\left[  f_{k,+}\left(\max_{n_{2k}\le m\le m_{2k}}\frac{ S_m- S_{m_{2k-1}}}{a_m}\right)\right]\Sbep\left[ f_{k,-}\left(\max_{n_{2k+1}\le m\le m_{2k+1}}\frac{ -(S_m- S_{m_{2k}})}{a_m}\right)\right] \\
\ge &  \prod_{k=\ell}^{N}\left[\upCapc(\bm X\in \widetilde{B}_{k,+}(\epsilon_{2k}))\cdot \upCapc(\bm X\in \widetilde{B}_{k,-}(\epsilon_{2k+1}))\right] \\
\ge &\prod_{k=\ell}^{\infty}\left[\upCapc(\bm X\in \widetilde{B}_{k,+}(\epsilon_{2k}))\cdot \upCapc(\bm X\in \widetilde{B}_{k,-}(\epsilon_{2k+1}))\right] \ge \prod_{k=2\ell}^{\infty}(1-\epsilon_{k-1}).
\end{align*}
 Let $\widetilde{A}_{\ell}=\bigcap_{k=\ell}^{\infty} \Big(\widetilde{B}_{k,+}(2\epsilon_{2k})\bigcap \widetilde{B}_{k,-}(2\epsilon_{2k+1})\Big)$. Then $\bigcap_{k=\ell}^N \Big(\widetilde{B}_{k,+}(2\epsilon_{2k})\bigcap \widetilde{B}_{k,-}(2\epsilon_{2k+1})\Big)$ is a closed subset of $\mathbb R^{\infty}$.  Notice that each $X_i$ is tight since $\Sbep[X_i^2]$ is finite. By \eqref{eqClose}, it follows that
\begin{align}\label{eqproofthLILnew1.10}   &   \Capc^{\mathscr{P}}\left(\bm X\in \widetilde{A}_{\ell} \right)=\lim_{N\to \infty}\Capc^{\mathscr{P}}\left(\bm X\in \bigcap_{k=\ell}^{N} \Big(\widetilde{B}_{k,+}(2\epsilon_{2k})\bigcap \widetilde{B}_{k,-}(2\epsilon_{2k+1})\Big)\right)\nonumber \\
 \ge & \prod_{k=\ell}^{\infty}\left[\Capc(\bm X\in \widetilde{B}_{k,+}(\epsilon_{2k}))\cdot \Capc(\bm X\in \widetilde{B}_{k,-}(\epsilon_{2k+1}))\right]
 \ge   \prod_{k=2\ell}^{\infty}(1-\epsilon_{k-1})\to 1.
\end{align}
Hence $  \Capc^{\mathscr{P}}(\bigcup_{\ell=1}^{\infty}\{\bm X\in \widetilde{A}_{\ell}\})=1$. On the event $\bigcup_{\ell=1}^{\infty}\{\bm X\in \widetilde{A}_{\ell}\}$ and $\{\limsup_{n\to \infty}\frac{|S_n|}{a_n}\le 1\}$,
 \begin{align*}
1\ge & \limsup_{n\to\infty}\frac{ S_n}{a_n}\ge \limsup_{k\to \infty}\max_{n_{2k}\le m\le m_{2k}}\frac{  S_m}{a_m} \\
\ge & \limsup_{k\to \infty}\max_{n_{2k}\le m\le m_{2k}}\frac{  S_m-  S_{m_{2k-1}}}{a_m}-\limsup_{k\to \infty} \frac{   |S_{m_{2k-1}}|}{a_{m_{2k-1}}}\frac{a_{m_{2k-1}}}{a_{n_{2k}}}\ge 1.
\end{align*}
Similarly, $ \limsup_{n\to\infty}\frac{- S_n}{a_n}=1$.
So, \eqref{eqthLILnew1.2} holds by noting \eqref{eqthLILnew1.1} and the proof is completed. $\Box$

\bigskip
{\bf Proof of Theorem \ref{thLILnew2}}. For part (a), as in the proof of  Theorems \ref{thLILiid1} and \ref{thLILiid2}, without loss of generality, we can assume  $0<\overline{\sigma}\le \infty$, and denote  $d_n=\sqrt{2n\log\log n}$, $X_n=(-d_n)\vee Y_n\wedge d_n$, $S_n=\sum_{i=1}^n X_i$, $s_n^2=\sum_{i=1}^n \Sbep[X_i^2]$,  $t_n=\sqrt{2\log\log s_n}$ and $a_n=t_ns_n$. By \eqref{eqproofLILiid2.5} and \eqref{eqproofLILiid2.6}, $\{X_n;n\ge 1\}$ satisfies the conditions in Theorem \ref{thLILnonid1}. So, by Theorem \ref{thLILnew1}, \eqref{eqthLILnew1.1} and \eqref{eqthLILnew1.2} hold, which, together with \eqref{eqproofLILiid2.8} and \eqref{eqlimvar}, imply
\eqref{eqthLILnew2.2} and
$$
 V\left(  \liminf_{n\to \infty} \frac{\sum_{i=1}^n  Y_i}{\sqrt{2n\log\log n}}-\overline{\sigma}\; \text{ and }\; \limsup_{n\to \infty} \frac{\sum_{i=1}^n  Y_i}{\sqrt{2n\log\log n}}=\overline{\sigma}\right)=1.
$$
\eqref{eqthLILnew2.3} follows from the above equality through standard arguments.

 Next, we consider the part (b).   Suppose
$$C_{\Capc}\left[\frac{Y_1^2}{\log\log|Y_1|}\right]=\infty, $$
which, as shown in the proof of Theorems \ref{thLILiid1} and \ref{thLILiid2}, implies
 that there exists a sequence $M_n\nearrow \infty$  for which
\begin{align*}
   \sum_{n=1}^{\infty} \Capc\big(|Y_n|\wedge (3M_nd_n)\ge 2M_nd_n \big)=\sum_{n=1}^{\infty} \Capc\big(|Y_n|\ge 2M_nd_n \big) =\infty.
\end{align*}
Let $X_n=|Y_n|\wedge (3M_nd_n)$. Then $X_n$ is tight because it is bounded. By Lemma \ref{lemBC} (ii), it follows that
$$ \lim_{N\to \infty} \Capc\Big(\max_{n\le m\le N } \frac{|X_m|}{M_md_m}\ge 1\Big)=1\; \text{ for all } n\ge 1. $$
Let $\epsilon_k=2^{-k}$. Choose $n_k\nearrow\infty$ such that
$$ \Capc\Big(\max_{n_k+1\le m\le n_{k+1} } \frac{|X_m|}{M_md_m}\ge 1\Big)\ge 1-\epsilon_k. $$
Notice the independence, the tightness of each  $X_n$, and that $\bigcap_{k=\ell}^N \big\{\max\limits_{n_k+1\le m\le n_{k+1} } \frac{|x_m|}{M_md_m}\ge 1-\epsilon_k\big\}$ is a closed set.  With the same arguments as in \eqref{eqproofthLILnew1.10}, we can show that
\begin{align*}
  & \Capc^{\mathscr{P}}\left( A_\ell\right)
\ge  \prod_{k=\ell}^{\infty} \Capc\left(  \max_{n_k+1\le m\le n_{k+1} } \frac{|X_m|}{M_md_m}\ge 1 \right)\to 1,\\
&\text{ with } A_\ell=\bigcap_{k=\ell}^{\infty} \big\{\max_{n_k+1\le m\le n_{k+1} } \frac{|X_m|}{M_md_m}\ge 1-\epsilon_k\big\}.
\end{align*}
Hence $ \Capc^{\mathscr{P}}\left(\bigcup_{\ell=1}^{\infty} A_\ell\right)=1$. On the event  $\bigcup_{\ell=1}^{\infty} A_\ell$, we have $\limsup_{n\to\infty}\frac{|Y_n|}{d_n}=\infty$ and hence $\limsup_{n\to \infty} \frac{|\sum_{i=1}^nY_i|}{d_n}=\infty$,
which contradicts \eqref{eqthLILiid2.6}. It follows that \eqref{eqLILmomentcondition1} holds.

Under \eqref{eqLILmomentcondition1}, $\breve{\mathbb E}[Y_1]$ and  $\breve{\mathbb E}[-Y_1]$ exist and are finite.     On the other hand, by Lemma \ref{lem3},
$$ \Capc\left(\frac{\sum_{i=m+1}^{m+n} Y_i-n \breve{\mathbb E}[Y_1]}{d_n}\ge -\epsilon \right)=1-\cCapc\left(\frac{\sum_{i=1}^{n} (-Y_i+ \breve{\mathbb E}[Y_1])}{d_n}>  \epsilon \right)\to 1 \;\;\text{ for all } \epsilon>0. $$
Let $\epsilon_k=2^{-k}$. We can choose a sequence $n_k\nearrow \infty$ with $n_{k-1}/n_k\to 0$ such that
$$ \Capc\left(\frac{\sum_{i=n_{k-1}+1}^{n_k} (Y_i- \breve{\mathbb E}[Y_1])}{d_{n_k-n_{k-1}}}\ge -4\epsilon_k \right) \ge 1-\epsilon_k. $$
Without loss of generality, we can assume $d_{n_k-n_{k-1}}/d_{n_k}\ge 3/4$. Then
$$ \Capc\left(\frac{\sum_{i=n_{k-1}+1}^{n_k} (Y_i- \breve{\mathbb E}[Y_1])}{d_{n_k}}\ge -3\epsilon_k \right) \ge 1-\epsilon_k. $$
Notice the independence. On the other hand, \eqref{eqLILmomentcondition1} implies that $Y_n$ is tight.  Again, with the same arguments as in \eqref{eqproofthLILnew1.10}, we have
$$ \Capc^{\mathscr{P}}\left(\bigcap_{k=\ell}^{\infty}A_k\right)\ge \prod_{k=\ell}^{\infty}\Capc\left(\frac{\sum_{i=n_{k-1}+1}^{n_k} (Y_i- \breve{\mathbb E}[Y_1])}{d_{n_k}}\ge -3\epsilon_k \right) \ge \prod_{k=\ell}^{\infty}( 1-\epsilon_k)\to 1, $$
where $A_k=\big\{\frac{\sum_{i=n_{k-1}+1}^{n_k} (Y_i- \breve{\mathbb E}[Y_1])}{d_{n_k}}\ge -2\epsilon_k\big\}$. It follows that
$\Capc^{\mathscr{P}}\left(\bigcup_{\ell=1}^{\infty}\bigcap_{k=\ell}^{\infty}A_k\right)=1$. On the event $\bigcup_{\ell=1}^{\infty}\bigcap_{k=\ell}^{\infty}A_k$ we have
$$\liminf_{k\to \infty}\frac{\sum_{i=n_{k-1}+1}^{n_k} (Y_i- \breve{\mathbb E}[Y_1])}{d_{n_k}}\ge 0. $$
On the other hand, on the event $\big\{\limsup_{n\to \infty}\frac{|\sum_{i=1}^nY_i|}{d_n}<\infty\big\}$, we have
$$ \limsup_{k\to \infty} \frac{\sum_{i=n_{k-1}+1}^{n_k} Y_i}{d_{n_k}}\le 2\limsup_{n\to \infty}\frac{|\sum_{i=1}^nY_i|}{d_n}<\infty. $$
It follows that
$$  \Capc^{\mathscr{P}}\left(\liminf_{k\to \infty}\frac{(n_k-n_{k-1})\breve{\mathbb E}[Y_1])}{d_{n_k}}<\infty \right)>0. $$
Therefore,
$ \breve{\mathbb E}[Y_1]\le  0. $
Similarly, $\breve{\mathbb E}[-Y_1]\le 0$. From the fact that $\breve{\mathbb E}[Y_1]+\breve{\mathbb E}[-Y_1]\ge 0$, it follows that \eqref{eqLILmeanzerocondition}  holds.

 Under \eqref{eqLILmomentcondition1}  and \eqref{eqLILmeanzerocondition}, we  still have  \eqref{eqthLILiid2.4}  which contradicts \eqref{eqproofLILiid2.14}   if $\overline{\sigma}=\infty$.
 Hence \eqref{eqthLILiid2.8}  holds.
  $\Box$

\bigskip
{\bf Proof of Theorem \ref{thLILnew3} and \ref{thLILnew4}.} Before we prove Theorem \ref{thLILnew3}, we first notice that, the condition (CC) in Proposition \ref{prop2} is satisfied for the new space $(\widetilde{\Omega}, \widetilde{\mathscr{H}},\widetilde{\mathbb E})$ and the family $\widetilde{\mathscr{P}}$ of probability measures. In fact,  $\widetilde{\mathbb V}(\bm x: |x_i|\ge c)\le \Capc(|X_i|\ge c/2)\to 0$ as $c\to \infty$ by \eqref{eqV-V} and the tightness of $X_i$
(which is implied by $\Sbep[X_i^2]<\infty$). The condition (c) in Lemma \ref{lem4.2} is satisfied with $\mathcal T=\{1,2,\ldots\}$.

For Theorem \ref{thLILnew3}, by noting that \eqref{eqV-V}, \eqref{eqthLILnonid2.4} and \eqref{eqthLILnonid2.5ad} hold under $V=\widetilde{\Capc}^{\mathscr{P}}$, $\widetilde{\mathbb C}^{\ast}$ or $\widetilde{\Capc}^{\ast}$ for $\{\tilde X_n;n\ge 1\}$, which imply \eqref{eqthLILnew3.1} and \eqref{eqthLILnew3.2} as shown in Theorem \ref{thLILnew1}.

For Theorem \ref{thLILnew4},  now, $\widetilde{\mathbb E}$ in \eqref{eqcopyE}  is defined with $\{Y_n;n\ge 1\}$ taking the place of $\{X_n;n\ge 1\}$, and
$\tilde Y_n(\widetilde{\omega})=x_n$ for $\widetilde{\omega}=(x_1,x_2,\ldots)$. Also, \eqref{thLILiid2} implies the tightness of $Y_n$. Hence, the condition (CC) in Proposition \ref{prop2} is satisfied for the new space $(\widetilde{\Omega}, \widetilde{\mathscr{H}},\widetilde{\mathbb E})$ and the family $\widetilde{\mathscr{P}}$.  Further, for the new sequence $\{\tilde Y_n;n=1,2,\ldots\}$, it is obvious that
$\lim_{c\to \infty}\widetilde{\mathbb E}[(-c)\vee (\pm \tilde Y_1)\wedge c]=\lim_{c\to \infty}\Sbep[(-c)\vee (\pm Y_1)\wedge c]$,
$\lim_{c\to \infty}\widetilde{\mathbb E}[(-c)\vee (\pm \tilde Y_1^2)\wedge c]=\lim_{c\to \infty}\Sbep[(-c)\vee (\pm Y_1^2)\wedge c]$, and
$$C_{V}\left[\frac{\tilde Y_1^2}{\log\log|\tilde Y_1|}\right]=C_{\Capc}\left[\frac{Y_1^2}{\log\log|Y_1|}\right]
$$
for $V=\widetilde{\Capc}^{\widetilde{\mathscr{P}}}$, $\widetilde{\mathbb C}^{\ast}$ or $\widetilde{\Capc}^{\ast}$, by \eqref{eqV-V3}. Now,
\eqref{eqthLILnew4.2}, \eqref{eqthLILnew4.3} and (b) follow from Theorem \ref{thLILnew2} immediately. It  remains to show \eqref{eqthLILnew4.4} and (c).

When $\sigma>\overline{\sigma}$, \eqref{eqthLILnew4.4} follows from \eqref{eqthLILnew4.2}. When   $\sigma<\underline{\sigma}$, \eqref{eqthLILnew4.4} follows from Theorem \ref{thLIL5}. Next, suppose $\sigma\in [\underline{\sigma},\overline{\sigma}]$.
Let $d_n=\sqrt{2n\log\log n}$, $Z_n=(-d_n)\vee \tilde Y_n\wedge d_n$, $S_n=\sum_{i=1}^n Z_i$.  Then
$V(Z_n\ne \tilde Y_n\;\; i.o.)=0$ due to the fact that
 $$ \sum_{n=1}^{\infty}V(Z_n\ne \tilde Y_n)\le \sum_{n=1}^{\infty}\widetilde{\Capc}(|\tilde Y_n|\ge a_n)<\infty, $$
 by  Lemma \ref{lem2}.
 So,  for \eqref{eqthLILnew4.4} it is sufficient to show that for any $\sigma\in [\underline{\sigma}, \overline{\sigma}]$,
\begin{equation}\label{eqproofnewcor5.3.1}
\widetilde{\Capc}^{\widetilde{\mathscr{P}}}\left(\liminf_{n\to\infty}\frac{S_n}{d_n}=-\sigma \text{ and } \; \limsup_{n\to\infty}\frac{S_n}{d_n}=\sigma\right)=1.
\end{equation}
By the expression \eqref{eqexpressbyP}, for each $i$, there are probability measures  $P^{(1)}, P^{(2)}\in\widetilde{\mathscr{P}}$   such that
 $$ P^{(1)}[Z_i^2]=\widetilde{\mathbb E}[Z_i^2], \;\; P^{(2)}[Z_i^2]=-\widetilde{\mathbb E}[-Z_i^2]. $$
 We consider a mixture of $P^{(1)}$ and $P^{(2)}$ as
 $$ P_i=\alpha_i  P^{(1)}+(1-\alpha_i)P^{(2)}\in\widetilde{\mathscr{P}}  \text{ such that } P_i[Z_i^2]\to \sigma^2. $$
The coefficient $\alpha_i$ is chosen as follows.  When $\sigma=\underline{\sigma}$ (finite or infinite), we choose $\alpha_i\equiv \alpha=0$. When $\sigma=\overline{\sigma}$ (finite or infinite), we choose $\alpha_i\equiv \alpha=1$.  When  $\sigma\in (\underline{\sigma}, \overline{\sigma})$ and $\overline{\sigma}$ is finite, there exists   $\alpha_i\equiv\alpha\in(0,1)$ such that $\sigma^2=\alpha_i\overline{\sigma}^2+(1-\alpha_i)\underline{\sigma}^2$.
 If $\underline{\sigma}<\sigma<\overline{\sigma}=\infty$, then
$$ -\widetilde{\mathbb E}[-Z_i^2]\to \underline{\sigma}^2 \text{ and }  \widetilde{\mathbb E}[Z_i^2]\to  \infty \; \text{ as } i\to \infty, $$
and  there exists $\alpha_i\in (0,1)$ such that
$$  \alpha_i\widetilde{\mathbb E}[Z_i^2]-(1-\alpha_i)\widetilde{\mathbb E}[-Z_i^2]\to \sigma^2. $$
 At any case, we have can choose $\alpha_i\in [0,1]$ such that
 $$P_i[Z_i^2]=\alpha_i\widetilde{\mathbb E}[Z_i^2]-(1-\alpha_i)\widetilde{\mathbb E}[-Z_i^2]\to \sigma^2 \; \text{ as } i\to\infty. $$

For each probability $P_i$, there exists a probability $Q_i$ on $\mathbb R$ such that $Q_i(A)=P_i(\bm y:y_i\in A)$, $A\in \mathscr{B}(\mathbb R)$. Then
\begin{equation}\label{eqQ<V1} Q_i[\varphi]=P_i[\varphi(\tilde Y_i)]\le \widetilde{\mathbb E}[\varphi(\tilde Y_i)]\; \text{ for all } \varphi\in C_{b,Lip}(\mathbb R).
\end{equation}
We define a probability measure $Q$ on $(\mathbb R^{\infty},\mathscr{B}(\mathbb R^{\infty}))$ to be a product probability measure:
$$ Q=Q_1 \times Q_2\times \cdots $$
in sense that
$$
  Q\left(\left\{\bm x:x_i\in C_i, i=1,\ldots, d\right\}\right)  =  \prod_{i=1}^d P_i(\tilde Y_i\in C_i), \;\; C_i \in \mathscr{B}(\mathbb R), d\ge 1.
$$
Such a probability exists and is unique by Kolmogorov's existence theorem.  Notice $\tilde Y_i(\widetilde{\omega})=x_i$ for $\widetilde{\omega}=(x_1,x_2,\ldots)$. We conclude that $\{\tilde Y_i;i=1,2,\ldots\}$ is a sequence of independent random variables under both  $Q$ and $\widetilde{\mathbb E}$.
For $\varphi\in C_{b,Lip}(\mathbb R^d)$,
denote
$$\varphi_1( y_1,\ldots, y_{d-1})=Q[\varphi(y_1,\ldots, y_{d-1}, \tilde Y_d],$$
$$\varphi_2(y_1,\ldots, y_{d-1})=\widetilde{\mathbb E}[\varphi(y_1,\ldots, y_{d-1}, \tilde Y_d]. $$
 Then
$$
\varphi_1(y_1,\ldots, y_{d-1})=P_i[\varphi(y_1,\ldots, y_{d-1}, \tilde Y_d]\le \varphi_2(y_1,\ldots, y_{d-1}),
$$
by noting $P_i\in \widetilde{\mathscr{P}}$ and \eqref{eqQ<V1}. By the independence under both $Q$ and $\widetilde{\mathbb E}$ and noting $P_{i-1}\in \widetilde{\mathscr{P}}$ and \eqref{eqQ<V1} again, we have
\begin{align*}
&Q[\varphi(y_1,\ldots, y_{d-2}, \tilde Y_{d-1},\tilde Y_d)]=Q[\varphi_1(y_1,\ldots, y_{d-2}, \tilde Y_{d-1})] \\
=&P_{d-1}[\varphi_1(y_1,\ldots, y_{d-2}, \tilde Y_{d-1})]\le P_{d-1}[\varphi_2(y_1,\ldots, y_{d-2}, \tilde Y_{d-1})]\\
\le & \widetilde{\mathbb E}[\varphi_2(y_1,\ldots, y_{d-2}, \tilde Y_{d-1})]=\widetilde{\mathbb E}[\varphi(y_1,\ldots, y_{d-2}, \tilde Y_{d-1},\tilde Y_d)].
\end{align*}
 By induction, we conclude  that
$$ Q[\varphi(\tilde Y_1, \ldots,   \tilde Y_d)] \le  \widetilde{\mathbb E}[\varphi(\tilde Y_1, \ldots,   \tilde Y_d)]. $$
It follows that, $Q[\varphi]\le \widetilde{\mathbb E}[\varphi]$, $\varphi\in \mathscr{H}_b$. Hence, $Q\in \widetilde{\mathscr{P}}$.  We conclude that $\{\tilde Y_i;i=1,2,\ldots\}$ is a sequence of independent random variables under $Q$ with
\begin{equation}\label{eqproofnewcor5.3.4} Q[\varphi(\tilde Y_i)]=P_i[\varphi(\tilde Y_i)]\le \widetilde{\mathbb E}[\varphi(\tilde Y_i)]\; \text{ for all } \varphi\in C_{b,Lip}(\mathbb R),
\end{equation}
\begin{equation}\label{eqproofnewcor5.3.5} Q[\varphi(\tilde Y_1,\ldots, \tilde Y_d)]=Q[\varphi\circ \pi_d]\le \widetilde{\mathbb E}[\varphi\circ \pi_d]= \widetilde{\mathbb E}[\varphi(\tilde Y_1,\ldots, \tilde Y_d))], \;\;  \varphi\in C_{b,lip}(\mathbb R^d),\;  d\ge 1.
\end{equation}
and
\begin{equation}\label{eqproofnewcor5.3.6}   Q(B)\le \sup_{P\in\widetilde{\mathscr{P}}}P(B)=\widetilde{\Capc}^{\widetilde{\mathscr{P}}}(B)\; \text{ for all } \; B\in \mathscr{B}(\mathbb R^{\infty}).
\end{equation}

Now, we show \eqref{eqproofnewcor5.3.1}. When $\sigma=0$, then
by \eqref{eqproofnewcor5.3.4}, $Q[Z_i^2]=\alpha_i\widetilde{\mathbb E}[Z_i^2]-(1-\alpha_i)\widetilde{\mathbb E}[-Z_i^2]\le \vSbep[Y_1^2]=0$. Hence $Q(Z_1=Z_2=\cdots=0)=1$ and so
$$ Q\left(\lim_{n\to \infty}\frac{S_n}{d_n}=0\right)=1. $$
It follows that \eqref{eqproofnewcor5.3.1}  holds by \eqref{eqproofnewcor5.3.6}.

Now, suppose $\sigma>0$, then
$Q[Z_i^2]=P_i[Z_i^2] \to \sigma^2>0. $
By \eqref{eqproofnewcor5.3.5},
$$ Q[Z_i]\le \widetilde{\mathbb E}[Z_i]=\Sbep[(-d_i)\vee Y_i\wedge d_i]\le \vSbep[Y_1]+ \vSbep[(|Y_1|-d_i)^+]\le \vSbep[(|Y_1|-d_i)^+]. $$
Similarly,
$$ Q[-Z_i]\le \widetilde{\mathbb E}[-Z_i]\le  \vSbep[-Y_1]+ \vSbep[(|Y_1|-d_i)^+]\le \vSbep[(|Y_1|-d_i)^+]. $$
Let $B_n=\sum_{i=1}^nQ[Z_i^2]$ and $\widetilde a_n=\sqrt{2B_n\log\log B_n}$. Then $\widetilde a_n\ge c_0d_n$.
By Lemma \ref{lem2} (iii), it follows that
$$\sum_{i=1}^n|Q[Z_i]|\le \sum_{i=1}^n\vSbep[(|Z_1|-d_i)^+]=\sum_{i=1}^no\left(\sqrt{\log\log i}/\sqrt{i}\right)
=o(d_n)=o(\widetilde a_n). $$
By \eqref{eqproofnewcor5.3.5}  again and \eqref{eqproofLILiid2.6},
\begin{align*}
\sum_{n=1}^{\infty}\frac{Q[|Z_n|^p]}{\widetilde a_n^p}\le c\sum_{n=1}^{\infty}\frac{Q[|Z_n|^p]}{d_n^p}\le c\sum_{n=1}^{\infty}\frac{\widetilde{\mathbb E}[|\tilde Y_n|^p\wedge d_n^p]}{d_n^p}=
c\sum_{n=1}^{\infty}\frac{\Sbep[|Y_n|^p\wedge d_n^p]}{d_n^p}<\infty.
\end{align*}
 By Theorem \ref{thLILnew1} for a probability $Q$,
$$ Q\left(\liminf_{n\to\infty}\frac{S_n}{\widetilde a_n}=-1 \text{ and } \; \limsup_{n\to\infty}\frac{S_n}{\widetilde a_n}=1\right)=1. $$
Notice $B_n/n\to \sigma^2$. We conclude that
$$ Q\left(\liminf_{n\to\infty}\frac{S_n}{d_n }=-\sigma \text{ and } \; \limsup_{n\to\infty}\frac{S_n}{  d_n }=\sigma\right)=1. $$
By \eqref{eqproofnewcor5.3.6}, \eqref{eqproofnewcor5.3.1}  is proved.

Now, suppose that \eqref{eqthLILnew4.6} holds for a constant $b$. Then
$$ \widetilde V\left(   \limsup_{n\to \infty}\frac{|\sum_{i=1}^n \tilde Y_i|}{\sqrt{2n\log\log n}}=+\infty \right)\le
\widetilde V\left(   \limsup_{n\to \infty}\frac{|\sum_{i=1}^n \tilde Y_i|}{\sqrt{2n\log\log n}}\ne b  \right)<1. $$
By (b), \eqref{eqLILmomentcondition1}, \eqref{eqLILmeanzerocondition}   and \eqref{eqthLILiid2.8} hold. Then by \eqref{eqthLILnew4.4}, for any $\sigma\in [\underline{\sigma},\overline{\sigma}]$,
$$ \widetilde V\left(   \limsup_{n\to \infty}\frac{|\sum_{i=1}^n \tilde Y_i|}{\sqrt{2n\log\log n}}=\sigma \right)=1. $$
It follows that
$$ \widetilde V\left(   \limsup_{n\to \infty}\frac{|\sum_{i=1}^n \tilde Y_i|}{\sqrt{2n\log\log n}}=\sigma \text{ and } =b\right)
\ge 1-\widetilde V\left(   \limsup_{n\to \infty}\frac{|\sum_{i=1}^n \tilde Y_i|}{\sqrt{2n\log\log n}}\ne b \right)>0. $$
Hence, $b=\sigma$ for all $\sigma\in [\underline{\sigma},\overline{\sigma}]$. We must have $\underline{\sigma}=\overline{\sigma} =b$.
  The proof is completed. $\Box$

\bigskip
 {\bf Acknowledgements}  Thanks to Professor Mingshang Hu for the constructive discussion which improved our original manuscript and the revision.   Special thanks to the anonymous referees for carefully reading the manuscript and constructive comments. An example given by the referees  led us to consider carefully the relationship between the capacity $\outCapc$ and the probability measure, and the properties of $\outCapc$.


\end{document}